\documentclass[11pt,a4paper]{amsart}

\usepackage{fullpage,amsxtra,amssymb,amsthm,amsmath}
\usepackage{supertabular}

\renewcommand{\leq}{\leqslant}
\renewcommand{\geq}{\geqslant}

\numberwithin{equation}{section}
\newcommand{\mysection}[1]{\section{{#1}}}

\def\stacksum#1#2{{\stackrel{{\scriptstyle #1}}
{{\scriptstyle #2}}}}

\newcommand{\Cc}{\mathbf{C}}

\newcommand{\Zz}{\mathbf{Z}}

\newcommand{\Rr}{\mathbf{R}}
\newcommand{\Gg}{\mathbf{G}}

\newcommand{\Qq}{\mathbf{Q}}
\newcommand{\Fp}{\mathbf{F}}
\newcommand{\mmu}{\boldsymbol{\mu}}
\newcommand{\Oc}{\mathcal{O}}
\newcommand{\Fbar}{\overline{\mathbf{F}}}

\newcommand{\mods}[1]{\,(\mathrm{mod}\,{#1})}

\newcommand{\ideal}[1]{\mathfrak{{#1}}}

\newcommand{\ra}{\rightarrow}
\newcommand{\lra}{\longrightarrow}
\newcommand{\surjecte}{\twoheadrightarrow}
\newcommand{\injecte}{\hookrightarrow}

\DeclareMathOperator{\spec}{Spec}

\DeclareMathOperator{\res}{Res}

\DeclareMathOperator{\Reel}{Re}

\DeclareMathOperator{\li}{li}
\DeclareMathOperator{\frob}{\sigma}
\DeclareMathOperator{\Gal}{Gal}

\DeclareMathOperator{\im}{Im}
\DeclareMathOperator{\Tr}{Tr}

\DeclareMathOperator{\End}{End}
\DeclareMathOperator{\Aut}{Aut}

\DeclareMathOperator{\disc}{disc}
\DeclareMathOperator{\Res}{Res}

\newcommand{\eps}{\varepsilon}

\newcommand{\dun}{d_1}
\newcommand{\dde}{d_2}
\newcommand{\duns}{D_1}

\newcommand{\oun}{\underline{O}}
\newcommand{\ono}{O}
\newcommand{\rhol}{\hat{\rho}_{\ell}}
\newcommand{\bark}{\bar{K}}

\newcommand{\dalpha}[1]{d_{\alpha}({#1})}
\newcommand{\ddalpha}[2]{d_{{#1}}({#2})}
\newcommand{\gale}{g}
\newcommand{\outside}{outside }
\newcommand{\inductive}{inductive }

\newcommand{\somg}{S}
\newcommand{\somp}{T}
\newcommand{\ejum}{j}
\newcommand{\ejumbis}{J}

\newcommand{\prim}{\mathcal{U}}
\newcommand{\nrmone}{T}
\newcommand{\bij}{\eta}
\newcommand{\grandm}{M}
\newcommand{\petitm}{m}
\newcommand{\frobs}{\Sigma}

\newcommand{\pmax}{60,000,000}

\DeclareMathSymbol{\gena}{\mathord}{letters}{"3C}
\DeclareMathSymbol{\genb}{\mathord}{letters}{"3E}

\def\sumb{\mathop{\sum \Bigl.^{\flat}}\limits}

\def\sums{\mathop{\sum \Bigl.^{*}}\limits}

\theoremstyle{plain}
\newtheorem{theorem}{Theorem}[section]

\newtheorem{lemma}[theorem]{Lemma}
\newtheorem{corollary}[theorem]{Corollary}

\newtheorem{problem}[theorem]{Problem}
\newtheorem{proposition}[theorem]{Proposition}

\theoremstyle{remark}
\newtheorem{remark}[theorem]{Remark}

\theoremstyle{definition}

\newtheorem*{definition}{Definition}
\newtheorem*{question}{Question}
\newtheorem{example}[theorem]{Example}

\begin{document}

\title{Analytic problems for elliptic curves}

\author{E. Kowalski}
\address{Universit\'e Bordeaux I - A2X\\
351, cours de la Lib\'eration\\
33405 Talence Cedex\\
France}
\email{emmanuel.kowalski@math.u-bordeaux1.fr}
\subjclass{Primary 11N99; Secondary 11G05, 11G20, 11F99}
\keywords{Elliptic curves, sieves, trace formula for Hecke operators,
  Chebotarev density theorem} 
\begin{abstract}
We consider some problems of analytic number theory for elliptic
curves which can be considered as analogues of classical questions
around the distribution of primes in arithmetic progressions to large
moduli, and to the question of twin primes. This leads to some local
results on the distribution of the group structures of elliptic curves
defined over a prime finite field, exhibiting an interesting
dichotomy for the occurence of the possible groups. (This paper was
initially written in 2000/01, but after a four year wait for a referee
report, it is now withdrawn and deposited in the arXiv).
\end{abstract}
\maketitle

\tableofcontents

\mysection{Introduction}
\label{sec-intro}

This paper introduces and discusses some problems of analytic number
theory which are related to the arithmetic of elliptic curves over
number fields. One can see them as analogues of some very classical
problems about 
the distribution of prime numbers, especially primes in arithmetic
progressions to large moduli. The motivation comes both from
these analogies and from the conjecture of Birch and Swinnerton-Dyer. 
\par
To explain this,
consider an elliptic curve $E$ defined over $\Qq$, given 
by a (minimal) Weierstrass equation (\cite[VII-1]{silv})
\begin{equation}
\label{eq-weirstrass-min}
y^2+a_1xy+a_3y=x^3+a_2x^2+a_4x+a_6
\end{equation}
with $a_i\in \Zz$. For all primes $p$ we can consider the reduced
curve $E_p$ modulo $p$, which for almost all $p$ will be an elliptic curve
over the finite field $\Fp_p=\Zz/p\Zz$.
We wish to study the behavior of sums of the type
\begin{equation}\label{eq-somme}
\sum_{p\leq X}{\iota_E(p)}
\end{equation}
as $X\ra +\infty$, where $\iota_E(p)$ is some invariant attached to the
reduced curve $E_p$ and to its finite group of $\Fp_p$-rational points
in particular. For example, taking
$$
\iota_E(p)=\frac{|E_p(\Fp_p)|}{p}
$$
one gets the sum
$$
\sum_{p\leq X}{\frac{|E_p(\Fp_p)|}{p}}
$$
which should be related to the behavior of the logarithmic derivative
of the Hasse-Weil zeta function of $E$ at $s=1$, and so conjecturally
to a global invariant of $E/\Qq$, the rank of its Mordell-Weil group
$E(\Qq)$.
\par
We wish to consider other sums of type~(\ref{eq-somme}) which are
natural from the point of view of analytic number theory. The hope is
to get precise enough asymptotics where global invariants of $E$
would enter, to gain an understanding of the local-global principles
which the Birch and Swinnerton-Dyer conjecture postulates.
\par
\medskip
The plan of this paper is as follows: in the first section
we state basic facts on elliptic curves that we'll use and
introduce some natural invariants $\iota_E(p)$. In
Section~\ref{sec-split}, we show how the study
of the sum~(\ref{eq-somme}) for one of them brings about questions
involving the equidistribution of Frobenius elements (in the extensions
of $\Qq$ generated by the torsion points of $E$) to uniform and
large moduli, especially on totally split primes in such
extensions. We analyze this problem on GRH and discuss the
new difficulties which arise in comparison
with the case of primes in arithmetic progressions. There are several
remarks here which may be of interest. One of the new phenomenon
(primes splitting completely in fields generated by $d$ torsion points
with $d$ very large) leads us to a notion of elliptic twins, analogues
of the classical twin primes that we again discuss in general terms.
At long last, non-trivial results are obtained in the next two
sections: for CM curves, in Section~\ref{sec-cm-resultats}, sieve
techniques in quadratic fields can be usefully applied, and in
Section~\ref{sec-local} the subject of totally split primes is
viewed from a different angle: now, given a prime $p$, and $d\geq 1$,
we ask whether or not there exists \emph{some} curve $E$ with $p$
totally split in $\Qq(E[d])$. This is done in two ways, adapting
results of Deuring, Waterhouse, Schoof, and using the trace formula
and modular curves. Finally, since the problems are amenable to
experimentation, we present in Section~\ref{sec-numerics} some
numerical data and further remarks.
\par
Most of the results presented here are not very strong and the overall
situation remains rather unsatisfactory.
The excuse for 
this is that the problems seem genuinely difficult. On the other hand,
to the author at least, their interest is very obvious. 
\par
\bigskip
\textbf{Notation.} The symbols $\ono()$ and $o()$ are used in the sense
of (for example) Bourbaki, so $f(x)=\ono(g(x))$ as $x\ra x_0$ means that
for $x$ in some neighborhood $U$ of $x_0$ we have $|f(x)|\leq Cg(x)$ for some
$C\geq 0$ (depending on $U$). On the other hand $f \ll g$ is used in the
sense that there exists $C\geq 0$ such that for
all $x$ (in some set to be described explicitly or implicitly) we have
$|f(x)|\leq Cg(x)$. The dependence of $C$ on other parameters is
indicated by subscripts $\ll_{\eps}$, etc.
\par
For notational convenience\footnote{ Many
papers in analytic number theory actually use the notation $O()$ in
this sense, and correspondingly speak of ``hidden constants'', as for
$\ll $.
} it is sometimes useful to use another
symbol $\oun()$ such that $f=\oun(g)$ is equivalent to $f\ll g$.
\par
It will be convenient in a number of places to use the following
notation: for every real number $x$, we let
\begin{equation}
\label{eq-x-pm}
x^-=(\sqrt{x}-1)^2,\quad x^+=(\sqrt{x}+1)^2.
\end{equation}
(defined for $x\geq 1$, $x\geq -1$, respectively).
Note that
$$
(x^+)^-=(x^-)^+=x,\text{ and }
x^+-x^-=4\sqrt{x}.
$$
\par
\bigskip
\textbf{Notice (2005).} Up to some updates of the numerical data and
typographical corrections, this is the
version of this text that was submitted to the Transactions of the
A.M.S on Oct. 10, 2001. After four years of wait, I have withdrawn the
paper to put it on \verb|arXiv| instead. The lengthy delay means that the
bibliography is not quite up to date; in particular, some papers of
A. Cojocaru (including collaborations with W. Duke, R. Murty)
are somewhat related to the topics presented here, see
Math. Ann. 329 (2004), 531--534; Math. Ann.  330  (2004),  601--625;
and Trans. A.M.S 355 (2003), 2651--2662 for instance.

\mysection{Some local invariants for elliptic curves}
\label{sec-pre}

In this section we want to define some of the invariants that
are of interest. First we recall some important facts about elliptic
curves.

\subsection{Elliptic curves}
\label{ssec-rappels}

Let $E/K$ be an elliptic curve defined over a field $K$. We will
mostly use ``old-fashioned'' language, identifying $E$ with its set of
$\bark$-valued points, where $\bark$ is a fixed algebraic closure
of $K$.
\par
The endomorphism ring of $E$ over $K$ is denoted by $\End(E)$, and the
endomorphism ring of $E$ over $\bark$ by $\End_{\bark}(E)$. The ring
$\End(E)$ contains the subring $\Zz$ corresponding to the morphisms $x
\mapsto nx$ for $n\in \Zz$. When $\End(E)$ is strictly bigger than
$\Zz$, the curve is said to be CM, or to have complex multiplication.
\par
To any $\varphi\in \End(E)$ is associated its \emph{dual}
$\overline{\varphi}\in \End(E)$ with the property that $\varphi\circ
\overline{\varphi}=\overline{\varphi}\circ \varphi=[\deg(\varphi)]$,
the multiplication by the degree of $\varphi$, as a morphism of
algebraic curves (\cite[III-6]{silv}).
\par
The various possibilities for $\End(E)$ have been studied extensively by
Deuring~\cite{deuring}. We are concerned with two cases. Let
$\Oc=\End(E)$.
\begin{itemize}
\item{
If $K$ is a finite field, $E$ is always CM, and $\Oc$ is either
an order in an imaginary quadratic field, in which case $E$ is said
to be ordinary, or an order in a quaternion algebra, in which case $E$
is said to 
be supersingular (see~\cite[V-3]{silv}). There are only finitely many
$j$-invariants $j\in \bark$ corresponding to supersingular curves,
all of degree $\leq 2$ over the prime field.
}
\item{
If $K$ is a number field, either $\Oc=\Zz$ or $\Oc$ is an order in an
imaginary quadratic field. In this case $j(E)$ is an algebraic
integer. For fixed $K$, there are only finitely many possible values
of $j\in K$ for which an elliptic curve over $K$ with $j(E)=j$ has
CM (see~\cite[II]{silv-2}, and for instance~\cite[App. A-3]{silv-2}
for a list of all CM curves over $\Qq$). The dual of an endomorphism
$\sigma$ is its (unique) conjugate over $\Qq$.
}
\end{itemize}
\par
\medskip
Let $E/K$ be an elliptic curve defined over a field $K$. For every
integer $d\geq 1$, the $d$-torsion points of $E$ form (depending on the point
of view) either a finite subgroup or a finite subgroup scheme of $E$,
denoted either $E[d]$ or $E[d](\bark)$ depending on the emphasis.
\par
The structure of this group depends on the characteristic $p$ of $K$
and is given as follows (\cite[III-6.4]{silv}):
\begin{itemize}
\item If $d_1$ and $d_2$ are coprime, then
$$
E[d_1d_2](\bark)=E[d_1](\bark)\oplus
E[d_2](\bark).
$$
\item If $(d,p)=1$ (in particular, if $K$ is of characteristic $0$),
  we have
$$
E[d](\bark)\simeq \Zz/d\Zz \oplus \Zz/d\Zz.
$$
\item If $K$ is a finite field, $d=p^v$ with $v\geq 1$ and $E$ is
  ordinary, then 
$$
E[d](\bark)\simeq \Zz/d\Zz.
$$
\item If $K$ is a finite field, $d=p^v$ with $v\geq 1$ and $E$ is
  supersingular, then 
$$
E[d](\bark)=0.
$$
\end{itemize}

In any case, $E[d]$ is a finite (and free) $\Zz/d\Zz$-module, and the
natural action of the Galois group $G_K=\Gal(\bark/K)$ induces
a Galois representation
$$
\rho_d(E)\,:\,G_K\lra \Aut(E[d]).
$$
Assume now that $(d,p)=1$, then by choosing a basis we get
$2$-dimensional representations, well-defined up to conjugacy
$$
\rho_d(E)\,:\, G_K\lra GL(2,\Zz/d\Zz).
$$
Those are compatible, meaning that if $e\mid d$, then we have
$$
\rho_e(E)=\rho_d(E)\mods{e}
$$
with obvious notations. In particular, taking a prime $\ell\not=p$ and
$d=\ell^v$ for all $v\geq 0$, we obtain a projective system of
representations into $GL(2,\Zz/\ell^v\Zz)$ which can be put together
into an integral $\ell$-adic representation
$$
\rhol\,:\,G_K\lra GL(2,\Zz_{\ell}).
$$
\par
Let now $K=\Fp_q$ be a finite field with $q$ elements, of
characteristic $p$ (this will be a standing convention). The group 
of $\Fp_q$-rational points on $E$ is finite. We write
$$
n(E)=|E(\Fp_q)|
$$
for its order. The most important invariant of $E/\Fp_q$ is the
integer $a(E)$ such that
\begin{equation}
n(E)=|E(\Fp_q)|=q+1-a(E).
\end{equation}
One knows that $a(E)$ characterizes the isogeny class of $E$ over
$K$ (see~\cite[Ex. 5.4]{silv}). Moreover, $E$ is supersingular if and
only if $p\mid a(E)$. In case $q=p$, this is equivalent with $a(E)=0$
(see~(\ref{eq-rh})), so 
there is a unique isogeny class of supersingular curves defined over
the base field $\Fp_p$.
\par
The integer $a(E)$ is also linked to $\End(E)$. The Frobenius automorphism
$\frob\,:\,x\mapsto x^q$ of $\Fbar_q$ is an
element of $\End(E)$. We have (\cite[V]{silv})
\begin{align}
a(E)&=\Tr(\sigma)=\sigma+\overline{\sigma}\\
n(E)&=N(\sigma-1)=(\sigma-1)(\overline{\sigma}-1).
\end{align}
\par
For any integer $d$ with $(d,p)=1$, $a(E)$ is further related to the
Galois representation $\rho_d(E)$ (\cite[V]{silv}) by
\begin{equation}
\label{eq-galois-rep1}
\det(\rho_d(\frob))=q\mods{d}\quad\text{ and }\quad
\Tr(\rho_d(\frob))=a(E)\mods{d}.
\end{equation}
\par
Hence the $\ell$-adic representation $\rhol$ satisfies the
fundamental property
\begin{equation}
\label{eq-galois-rep}
\det(\rhol(\frob))=q\quad\text{ and }\quad
\Tr(\rhol(\frob))=a(E).
\end{equation}

Hasse proved (the Riemann Hypothesis for elliptic curves over finite
fields, see~\cite[V-1.1]{silv}) that 
\begin{equation}
\label{eq-rh}
|a(E)|\leq 2\sqrt{q}.
\end{equation}

If $K$ is a number field, then for any prime ideal $\ideal{p}$ of $K$
where $E$ has good reduction, the above theory applies to the reduced
curve $E_{\ideal{p}}$ modulo $\ideal{p}$. For instance, the Galois
representation 
$\rhol(E)$ (for any $\ell$ not dividing $\ideal{p}$) satisfies
\begin{equation}
\label{eq-galois-rep-nf}
\det(\rhol(\frob_{\ideal{p}}))=N\ideal{p}\quad\text{ and }\quad
\Tr(\rhol(\frob_{\ideal{p}}))=a_{\ideal{p}}(E)
\end{equation}
where $\frob_{\ideal{p}}$ is a Frobenius element at $\ideal{p}$ and
$a_{\ideal{p}}=a(E_{\ideal{p}})$.
\par
\medskip
For an elliptic curve $E/K$, and an integer $d\geq 1$, we let
$K(E[d])$ denote the finite extension of $K$ obtained by adjoining the
coordinates of the $d$-torsion points of $E$, or in other words the
smallest extension $L/K$ such that $E[d](\bark)\subset E(L)$. 
This is a Galois extension and
in fact $K(E[d])$ is the extension of $K$ corresponding to the closed
subgroup $\ker(\rho_d)$ of $G_K$, i.e. $K(E[d])=\bark^{\ker\rho_d}$,
so that there is a canonical isomorphism
\begin{equation}\label{eq-gal-torsion}
\Gal(K(E[d])/K)\simeq \im \rho_d \subset \Aut(E[d]).
\end{equation}
We will denote $G_d=\Gal(K(E[d])/K)$ when $E$ and $K$ are clear in the
context.
\par
In the case $d=2$, and $K$ of characteristic $\not=2$, if $E/K$ is
given by an equation 
$$
y^2=f(x)
$$
for some cubic polynomial $f\in K[X]$, the $2$-division points of $E$
are the origin, and the points $(\alpha,0)$ where $\alpha$ runs over the
three distinct roots of $E$ in $\bark$. In particular, $E[2]\subset
E(K)$ if and only if $f$ splits into linear factors in $K[X]$.
\par
We let $\mmu_d$ denote the group (scheme) of the $d$-th roots of
unity. It is known (\cite[III-8.11]{silv}) that $K(\mmu_d)\subset
K(E[d])$ where 
$K(\mmu_d)$ is the field obtained by adjoining all $d$-th roots of
unity to $K$. In the case of number fields (resp. finite fields), this
can be seen from~(\ref{eq-galois-rep-nf})
(resp.~(\ref{eq-galois-rep1})): the determinant condition 
implies that primes totally split in $K(E[d])$ are totally split in
$K(\mmu_d)$, which implies that $K(E[d])$ contains $K(\mmu_d)$ (see
e.g. \cite[V-6.8]{neu}).
\par
If $K$ is a number field and $\ideal{p}$ is a prime ideal in $K$ where
$E$ has good reduction, 
the residue field extension of $K(E[d])$ at $\ideal{p}$ is 
isomorphic to $\Fp_{\ideal{p}}(E_{\ideal{p}}[d])$. Indeed the reduction
map $E[d](\bar{K}_{\ideal{p}})\ra
E_{\ideal{p}}(\bar{\Fp}_{\ideal{p}})$ is surjective
(see e.g.~\cite[VII-3.1]{silv} if $(d,\ideal{p})=1$, which will be the
case we need, and adapt~\cite[Ex. IV-4.4]{silv} for the general case).
\par
In this case of a number field, the Galois groups $G_d$ are known ``up
to finite index''.

\begin{theorem}\label{th-galois-nf}
Let $K$ be a number field, $E/K$ an elliptic curve. Then
\par
\emph{1. (Deuring, see~\cite[4.5]{serre-div})} If $E$ has complex
multiplication and $\Oc=\End_{\bar{K}}(E)$, with $\Oc\subset K$, then
$\rho_d$ induces a group homomorphism
$$
\rho_d\,:\, G_K\ra (\Oc/d\Oc)^{\times}
$$
with the property that as $d$ ranges over all integers $d\geq 1$, the
index of $G_d$ in the finite group $(\Oc/d\Oc)^{\times}$ is bounded by
a constant $i(E)$.\footnote{
If $K$ does not contain the endomorphism ring, $G_d$ is at most an
extension of $(\Oc/d\Oc)^{\times}$ by $\Zz/2\Zz$.
}
\par
\emph{2. (Serre~\cite{serre-div})} If $E$ does not have complex
multiplication, then the index of $G_d$ in the finite group
$\Aut(E[d])\simeq GL(2,\Zz/d\Zz)$ is bounded by a
constant $i(E)$.
\end{theorem}

Note that
\begin{align}
|GL(2,\Zz/d\Zz)| &= d\psi(d)\varphi(d)^2
\label{eq-ord-ncm}
\end{align}
where $\varphi$ is Euler's function
and
\begin{equation}\label{eq-psi}
\psi(d)=d\prod_{p\mid d}{\Bigl(1+\frac{1}{p}\Bigr)}.
\end{equation}
Since $\Oc$ is not a Dedekind ring in general, hence does not have
unique factorization into ideals, the order of $(\Oc/d\Oc)^{\times}$
is not a multiplicative function of $d$. If $\Oc$ is the full ring of
integers of its fraction field $k$, or if $d$ is coprime with the
discriminant of $\Oc$, then $|(\Oc/d\Oc)^{\times}|$ is the analogue of
the Euler function for ideals in $k$:
\begin{equation}
  \label{eq-ord-cm}
|(\Oc/d\Oc)^{\times}|
=d^2\prod_{\ideal{p}\mid (d)}
{\Bigl(1-\frac{1}{N\ideal{p}}\Bigr)}.
\end{equation}
\par
Informally, we say that in the CM case, $|G_d|$ is of order of magnitude
$d^2$, and in the non-CM case, $|G_d|$ is of order of magnitude
$d^4$. This difference will be important later on so we define
the \emph{Galois dimension} $\gale=\gale(E)$ of $E$ to be $2$ if $E$
has CM and $4$ if not (it is the dimension of the $\ell$-adic Lie
group $\im(\rhol(G_K))\subset GL(2,\Zz_{\ell})$, or of its Lie
algebra for $\ell$ large enough~\cite{serre-mg}).

\subsection{Local invariants}

First we describe the group structure of the rational
points of an elliptic curve defined over a finite field. This is well
known.

\begin{lemma}
Let $E/\Fp_q$ be an elliptic curve defined over a finite field with
$q$ elements. There exist unique integers $\dun$ and $\dde$ such that
\begin{equation}
\label{eq-grp-struct}
E(\Fp_q)\simeq \Zz/\dun \Zz\oplus \Zz/\dun\dde\Zz.
\end{equation}
\end{lemma}

\begin{proof}
The group $E(\Fp_q)$ is finite, hence of finite exponent, so for some
$d\geq 1$ we have
$$
E(\Fp_q)\subset E[d](\Fbar_q).
$$
As we recalled in Section~\ref{ssec-rappels}, the group on the right
has a system of generators with at most two elements. By the structure
theorem of 
finite abelian groups, the same is true for any subgroup, and they are
all of the form stated.
\end{proof}

The integers $\dun$, $\dde$ are very interesting invariants of $E$. We
will denote them by $\dun(E)$ (resp. $\dde(E)$) or $\dun(\ideal{p})$
(resp. $\dde(\ideal{p})$) when $E$ is obtained by reducing  a curve
over a number field modulo a prime ideal $\ideal{p}$.

\begin{lemma}\label{lm-def-dun}
Let $E/\Fp_q$ be an elliptic curve over a finite field with $q$
elements. Then
\par
\emph{(1)} We have
$$
\dun=\dun(E)=\max\{ d\geq 1\,\mid\, (d,p)=1\text{ and
}E[d](\Fbar_q)\subset E(\Fp_q) \}
$$
i.e $\dun(E)$ is the largest integer $d$ prime to $p$ for which all
of the $d$-torsion is rational over $\Fp_q$. The $\max$ can be taken
with respect to the order by divisibility or the ``linear'' order on
$\Zz$.
\par
\emph{(2)} We have
$$
\dun(E)^2\dde(E)=n(E)=q+1-a(E).
$$
\indent
\emph{(3)} We have
$$
q+1-a(E)=0 \mods{\dun^2},\quad
q=1\mods{\dun},\quad a(E)=2\mods{\dun}.
$$
\indent
\emph{(4)} We have 
\begin{equation}\label{eq-a-priori}
\dun(E)\leq \sqrt{q}+1.
\end{equation}
(see~(\ref{eq-x-pm}) for the definition).
\end{lemma}

\begin{proof}
All this is easy from the structure of the $d$-torsion points. For
(1), observe that the finite abelian group 
$$
\Zz/\dun\Zz\oplus \Zz/\dun\dde\Zz
$$
contains $\dun^2$ points of order $\dun$, namely $\Zz/\dun\Zz\oplus
\dde\Zz/\dun\dde\Zz$. Since it is known a priori
that $E(\Fbar_q)$ contains at most $d^2$ points of order $d$ for any
$d\geq 1$, all the $\dun$-torsion is $\Fp_q$-rational. Moreover,
if there exists $d>\dun $ for which $E[d](\Fbar_q)\subset E(\Fp_q)$,
we can write $d=\dun d'$ for some $d'>1$. Then $d'$ must be of the
form $d'=p^v$ for some $v\geq 1$, since otherwise 
there would be $e^2$ points of order $e$ which are $\Fp_q$-rational,
for some $e>d$, which the group structure~(\ref{eq-grp-struct})
forbids. 
\par
The second point is obvious, and gives the first congruence in~(3),
while~(\ref{eq-galois-rep1}) gives the other congruences.
\par
For~(\ref{eq-a-priori}), since $\dun^2\mid q+1-a(E)$, and
$q+1-a(E)>0$, it follows that $\dun^2\leq q+1-a(E)\leq (\sqrt{q}+1)^2$
by~(\ref{eq-rh}).
\end{proof}

\begin{remark}
The congruence $n(E)=0\mods{\dun^2}$ can also be obtained from the
Galois representations without referring to the points of the
elliptic curve: let $\gamma=\rho(\frob)\in GL(2,\Zz/\dun^2\Zz)$. We
know that $\gamma\equiv 1\mods{\dun}$ by definition. Now writing
$\gamma=1+\dun \gamma'$ and expanding the trace and determinant, we
obtain using~(\ref{eq-galois-rep1}) (both for $d=\dun$ and $d=\dun^2$)
$$
2+d\Tr(\gamma')=a(E)\mods{\dun^2},\quad
1+d\Tr(\gamma')=p\mods{\dun^2}.
$$
Then observe that $\Tr(\gamma')=a(E)\mods{\dun}$ and subtract to get
$1=a-p\mods{\dun^2}$. (This remark is due to N. Katz).
\end{remark}

\begin{remark}\label{rem-base}
If $q=p\geq 3$ the condition $(d,p)=1$ in the
characterization~(1) of $\dun$ can be omitted \emph{unless} either $E$
is supersingular or $a(E)=1$. In the first case, of course,
$0=E[p^n]\subset E(\Fp_p)$ for 
all $n\geq 1$, while in the second case we have $|E(\Fp_p)|=p$ so
$E(\Fp_p)$ is cyclic of order $p$ and must equal $E[p]$. Conversely,
for $E$ ordinary, if $d=p^ne$ with $(e,p)=1$ and $E[d]\subset
E(\Fp_p)$, we get $p^ne^2\mid p+1-a(E)$ and by the Riemann
Hypothesis~(\ref{eq-rh}) it follows immediately that $n=e=1$.
\par
Note that curves with $a(E)\equiv 1\mods{p}$ occur in other contexts.
If $E$ arises by reduction modulo $p$ of a curve over $\Qq$, the prime
$p$ is called \emph{anomalous}~\cite{mazur}. When $p\geq 7$, $a(E)=1$
is the same as $a(E)\equiv 1\mods{p}$, so those curves, and the
supersingular curves, form two isogeny classes of curves over
$\Fp_p$. 
\end{remark}

The next lemma is equally simple.

\begin{lemma}\label{lm-cong-end}
Let $E/\Fp_q$ be an elliptic curve defined over a finite field and
$d\geq 1$ an integer with $(d,p)=1$. Then $E[d]\subset E(\Fp_q)$ if
and only if $\frob \equiv 1\mods{d}$ in $\End(E)$, where $\frob$ is
the Frobenius endomorphism of $E$.
\end{lemma}

\begin{proof}
Let $K=\End(E)\otimes_{\Zz} \Qq$, which is either a quadratic field or
a quaternion algebra over $\Qq$, and let $\frob'=(\frob-1)/d\in
K$. The congruence in the statement of the Lemma means $\frob'\in
\End(E)$; since $d$ is central in $K$, there is no ambiguity in the
side on which $d^{-1}$ is put in the definition of $\frob'$.
\par
Now if $\frob'\in \End(E)$, we have $\frob=1+d\frob'$, so for any
$x\in E[d]$ it follows that $\frob(x)-x=\frob'(dx)=0$, hence $x$ is
$\Fp_q$-rational. Conversely, if $E[d]\subset E(\Fp_q)$, the $\Fp_q$-isogeny
$\phi=\frob-1$ of $E$ satisfies $\ker(d)\subset \ker(\phi)$; since
$(d,p)=1$, multiplication by $d$ is separable, hence
(\cite[III-4.11]{silv}) $\phi$ factorizes by $d\,:\, E\ra E$, which means
$\frob\equiv 1\mods{d}$.
\end{proof}


Here is the global interpretation of $\dun$.
\begin{lemma}\label{lm-dun-glob}
Let $E/K$ be an elliptic curve over a number field, $\ideal{p}$ a
prime ideal such that $E$ has good reduction modulo $\ideal{p}$. For
any integer  $d\geq 1$, we have $d\mid \dun(E_{\ideal{p}})$ if and
only if
$\ideal{p}$ is totally split in the field $K(E[d])$. 
\end{lemma}

\begin{proof}
Both statements imply that $(d,\ideal{p})=1$: this is by definition
for $\dun$ and because if $\ideal{p}$ is totally split, it is
unramified in $K(E[d])$, hence in $K(\mmu_d)$, which implies
$N\ideal{p}=1\mods{d}$.
\par
We know that the residue field extension of
$K(E[d])$ at 
$\ideal{p}$ is $\Fp_{\ideal{p}}(E_{\ideal{p}}[d])$.
If $\ideal{p}$ is totally split, this extension is trivial, so all the
$d$-torsion is rational, i.e. $d\mid \dun(E_{\ideal{p}})$.
\par
Conversely, if $d\mid \dun(E_{\ideal{p}})$, the condition
$(d,\ideal{p})=1$ implies that $\ideal{p}$ is unramified in $K(E[d])$
(\cite[4.1]{silv}). Then the residue field extension being trivial
means that $\ideal{p}$ is totally split.
\end{proof}





\mysection{Totally split primes}
\label{sec-split}

\subsection{The splitting problem for elliptic curves}
\label{ssec-splitting}

Let $E/K$ be an elliptic curve over a number field. Apart from the
number of points $N\ideal{p}+1-a_{\ideal{p}}(E)$ on $E$ modulo a prime
ideal, one of the most natural invariant to insert in a
sum~(\ref{eq-somme}) is $\iota(\ideal{p})=\dun(\ideal{p})$. Thus we
define for $X\geq 1$
\begin{equation}\label{eq-sum-dun}
S_E(X;\dun)=\sum_{N\ideal{p}\leq X}{\dun(\ideal{p})}
\end{equation}
where as before $\dun(\ideal{p})=\dun(E_{\ideal{p}})$ (we define,
rather arbitrarily, $\dun(\ideal{p})=0$ for ramified primes).

\begin{problem} 
\label{pr-splitting}
What is the asymptotic behavior of $S_E(X;\dun)$ as $X\ra
+\infty$ ?
\end{problem}

Because of the following link with primes totally split in division
fields of $E$, we call this the \emph{elliptic splitting problem} for
$E$. 

\begin{lemma}
Let $E/K$ be an elliptic curve over a number field. We have
\begin{equation}
S_E(X;\dun)=\sum_{d\leq \sqrt{X}+1}
{\varphi(d)\pi_E(X;d,1)}
\label{eq-dun-som}
\end{equation}
for $X\geq 1$, where 
\begin{equation}\label{eq-pi-E}
\pi_E(X;d,1)=|\{\ideal{p}\,\mid\, N\ideal{p}\leq X,\text{ and 
$\ideal{p}$ is totally split in $K(E[d])$}\}|.
\end{equation}
\end{lemma}

\begin{proof}
Using the convolution formula
$$
n=\sum_{ab=n}{\varphi(a)}
$$
and~(\ref{eq-a-priori}), we have
\begin{align*}
S_E(X;\dun) &= \sum_{N\ideal{p}\leq X}{\dun(\ideal{p})}\nonumber\\
&= \sum_{N\ideal{p}\leq X}{\sum_{d\mid \dun(\ideal{p})}
{\varphi(d)}}\\
&=\sum_{d\leq \sqrt{X}+1}
{\varphi(d)\sum_{\stacksum{N\ideal{p}\leq X}{d\mid \dun(\ideal{p})}}
{1}}\\
&=\sum_{d\leq \sqrt{X}+1}
{\varphi(d)\pi_E(X;d,1)},\quad\text{by Lemma~\ref{lm-dun-glob}}.
\end{align*}

\end{proof}

\begin{remark}
This lemma shows that the elliptic splitting problem is quite
analogous to the classical \emph{Titchmarsh divisor problem}
(first considered in~\cite{titch-div}) which
concerns the asymptotic behavior of the sum
$$
S(X,d)=\sum_{p\leq X}{d(p-1)}
$$
where $d(n)$ is the number of ($>0$) divisors of $n$. This was solved
by Linnik\footnote{ Titchmarsh had shown
  the result on the Riemann Hypothesis (see also below).}
(see~\cite{linnik}): 
\begin{theorem}
\label{th-linnik}(Linnik)
We have
\begin{equation}\label{eq-tit-asy}
S(X,d)\sim cx\quad\text{with }
c=\prod_p{\Bigl(1+\frac{1}{p(p-1)}\Bigr)}
=\frac{\zeta(2)\zeta(3)}{\zeta(6)}=1.943596\ldots
\end{equation}
as $X\ra +\infty$.
\end{theorem}
Linnik proved this by a very difficult argument using the dispersion
method, although now it is easy to derive from the Bombieri-Vinogradov
theorem and the Brun-Titchmarsh theorem (see e.g.~\cite[\S
3.5]{halb-rich}; we will essentially redo this argument later on).
Although this will not matter here, we mention that Bombieri,
Friedlander, Iwaniec~\cite{bom-fried-iwa} and independently
Fouvry~\cite{fouvry}, have  
proved a more precise formula, with a second term of magnitude
$X/\log X$, using their deep results about primes in arithmetic
progressions to moduli $d>\sqrt{X}$.
\par
The first step in this proof is to write
\begin{align}
d(n)&=\sum_{ab=n}{1}\nonumber\\
&=2\sum_{\stacksum{d\mid n}{d<\sqrt{n}}}{1}+
\begin{cases}
1&\text{if $n$ is a square,}\\
0&\text{otherwise}
\end{cases}
\text{ (Dirichlet's divisor-switching trick)}
\label{eq-trick}
\end{align}
which leads immediately to
\begin{align}
S(X,d)&=\sum_{d\leq X}{\pi(X;d,1)}\nonumber\\
&=2\sum_{d<\sqrt{X}}{\bigl(\pi(X;d,1)
-\pi(d^2+1;d,1)\bigr)}
+\sqrt{X}+\oun(1).
\label{eq-div-som1}
\end{align}
where for any integer $a$, $\pi(X;d,a)$ is the classical counting
function for primes $p\equiv a\mods{d}$. By the elementary
theory of cyclotomic fields, this is also the number of primes $p\leq X$
such that the Frobenius at $p$ acts on $d$-th roots of unity by
$\zeta\mapsto \zeta^a$, so that $\pi(X;d,1)$ is the number of $p\leq
X$ totally split in the cyclotomic field generated by $d$-th roots of
unity. 
\par
Theorem~\ref{th-linnik}, via the formula~(\ref{eq-div-som1}), will
actually be used in Section~\ref{sec-local}, reinforcing the
connection between this classical result and
Problem~\ref{pr-splitting}. We may also remark that another connection
arises if one interprets $d(p-1)$ as the number of subgroups of the
cyclic group $(\Zz/p\Zz)^{\times}$. Indeed, the number of subgroups of
the finite abelian group $E_{\ideal{p}}(\Fp_{\ideal{p}})$ with
$$
E_{\ideal{p}}(\Fp_{\ideal{p}})\simeq
\Zz/\dun\Zz \oplus \Zz/\dun\dde\Zz
$$
is ``essentially'' dominated by $\dun$ (see Birkhoff's description of
the subgroups of a finite abelian group, \cite[Th. 8.1]{birkhoff},
or~\cite[4.1.10]{cohen-2}), so $S_E(X;\dun)$ is closely related to the
sum 
$$
\sum_{N\ideal{p}\leq X}{T(\ideal{p})}
$$
where $T(\ideal{p})$ is the number of subgroups of
$E_{\ideal{p}}(\Fp_{\ideal{p}})$. The analogy between $S(X,d)$ and
$S_E(X;\dun)$ seems however deeper using the Galois-theoretic
interpretation. 
\end{remark}

More generally, for $\mathcal{C}\subset GL(2,\Zz/d\Zz)$ a set of
conjugacy classes, we will let
\begin{equation}\label{eq-pi-gen}
\pi_E(X;d,\mathcal{C})
=|\{\ideal{p}\,\mid\, N\ideal{p}\leq X,
\text{ and }\frob_{\ideal{p}}\mods{d} \in \mathcal{C}\}|.
\end{equation}
It is also convenient in many situations to weigh primes by $\log p$,
so we define also\footnote{ It would be better to consider here the
  partial sum of coefficients of the logarithmic derivative of the
  Artin $L$-function associated to the character of $G_d$
  which has trace equal to the characteristic function of
  $\mathcal{C}$.
}
\begin{equation}\label{eq-theta-gen}
\theta_E(X;d,\mathcal{C})
=\sum_{\stacksum{N\ideal{p}\leq X}
{\frob_{\ideal{p}}\in\mathcal{C}}}{\log N\ideal{p}}.
\end{equation}
Since we deal with all fields $K(E[d])$ at the same time, we 
use the shorthand notation $\frob_{\ideal{p}}\mods{d}$ to denote a
Frobenius element at $\ideal{p}$ for the field $K(E[d])$, so
$\frob_{\ideal{p}}\in G_d$; by convention, writing this implies also
that $\ideal{p}$ is unramified in $K(E[d])$. This notation is
compatible, in the case of the 
cyclotomic fields $\Qq(\mmu_d)$, with the usual meaning of
congruences and the isomorphism
$\Gal(\Qq(\mmu_d)/\Qq) \lra (\Zz/d\Zz)^{\times}$ which sends
$\frob_p$ to $p\mods{d}$.
\par
We see that~(\ref{eq-div-som1}) and~(\ref{eq-dun-som}) are comparable
in that both involve the average distribution of Frobenius elements in
the extensions generated by $d$-torsion points of some algebraic group
(either $E$ or the multiplicative group), uniformly for $d$ quite
large. However, there are a number of important qualitative
differences, as will be explained later on. Here we only mention that
the factor $\varphi(d)$ in~(\ref{eq-dun-som}) makes it impossible to
switch divisors there as in~(\ref{eq-trick}), making the contribution
of the very large moduli very hard to handle.
\par
The estimation of~(\ref{eq-dun-som}) seems to be a
much harder problem than the Titchmarsh divisor problem.

\begin{remark}
I have not seen any mention of the problem of estimating
$S_E(X;\dun)$ in the literature; however, there are a number of not
unrelated works concerning the question of counting primes $p\leq X$
such that $E_p(\Fp_p)$ is cyclic (i.e. $\dun(E_p)=1$) for an
elliptic curve $E/\Qq$, for instance~\cite{gupta-murty}. Also
Serre~\cite{serre-tch}, for counting supersingular primes $p\leq X$,
uses the fields of $\ell$-torsion with $\ell$ prime and quite large with
respect to $X$; however $\ell$ is fixed for a given $X$, so the
question of uniformity with respect to the modulus occurs in somewhat
attenuated form.
\end{remark}

\subsection{Analysis of the elliptic splitting problem on GRH}
\label{ssec-splitting-grh}

For \emph{fixed} $d\geq 1$, the asymptotic behavior of
$\pi_E(X;d,1)$ is given by the Chebotarev Density Theorem. Under GRH,
it can be stated in a sharp form. First 
we introduce some 
notation. As before, $E/K$ is an elliptic curve over a number field,
$d\geq 1$ an integer, $G_d$ is the Galois group of $K(E[d])$ over
$K$. Let $\Delta_d$ be the absolute value of the discriminant of
$K(E[d])/\Qq$, $n_1$ the degree $[K:\Qq]$, so
$[K(E[d]):\Qq]=|G_d|n_1$.
We let $N_E$ be the
norm of the conductor of $E/K$ (\cite[IV-10]{silv-2}).

\begin{proposition}\label{pr-tch}
Assume GRH for the Artin $L$-functions. With the above notation,
we have
\begin{equation}
\label{eq-tch-grh}
\pi_E(X;d,1)=\frac{1}{|G_d|}\li(x)+
\oun\bigl(\sqrt{X}\log (\Delta_1(d|G_d|N_EX)^{n_1})\bigr)
\end{equation}
for $X\geq 2$, with an absolute implied constant, and
\begin{equation}
\label{eq-tch-grh-t}
\theta_E(X;d,1)=\frac{X}{|G_d|}+
\oun\bigl(\sqrt{X}(\log x)(\log (\Delta_1(d|G_d|N_EX)^{n_1}))\bigr).
\end{equation}
for $X\geq 2$, with an absolute implied constant.
\end{proposition}

\begin{proof}
This is just making explicit the version given by
Serre~\cite{serre-tch}, based on that of Lagarias--Odlyzko, and
is well-known: we include the proof for completeness. Th\'eor\`eme~4
of~\cite{serre-tch} reads in this case
$$
\pi_E(X;C,1)=\frac{1}{|G_d|}\li(x)
+r_E(X;d)
$$
with the estimate
\begin{equation}\label{eq-rem-est}
r_E(X;d)\ll \frac{1}{|G_d|}\sqrt{X}(\log (\Delta_d)+
n_1|G_d|\log X),
\end{equation}
with an absolute implied constant. We have (see
e.g.~\cite[III]{serre-cl})
$$
\log \Delta_d=|G_d| \log \Delta_1 +\log N\ideal{d}_d
$$
where $\ideal{d}_d$ is the relative discriminant of $K(E[d])/K$. Then 
Proposition~5 of~\cite{serre-tch} gives an upper bound
$$
\log \Delta_d\leq |G_d|\log \Delta_1 + n_1|G_d|\Bigl(1-\frac{1}{|G_d|}
\Bigr)\log P_d+n_1|G_d|\log |G_d|,
$$
where $P_d$ is the product of the primes $p$ which are residue
characteristics of primes of $K$ ramified in $K(E[d])$. If $\ideal{p}$
is a prime of good reduction of $E$ and $(d,\ideal{p})=1$, 
$\ideal{p}$ is unramified in $K(E[d])$. It follows easily that
$$
P_d\mid dN_E.
$$
Thus we get
$$
\frac{\log \Delta_d}{|G_d|}\leq
\log \Delta_1 + n_1 \log dN_E + n_1 \log |G_d|.
$$
The first term in~(\ref{eq-rem-est}) is thus
$$
\frac{1}{|G_d|}\sqrt{X}\log (\Delta_d)\leq
\sqrt{X}\log (\Delta_1 (d|G_d|N_E)^{n_1}),
$$
so that we obtain
$$
r_E(X;d)\ll \sqrt{X}\log (\Delta_1 (d|G_d|N_EX)^{n_1})
$$
with an absolute implied constant.
\par
The proof for $\theta_E$ is similar or deduced by partial summation.
\end{proof}

\begin{remark}
If $K=\Qq$, this can be written
\begin{equation}\label{eq-tch-grh-q}
\pi_E(X;d,1)=\frac{1}{|G_d|}\li(x)+
\oun(\sqrt{X}\log (dN_EX)),
\end{equation}
(with an absolute implied constant) by observing that $|G_d|\leq d^4$
(for example), and one can replace $N_E$ by the absolute 
value of the discriminant of $E$, which it divides.
\end{remark}

For comparison, it is classical that GRH for Dirichlet $L$-functions
implies
\begin{equation}\label{eq-dir-grh}
\pi(X;d,a) = \frac{1}{\varphi(d)}\li(x)
+\oun(\sqrt{X}\log (dX))
\end{equation}
with an absolute implied constant.
\par
\medskip
Recall from Theorem~\ref{th-galois-nf} and the Remark following, that
as $d\ra +\infty$ the 
order of $G_d$ is comparable with $d^{\gale}$ where $\gale=2$ if $E$
has CM and with $\gale=4$ if not. Comparing the error term
in~(\ref{eq-tch-grh}) with $|G_d|$, it follows that~(\ref{eq-tch-grh})
gives the asymptotic behavior
$$
\pi_E(X;d,1)\sim \frac{1}{|G_d|}\li(x)\quad\text{as $X\ra +\infty$}
$$
\emph{uniformly} for $d$ up to $X^{1/(2\gale)-\eps}$ for any
$\eps>0$, whereas~(\ref{eq-dir-grh}) implies the corresponding
asymptotic for primes in arithmetic progression to moduli $d\leq
X^{1/2-\eps}$. Hence, since $1/(2\gale)=1/4$ (in the CM case) or
$=1/8$ (otherwise), we see a great difference for the purpose of
applying the estimates~(\ref{eq-tch-grh}) or~(\ref{eq-dir-grh}) to the
sums~(\ref{eq-dun-som}) and~(\ref{eq-div-som1}). In the case of the
Titchmarsh divisor problem, GRH provides an asymptotic formula valid
for ``almost all'' the moduli $d$ involved in~(\ref{eq-div-som1}),
leaving only those $d$ very close to $X^{1/2}$ to be dealt with; but
for an elliptic curve, a whole range of $d$ remains for which GRH does
not give anything, namely
$$
\begin{cases}
X^{1/4-\eps}\leq d\leq X^{1/2}+1&\text{if $E$ has CM}\\
X^{1/8-\eps}\leq d\leq X^{1/2}+1&\text{if $E$ does not have CM.}
\end{cases}
$$
(it is certainly not surprising that the non-CM case appears
superficially to be worse than the other, although whether it should
really be is open to question...)
\par
However, we can at least state what this gives
for~(\ref{eq-dun-som}).

\begin{proposition}\label{pr-part}
Let $E/K$ be an elliptic curve over a number field. Assume GRH for
Artin $L$-functions. Then we have
\begin{align}\label{eq-dun-part-cm}
\sum_{d\leq \frac{X^{\scriptscriptstyle{1/4}}}{\log X}}
{\varphi(d)\pi_E(X;d,1)}= c(E)X
+\oun\Bigl(\frac{X}{(\log X)^2}
\log(\Delta_1N_E^{n_1}X^{3n_1+1})
\Bigr)
\quad\text{if $E$ has CM,}\\
\sum_{d\leq \frac{X^{\scriptscriptstyle{1/4}}}{(\log X)^2}}
{\varphi(d)\pi_E(X;d,1)}=c(E)\li(X)
+\oun\Bigl(\frac{X}{(\log X)^4}
\log(\Delta_1N_E^{n_1}X^{5n_1+1})
\Bigr)
\quad\text{otherwise}
\label{eq-dun-part-noncm}
\end{align}
for $X\geq 2$, with absolute implied constants, where
\begin{align}
c(E)&=\res_{s=0}
{\sum_{d\geq 1}{\frac{\varphi(d)}{|G_d|}d^{-s}}}\quad\text{if $E$ has
  CM}\label{eq-const-cm}\\ 
c(E)&=\sum_{d\geq 1}{\frac{\varphi(d)}{|G_d|}}\quad
\text{otherwise.}\label{eq-const-noncm}
\end{align}
Unconditionally, we have a lower bound
\begin{equation}\label{eq-uncond-lower}
S_E(X;\dun)\gg_K \frac{X}{\log X},
\end{equation}
where the implied constant depends only on $K$.
\end{proposition}


\begin{proof}
This is an immediate corollary of Proposition~\ref{pr-tch}.
Take the non-CM case for example: we have
$$\
\frac{\varphi(d)}{|G_d|}\leq \frac{1}{d^2\varphi(d)}
$$
so the series defining $c(E)$ is absolutely convergent, and the main
term of~(\ref{eq-tch-grh}) gives
$$
\li(X)\sum_{d\leq X^{1/4}/(\log X)}{\frac{\varphi(d)}{|G_d|}}
=c(E)\li(X)+\oun_{\eps}(X^{1/2+\eps})
$$
(for any $\eps>0$, say we take $\eps=1/4$),
while for the error term we have
$$
\sqrt{X}\sum_{d\leq X^{1/4}/(\log X)^2}{\varphi(d)
\log (\Delta_1(d|G_d|N_EX)^{n_1})}
\ll \frac{X}{(\log X)^4}
\log(\Delta_1N_E^{n_1}X^{5n_1+1})
$$
by trivial summations (using $|G_d|\leq d^4$).
The CM case is exactly similar, except that the series over
$d$ has logarithmic growth, hence the different formula for $c(E)$.
\par
The lower bound~(\ref{eq-uncond-lower}) is an immediate consequence of
the Prime Ideal Theorem in $K$ since by~(\ref{eq-dun-som})
$$
S_E(X;\dun)\geq \pi_E(X;1,1)=\pi_K(X) \gg_K \frac{X}{\log X}
$$
where $\pi_K(X)$ is the number of prime ideals with norm $\leq X$.
\end{proof}

\begin{remark}
Note that the restriction to $d\leq X^{1/4}$ comes from the occurrence
of $\varphi(d)$ in~(\ref{eq-dun-som}). The exponent is thus
independent of the Galois dimension of $E$, and so of the actual range
where~(\ref{eq-tch-grh}) gives an asymptotic formula for
$\pi_E(X;d,1)$. In other words, in the non-CM case, in part of the
summation range in~(\ref{eq-dun-part-noncm}), the estimated term in
the Chebotarev density theorem dominates over the main term.
\par
Note that the constant $c$ in~(\ref{eq-tit-asy}) is also
$$
c=\res_{s=0}\sum_{d\geq 1}{\frac{1}{\varphi(d)d^s}}.
$$
and the same argument gives
$$
\sum_{d\leq \sqrt{X}/(\log X)}{\pi(X;d,1)}
\sim cX\quad\text{as $X\ra +\infty$},
$$
leaving only the range $\sqrt{X}/(\log X)\leq d\leq \sqrt{X}$ to
handle to solve (under GRH) the Titchmarsh divisor problem.
\end{remark}

It is reasonable to expect that the sum in Proposition~\ref{pr-part}
could be extended to all $d\leq \sqrt{X}+1$, giving the desired
asymptotic formula for the average of $\dun(\ideal{p})$ over
$\ideal{p}$. 



\subsection{Computation of $c(E)$}
\label{sec-comp-ce}

In Section~\ref{sec-numerics} below we perform numerical experiments
for the elliptic splitting problem, and it is therefore useful to
be able to explicitly evaluate the constant $c(E)$, at least for some
elliptic curves $E$. This requires some knowledge of the Galois
groups $G_d$, which is available in the case of what Lang-Trotter call
\emph{Serre curves} (\cite[I, \S 5-6-7]{lang-trotter}). Serre~\cite[\S
5]{serre-div} has indeed given concrete examples of such curves, and
we will use his examples in Section~\ref{sec-numerics}. Throughout
this section, all curves are over $\Qq$.
\par
The difficulty in computing $|G_d|$, and hence $c(E)$, is that
although the index between them is bounded, it is never the case that
$G_d=GL(2,\Zz/d\Zz)$ for all $d\geq 1$, as shown by
Serre. More precisely, let
$E[\infty]$ be the set of all torsion points of $E$, and 
$$
\rho_{\infty}\,:\,G_{\Qq}\lra \Aut(E[\infty])
$$
the natural Galois representation, so that
$\rho_{\infty}\mods{d}=\rho_d$ for 
all $d\geq 1$. Recall that
$$
\Aut(E[\infty])=\prod_{\ell}{GL(2,\Zz_{\ell})}
$$
and the $\ell$-th component of $\rho_{\infty}$ is the $\ell$-adic
representation $\rhol$.
\par
Define an index $2$ subgroup $H_E$ of $\Aut(E[\infty])$ as follows:
let $\eps\,:\,GL(2,\Zz_2)\ra \{\pm 1\}$ be the map given by
composition
$$
GL(2,\Zz_2)\ra GL(2,\Zz/2\Zz)\simeq \mathfrak{S}_3
\stackrel{\eps}{\lra} \{\pm 1\}
$$
where $\eps$ is the signature on $\mathfrak{S}_3$. Let $\chi$ be the
Kronecker symbol of the quadratic extension $\Qq(\sqrt{\Delta})$,
where $\Delta$ is the discriminant of $E$, and $m$ its conductor. The
subgroup in question is defined by
$$
H_E=\{g=(g_{\ell})\in \Aut(E[\infty])\,\mid\,
\eps(g_{2})=\chi(g\mods{m})\}
$$
Then the precise form of Serre's result (\cite[Prop. 22]{serre-div})
is: 
\begin{proposition}(Serre)
For any elliptic curve $E/\Qq$ we have $\rho_{\infty}(G_{\Qq})\subset
H_E$.
\end{proposition}

By definition, a \emph{Serre curve} is an elliptic curve $E/\Qq$ such
that $\rho_{\infty}(G_{\Qq})=H_E$ (see~\cite[I, \S 5]{lang-trotter}
for a more detailed discussion, Section~\ref{sec-numerics} for
concrete examples).

\begin{proposition}\label{pr-galois-serre}
Let $E/\Qq$ be a Serre curve, and let $m$ be as above. We have
$$
[GL(2,\Zz/d\Zz):G_d]=
\begin{cases}
2&\text{ if $2m\mid d$}\\
1&\text{ otherwise.}
\end{cases}
$$
\end{proposition}

\begin{proof}
Clearly we have
$$
H_E=H_{2m} \times \prod_{(\ell,m)=1}{GL(2,\Zz_{\ell})},
$$
where $H_{2m}$ is the obvious subgroup (the definition of $H_E$ only
involves the components of $g$ at $\ell\mid 2m$). 
Let $g=(g_{\ell})$ be a representative of the
non-trivial coset of $H_{2m}$. Correspondingly, if $d=d_1d_2$ with
$d_1\mid (2m)^{\infty}$ and $(d_2,2m)=1$, we have
$$
G_d=G_{d_1}\times GL(2,\Zz/d_2\Zz).
$$
So it is enough to compute the index of $G_{d_1}$.
Since $H_E$ is of index $2$ in $\Aut(E[\infty])$, it is either $1$ or
$2$. Now if $2m\mid d_1$, the reduction modulo $d$ of $g$ is an element
in $G_{d_1}$ which is not in $H_{2m}\mods{d_1}$, so the index is $2$
in this case.
\par
Conversely, if $2m$ does not divide $d_1$, let $\ell$ be a prime
dividing $2m$ but not $d_1$. For any $g\in GL(2,\Zz/d_1\Zz)$, we can
lift it to 
$$
\prod_{\ell'\mid d_1}{GL(2,\Zz_{\ell'})}
$$
and then change the component at $\ell$ so that the resulting $\hat{g}$ is
in $H_{2m}$; this element reduces to $g$ modulo $d_1$, so the index of
$G_{d_1}$ is $1$ in this case.
\end{proof}

\begin{lemma}\label{lm-arith-sum}
Let $f$ and $g$ be arithmetic functions with $g$ multiplicative such
that
$$
f(d)=
\begin{cases}
\alpha g(d)&\text{ if $n\mid d$}\\
g(d)&\text{otherwise}
\end{cases}
$$
for some integer $n\geq 1$ and some $\alpha\in \Rr$. Assume moreover
that \begin{equation}\label{eq-assump}
g(nd)=d^{-\kappa}g(n)
\end{equation}
for all $d\mid n^{\infty}$ and some real number $\kappa$. Assume that
the series $\sum{g(d)}$ converges absolutely. Then we have
$$
\sum_{d\geq 1}{f(d)}=c\prod_{p}{g_p}
$$
where
$$
g_p=\sum_{k\geq 0}{g(p^k)},
$$
and
$$
c=1+(\alpha-1)g(n)\prod_{p\mid n}{g_p^{-1}(1-p^{-\kappa})^{-1}}.
$$
\end{lemma}

\begin{proof}
We compute, from the assumption:
\begin{align*}
\sum_{d\geq 1}{f(d)}&=\alpha\sum_{n\mid d}
{g(d)}+\sum_{n\nmid d}{g(d)}\\
&=\alpha\sum_{n\mid d}{g(d)}+
\sum_{d\geq 1}{g(d)}-\sum_{n\mid d}{g(d)}\\
&=\sum_{g\geq 1}{g(d)}+(\alpha-1)\sum_{n\mid d}{g(d)}.
\end{align*}
By multiplicativity we have
$$
\sum_{d\geq 1}{g(d)}=\prod_p{g_p}.
$$
Factorizing uniquely $d=d_1d_2$ with $d_1\mid n^{\infty}$ and
$(d_2,n)=1$, we have further
\begin{align*}
\sum_{n\mid d}{g(d)}&=\sum_{d\geq 1}{g(nd)}\\
&=\sum_{\stacksum{d_1\mid n^{\infty}}{(d_2,n)=1}}
{g(nd_1d_2)}\\
&=\Bigl(\sum_{(d_2,n)=1}{g(d_2)}\Bigr)
\Bigl(\sum_{d_1\mid n^{\infty}}{g(nd_1)}\Bigr)\\
&=g(n)\Bigl(\prod_{(p,n)=1}{g_p}\Bigr)
\Bigl(\sum_{d_1\mid n^{\infty}}
{d_1^{-\kappa}}\Bigr)
\text{ (by multiplicativity and~(\ref{eq-assump}))}\\
&=g(n)\prod_{p\mid n}{(1-p^{-\kappa})^{-1}}
\prod_{(p,n)=1}{g_p},
\end{align*}
whence the result follows.
\end{proof}

\begin{corollary}\label{cor-const-serre}
Let $E/\Qq$ be a Serre curve. We have
$$
c(E)=\sum_{d\geq 1}{\frac{\varphi(d)}{|G_d|}}
=c'(E)\zeta(2)\zeta(3)
\prod_{p}{(1-p^{-2}+p^{-5})}
$$
with
$$
c'(E)=1+\frac{1}{(2m)^{3}}
\prod_{p\mid 2m}{(1-p^{-2}+p^{-5})^{-1}}.
$$
\end{corollary}

\begin{proof}
In view of Proposition~\ref{pr-galois-serre}, we can apply
Lemma~\ref{lm-arith-sum} with $n=2m$, $\alpha=2$ and
\begin{align*}
f(d)&=\frac{\varphi(d)}{|G_d|}\\
g(d)&=\frac{\varphi(d)}{|GL(2,\Zz/d\Zz)|}.
\end{align*}
Indeed~(\ref{eq-assump}) holds with $\kappa=3$ since more generally we
have by~(\ref{eq-ord-ncm}),~(\ref{eq-psi})
$$
g(dd_1)=(dd_1)^{-3}\prod_{p\mid dd_1}{p^{-1}(1-p^{-2})^{-1}}
=d_1^{-3}g(d)
$$
if $d_1$ has no prime divisor outside $d$ (this formula explains where
functions satisfying~(\ref{eq-assump}) arise naturally).
\par
We have by~(\ref{eq-ord-ncm})
\begin{align*}
g_p&=1+\sum_{k\geq 1}{\frac{1}{p^{3k}(1-p^{-1})(1+p^{-1})}}\\
&=1+\frac{1}{p^3(1-p^{-2})(1-p^{-3})}\\
&=\frac{p^5-p^3+1}{(p^2-1)(p^3-1)},
\end{align*}
hence the result after some rearranging of terms.
\end{proof}

\begin{remark}
Note that the correction factor $c'(E)$ is usually very close to $1$, so
the value of $c(E)$ for a Serre curve is close to
\begin{equation}\label{eq-c0}
c_0=
\zeta(2)\zeta(3)
\prod_{p}{(1-p^{-2}+p^{-5})}=1.25845\ldots
\end{equation}
This means in particular that if the expected asymptotic formula for
$S_E(X;\dun)$ holds, by itself it does not carry much \emph{global}
information about $E$, except for distinguishing between CM curves and
non-CM curves.
\end{remark}

\begin{remark}
One may hope that this gives the ``generic'' value of $c(E)$. More
precisely, recall that Duke~\cite{duke-nex}
has shown that for ``almost all'' elliptic curves over $\Qq$ (in the
sense of almost all coefficients of Weierstrass equations), there
are no ``exceptional primes'', i.e. we have
$$
G_{p}=GL(2,\Zz/p\Zz)
$$
for all primes $p$. It may be possible to refine this statement
to show that almost all $E/\Qq$ (in the same sense) are Serre curves.

\end{remark}


\subsection{Outside primes}

The simple-minded analysis based on GRH of the previous section points
to a striking difference between the distribution of totally split
primes in $K(E[d])$ for large modulus $d$ and the case of arithmetic
progressions. This is best made explicit using
$$
\psi(X;d,1)=\sum_{\stacksum{n\leq X}{n\equiv a\mods{d}}}
{\Lambda(n)}
$$
where $\Lambda(n)$ is the von Mangoldt function, equal to $\log p$ if
$n=p^m$ for some prime $p$ and $m\geq 1$, and to $0$ otherwise. As for
$\theta_E$, we have on GRH
\begin{equation}\label{eq-pr-arith-grh}
\psi(X;d,a)=\frac{X}{\varphi(d)}+\oun(\sqrt{X}(\log X)(\log dX))
\end{equation}
for $X\geq 2$.
\par
Now, consider the smallest prime $\equiv a$ modulo $d$, or the smallest
$X$ for which $\psi(X;d,a)>0$. Since
$p\equiv a\mods{d}$ implies $d\leq p-a$, it follows that 
$p\geq d+a>d$, in particular the main term of~(\ref{eq-pr-arith-grh})
is $>1$, i.e. we have
$$
\psi(X;d,1)>0 \Rightarrow X>d \Rightarrow \frac{X}{\varphi(d)}>1.
$$
We restate this as follows: all primes in arithmetic progression
can be ``accounted for'' by the main term in the Chebotarev density
theorem. Such is still the case of CM elliptic curves, since the \emph{a
priori} estimate~(\ref{eq-a-priori}) shows that
\begin{equation}\label{eq-small}
\theta_E(X;d,1)>0 \Rightarrow X\geq (d-1)^2
\end{equation}
which is (roughly) compatible with the density $1/|G_d|\geq 1/d^2$ in
this case. 
\par
Non CM curves are different: the estimate~(\ref{eq-a-priori}) is the best
general bound 
(as shown below), but now the density of totally split primes is roughly
$1/d^{\gale}=1/d^4$. If $\ideal{p}$ splits in $K(E[d])$ with $N\ideal{p}<
|G_d|$, the main term in the Chebotarev density theorem is $<1$, and
this may be the case for values of $d$ as large as
$\sqrt{N\ideal{p}}+1$. Such a prime is \emph{not} accounted for by the
main term of the Chebotarev density theorem.
\begin{definition}
Let $E/K$ be a non-CM elliptic curve over a number field $K$. A prime
ideal $\ideal{p}$ which splits totally in $K(E[d])$ with
$N\ideal{p}<|G_d|$ is called an \outside prime of
$E$. If $\ideal{p}$ satisfies the weaker inequality $N\ideal{p}<d^4$,
it is called a weak \outside prime.
\end{definition}
Equivalent formulations are $|G_{\dun(\ideal{p})}|> N\ideal{p}$ and
$\dun(\ideal{p})>(N\ideal{p})^{1/4}$ respectively.
\par
The existence of \outside primes is understandable: since the
invariant $\dun(\ideal{p})$ only depends on the reduction of $E$
modulo $\ideal{p}$, it follows that for given $\ideal{p}$, $E$ being
globally CM or not does not matter. The results on the possible group
structures of elliptic curves over finite fields (see
Section~\ref{sec-local}) show that the a priori
bound~(\ref{eq-a-priori}) is always best possible.
\par
We give here a simple illustrative example.
\begin{example}\label{ex-cm-curve}
Let $A/\Qq$ be the classical CM curve given by the Weierstrass
equation
\begin{equation}
  \label{eq-cm}
y^2=x^3-x
\end{equation}
which has $j(A)=1728$, conductor $N_A=32$ and endomorphism ring
$\End_{\bar{\Qq}}(A)=\Zz[i]$, the ring of Gaussian integers.
\par
The determination of the local Frobenius endomorphism of $A$ modulo
$p$, up to conjugation, is classical (see
e.g.~\cite[18.4]{ire-ros}). 
If $p\equiv 3\mods{4}$, $A$ is supersingular at $p$ and $a_p(A)=0$. If
$p\equiv 
1\mods{4}$, on the other hand, $p$ splits in $\Qq(i)$, say
$p=\pi\bar{\pi}$ for some prime element $\pi$, and the Frobenius at
$p$ is one of the elements $\pm \pi$, $\pm i\pi$, $\pm\bar{\pi}$, $\pm
i\bar{\pi}$. Which one it is, up to conjugation,
is settled by a congruence modulo $2(1+i)$, namely
\begin{equation}\label{eq-congr-prim}
\frob_p\equiv 1\mods{2(1+i)}
\end{equation}
(a Gaussian integer $z\equiv 1\mods{2(1+i)}$ is called primary).
To see this, one can either express $a_p$ in terms of
Jacobstahl sums and reduce modulo $2(1+i)$ (see
e.g.~\cite[8.2]{iwaniec-classical}) or observe that the
$2(1+i)$-torsion of $A$ is 
rational over $\Qq(i)$, hence over $\Fp_p$ for $p$ split in $\Qq(i)$,
so that~(\ref{eq-congr-prim}) follows ($A[2(1+i)]$ is
generated by the two-torsion points $(0,0)$, $(\pm 1,0)$ and by
$(i,\pm (i-1))$; see e.g.~\cite[Ex. 12.3]{rubin}).
\par
\par
Now if $\pi$ is a Gaussian prime of the form $\pi=1+ni$ such that
$\pi\equiv 1\mods{2(1+i)}$, then $p=n^2+1$ is prime and $\pi$ is the
Frobenius at $p$. But tautologically we have $\pi\equiv 1\mods{n}$ in
$\End(A_p)$, so that (Lemma~\ref{lm-cong-end}) $\dun(A_p)=n$, and in
fact, since $N(\pi-1)=n^2$, $A_p(\Fp_p)\simeq (\Zz/n\Zz)^2$
(compare~\cite[2.5]{schoof-1}). Obviously $n=[\sqrt{p}+1]$. 
\par
In terms of $p$ the condition is that $p=n^2+1$ and $4\mid n$ (i.e
$p=16n^2+1$). It is
expected that there exist infinitely many primes $p$ of this form,
but this is not known (see~\cite{iwa-p2} for the best ``almost prime''
results). The first few are $p=17$, $257$, $401$, $577$,..., $739601$,
...
\par
Now if $p$ is a prime of this type, any curve $E'/\Qq$ with the same
reduction modulo $p$ as $E$ will also have $\dun(E')=n$. For
instance, take
$$
E'\,:\, y^2=x^3-x-p
$$
which for all $p$ does not have CM and for $p=16n^2+1$ will
have $\dun(E'_p)=[\sqrt{p}+1]$ by construction.
\end{example}

Obviously, for the purpose of finding an asymptotic evaluation
of~(\ref{eq-sum-dun}), a few prime ideals with $\dun(\ideal{p})$ close
to $\sqrt{N\ideal{p}}$ do not matter much. One might expect that in
general \outside primes are rare, and the presence of ``too many'' of
them should mean that $E$ has CM.
\par
A partial clue in this direction is implicit
in~\cite[p. 330]{schoof-1}. We state the following simple result as an
illustration: it shows that Example~\ref{ex-cm-curve} is basically the only
possibility in the most extreme case.

\begin{proposition}
\label{pr-cm-cong}
Let $E/\Qq$ be an elliptic curve with $j$-invariant $j$, $p\geq 11$ a
prime of good reduction of $E$ such that
$$
\dun(p)\geq \sqrt{\frac{p}{2}}.
$$
Then $j\equiv j_0\mods{p}$, where
$$
j_0\in \mathcal{J}=\{
0, 1728, -3375, 8000, -32768, 54000 \}.
$$
In particular, there are only finitely many such $p$ unless $j\in
\mathcal{J}$. In this case $E$ is a CM curve.
\end{proposition}

\begin{proof}
First observe that the reduced curve $E_p/\Fp_p$ is ordinary. Let
$\pi\in \Oc=\End_{\bar{\Fp}_p}(E_p)$ be the Frobenius endomorphism. We
have (Lemma~\ref{lm-cong-end}) $\pi=1+\dun(E_p)\pi'$ for some $\pi'\in
\Oc$, and
$$
|E_p(\Fp_p)|=N(\pi-1)=\dun(E_p)^2N\pi'.
$$
Moreover, since $\pi\not\in \Zz$, $\pi'$ is not in $\Zz$ either.
Let $D$ be the discriminant of the quadratic imaginary order
$\Oc$. For any $z\in\Oc$, $z\not\in \Zz$, we have
$$
Nz\geq \frac{|D|}{4},
$$
and applying this to $\pi'$ we get
\begin{align*}
|D|&\leq 4\frac{|E_p(\Fp_p)|}{\dun(E_p)^2}\\
&\leq 8\frac{|E_p(\Fp_p)|}{p}\text{ (by assumption)}\\
&\leq 8\Bigl(1+\frac{1}{\sqrt{p}}\Bigr)<15
\text{ (since $p\geq 11$).}
\end{align*}
But all quadratic imaginary orders of discriminant $<15$ have class
number one (see e.g.~\cite[Th. 7.30]{cox}). Now Deuring~\cite{deuring}
has shown 
that an ordinary elliptic curve $A$ over a finite field $\Fp_q$
``lifts to characteristic 
$0$''. This means that there exists a number field $K$, a prime ideal
$\ideal{p}$ of $K$ with $\Fp_{\ideal{p}}=\Fp_q$, and an
elliptic curve $\tilde{A}/K$ with CM by $\End(A)$ such that
$\tilde{A}_{\ideal{p}}\simeq A$.
\par
Let $\tilde{E}$ be such a lift of $E_p$. It has CM by the
order $\Oc$ with class number one, hence (see e.g.~\cite[II-2]{silv-2}) is
defined over $\Qq$, and actually a
table such as that in~\cite[12-C]{cox} or~\cite[App. A-3]{silv-2},
shows that $j(\tilde{E})\in \mathcal{J}=\{
0, 1728, -3375, 8000, -32768, 54000 \}$. Since
$\tilde{E}_p\simeq E_p$, we have $j\equiv j(\tilde{E})\mods{p}$.
\end{proof}

Obviously this argument can be extended somewhat, but it seems hard to
make interesting conclusions in greater generality. The difficulty is
roughly as follows: say we want to estimate the number of $\ideal{p}$
with $\dun(\ideal{p})\geq N\ideal{p}^{\theta}$ for some $\theta>0$
(for example, $\theta=1/4$, corresponding essentially to \outside
primes). As above one derives 
$$
|D|\leq 8p^{1-2\theta}
$$
where $D$ is the discriminant of the quadratic order
$\End(\Fp_{\ideal{p}})$. This implies
$$
j(E)\mods{\ideal{p}}\in
\Omega(\ideal{p})
$$
for some finite set $\Omega(\ideal{p})$ with
$$
|\Omega(\ideal{p})|\ll p^{1-2\theta}
$$
for all such $\ideal{p}$. However since the cardinality of
$\Omega(\ideal{p})$ is not bounded anymore, it is hard to go
further.
\par
Indeed, compare this to the analogue approach to the study of
supersingular primes of $E$: if $\ideal{p}$ is a prime of
supersingular reduction, we have $j(E)\mods{\ideal{p}}\in
\Omega'(\ideal{p})$, where $\Omega'(\ideal{p})$ is the finite set of
supersingular $j$-invariants. Lang and Trotter, who initiated the
study of the set of supersingular primes of elliptic curves,
explicitly mention this idea and state~\cite[p. 7]{lang-trotter} that it
doesn't seem to bring useful results.
\par
We thus have the following problem:
\begin{problem}\label{prob-outside}
Let $E/K$ be an elliptic curve without CM. What can one say about the
distribution of \outside primes of $E$? Are there infinitely many of
them? If yes, how many are there $\leq X$? Is is true that the series
$$
\sum_{\ideal{p}}{\frac{1}{N\ideal{p}}}
$$
over \outside primes of $E$ converges?
\end{problem}

The first guess, for $E/\Qq$, might be that there are
infinitely many \outside primes. Heuristically from
Proposition~\ref{pr-large-poids}, one would expect that there
are at most about $X^{1/4}$ \outside primes $\leq X$.
See Section~\ref{sec-numerics} for some numerical data: \outside
primes appear to be extremely scarce and Section~\ref{sec-twins} below
for a first idea.


\begin{remark}
Another seemingly simpler situation where ``outside'' primes can
occur, which throws some light on the situation, is that of Kummer
extensions. For simplicity, let $a\in\Zz$ be a squarefree number. For
$d\geq 1$, let $K_d=\Qq(\mmu_d,a^{1/d})$ be the Kummer extension
generated by $d$-th roots of $a$. As is well-known, we have in this
case an isomorphism
$$
\Gal(K_d/\Qq)\simeq (\Zz/d\Zz)^{\times} \rtimes (\Zz/d\Zz).
$$
\par
The order of the Galois group is thus $d\varphi(d)$, and one can define 
an \emph{outside} prime for $a$ to be $p$ such that $p$ splits
completely in $K_d$ with $p<d\varphi(d)$. 
\par
It is easy to see that, given $p$, the largest $d$ for which $p$
splits completely in $K_d$ is $d=(p-1)/o_p$ where $o_p$ is the
multiplicative order of $a$ modulo $p$ : indeed we have $p\equiv
1\mods{d}$, and $a^{(p-1)/d}=a^{o_p}=1\mods{p}$ so $a$ is a $d$-th power
modulo $p$.
\par
Hence $p$ is an outside prime if and only if 
$$
p<\frac{p-1}{o_p}\varphi\Bigl(\frac{p-1}{o_p}\Bigr).
$$
\par
Roughly speaking this is true if $o_p\ll \sqrt{p}$ (or equivalently if
$p\mid a^{j}-1$ with $j\ll \sqrt{p}$). Thus the question is clearly
related to Artin's conjecture about primitive roots and is currently
much of a mystery. Getting non-trivial results seems extremely
difficult, and one might expect the (non-CM) elliptic curve case to be
also very hard.
\end{remark}

\subsection{Brun-Titchmarsh problems}\label{ssec-bt}

In the study of the Titchmarsh divisor problem, 
to obtain a proof of~(\ref{eq-tit-asy}) requires dealing with the
large moduli 
$\sqrt{X}/(\log X)\leq d\leq \sqrt{X}$. 
Asymptotic formulae are not known in this range and do not follow from
GRH (although they are conjectured to hold for $d\leq X^{1-\eps}$, see
e.g.~\cite{granville}), but one can prove by
sieve methods upper bounds of the correct 
order of magnitude which are sufficient to derive the asymptotic
formula from that given by GRH (or, unconditionally, by the
Bombieri-Vinogradov Theorem). This was first done by
Titchmarsh~\cite{titch-div}).

\begin{theorem}\label{th-bt-classic}
For all $d\geq 1$, all $a$ with
$(a,d)=1$, and any $\eps>0$ we have
\begin{equation}\label{eq-bt-classic}
\pi(X;d,a)\ll_{\eps} \frac{X}{\varphi(d)\log X}
\end{equation}
for $d\leq X^{1-\eps}$, the implied constant depending only on $\eps$.
\end{theorem}

A sharp version has been proved by Montgomery-Vaughan~\cite{mon-vau}:
\begin{equation}\label{eq-bt-strong}
\pi(X;d,a)\leq \frac{2X}{\varphi(d)\log X/d}
\end{equation}
for all $d<X$.
\par
We recall for convenience how, using~(\ref{eq-bt-classic}), one can now
finish the proof of~(\ref{eq-tit-asy}) on GRH
from~(\ref{eq-div-som1}). Indeed, one has
\begin{align}
\sum_{X^{1/2-\delta}\leq d \leq\sqrt{X}}{\pi(X;d,1)}
&\ll \frac{X}{\log X}\sum_{\sqrt{X}/(\log X)\leq d\leq X}
{\frac{1}{\varphi(d)}}\nonumber\\
&\ll \frac{X\log\log X}{\log X}\label{eq-fin-tit}
\end{align}
for $X\geq 2$, and similarly
$$
\sum_{d<\sqrt{X}}{\pi(d^2+1;d,1)}\ll 
\sum_{d<\sqrt{X}}{\frac{d^2}{\varphi(d)\log d}}\ll \frac{X}{\log X},
$$
for $X\geq 2$.
\par
\medskip
This naturally suggests the following problem:

\begin{problem}\label{pr-bt-ell}
Let $E/K$ be an elliptic curve over a number field $K$. Is it true
that for any $\eps>0$ there exists $C(E,\eps)>0$ such that
\begin{equation}\label{eq-bt-ell}
\pi_E(X;d,1)\leq C(E,\eps)\frac{X}{|G_d|(\log X)}
\end{equation}
for all $d\leq X^{1/g-\eps}$?
\end{problem}

Note that the restriction to $d\leq X^{1/g}$ is certainly necessary,
since for larger $d$ the ``main term'' of the Chebotarev density
theorem is $<1$ (for $X$ large enough). See below for further
discussion of this point.
\par
There's a remark that arises in writing such an inequality: should one
write $|G_d|$, in the denominator, or instead, assuming that $E$ is non-CM,
$|GL(2,\Zz/d\Zz)|$? Both forms are equivalent, because of
Serre's result that the index of $G_d$ in $GL(2,\Zz/d\Zz)$ is bounded.
But in fact an inequality
$$
\pi_E(X;d,1)\leq C(E,\eps)\frac{X}{|GL(2,\Zz/d\Zz)|(\log X)}
$$
for $d\leq X^{1/g-\eps}$ \emph{implies} Serre's result: fix
$d$, take $\eps=1/(2g)$ (say), so for all $X\geq d^{2g}$ we have
$$
\pi_E(X;d,1)\leq C(E,(2g)^{-1})  \frac{X}{|GL(2,\Zz/d\Zz)|(\log X)},
$$
whereas by the Chebotarev density theorem
$$
\pi_E(X;d,1)\sim \frac{X}{|G_d|(\log X)}
$$
as $X\ra +\infty$. Comparing implies
$$
[GL(2,\Zz/d\Zz):G_d]\leq C(E,(2g)^{-1}).
$$
Now it is interesting to note that the Brun-Titchmarsh
inequality~(\ref{eq-bt-classic}) is proved, with $\varphi(d)$ in the
denominator, without any mention of cyclotomic fields! The same
argument backwards then deduces from~(\ref{eq-bt-classic}) that
the index of the Galois group of $\Qq(\mmu_d)$ in
$(\Zz/d\Zz)^{\times}$ is bounded (by $2$,
using~(\ref{eq-bt-strong})). Of course, it is not hard to 
prove that it is $1$ for all $d$ (i.e. the cyclotomic polynomials are
irreducible).\footnote{ Any constant $<2$ 
in~(\ref{eq-bt-classic}) would reprove this, but it is
well-known (see references in~\cite[p. 123]{halb-rich}) that such a result
would bring much richer rewards, as it would eliminate the possibility
that the so-called Landau-Siegel zeros of quadratic Dirichlet
$L$-functions exist.}

\begin{proposition}
Let $E/K$ be an elliptic curve over a number field. Assume
that~(\ref{eq-bt-ell}) holds for $E$ in the range stated. Then we have
\begin{align}
\sum_{d\leq X^{1/g-\eps}}{\varphi(d)\pi_E(X;d,1)}&
\ll X\quad \text{(if $E$ has CM)}\\
\sum_{d\leq X^{1/g-\eps}}{\varphi(d)\pi_E(X;d,1)}&
\ll \frac{X}{\log X}\quad\text{(otherwise)}
\end{align}
for any $\eps>0$ and any $X\geq 2$, the implied constant depending
only on $E$ and $\eps$.
\end{proposition}
Note this is weaker than what GRH implies (Proposition~\ref{pr-part}),
but it may be the case
that~(\ref{eq-bt-ell}) is easier to prove, as in the cyclotomic case.
The proof is immediate, and the statement is given only for
completeness.
\par
\medskip
It is clear that the Brun-Titchmarsh problem for $K(E[d])$ 
can be much generalized. Let us consider the following rather general
context (compare~\cite[]{serre-tch}): let $K$ be a number field and
$K'/K$ an infinite Galois extension which is unramified outside a
finite set of primes of $K$, and has Galois 
group $\Gal(K'/K)$ which is (isomorphic to) a 
finite index subgroup of $G(\hat{\Zz})$ for 
some smooth algebraic group $G$ of finite type over $\Zz$. For $d\geq 1$, let
$K_d/K$ be the fixed field of the kernel of the reduction modulo $d$
map
$$
\Gal(K'/K)\injecte G(\hat{\Zz}) \surjecte G(\Zz/d\Zz),
$$
a Galois extension of $K$ with $\Gal(K_d/K)=\Gal(K'/K)\mods{d}$,
with obvious notation. The Galois groups $\Gal(K_d/K)$ are, by the map
above, subgroups of $G(\Zz/d\Zz)$ with index bounded for $d\geq 1$.
Let $\gale$ be the (relative) dimension of $G/\Zz$. 

\begin{definition}
With notation as above, the field $K'$ is a \emph{Brun-Titchmarsh}
field if and only if for any $\eps>0$ we have
\begin{equation}\label{eq-bt-generale}
\pi_{K'}(X;d,1)\ll_{K',\eps} \frac{1}{|G(\Zz/d\Zz)|}
\frac{X}{\log X}
\end{equation}
if $d\leq X^{1/\gale-\eps}$, the implied constant depending only on
$K'$ and $\eps$, where $\pi_{K'}(X;d,1)$ is the number of prime ideals
of $K$ with norm $\leq X$ which are totally split in $K_d$.
\end{definition}

So the cyclotomic extension $\Qq^{ab}/\Qq$ is a Brun-Titchmarsh
field, and Problem~\ref{pr-bt-ell} can be re\-phrased as asking whether
the field $K(E[\infty])=\bigcup_d K(E[d])$ is a Brun-Titchmarsh
field. Other examples arise naturally: for the same $E/K$, not CM, 
let $K'\subset K(E[\infty])$ be the subextension corresponding to the
closed subgroup $Z(\hat{\Zz})\cap \Gal(K(E[\infty]))$, where $Z$ is
the center of $GL(2)$. It has Galois group $G$ which is of finite
index in $PGL(2,\hat{\Zz})$, hence $\gale=3$ in this case.
If $K'/K$ were a Brun-Titchmarsh field, and assuming GRH for Artin
$L$-functions, the asymptotic formula
$$
\sum_{d\leq X^{1/3-\eps}}{\varphi(d)\pi_E(X;d,1)}
\sim c(E)\li(X)
$$
(as $X\ra +\infty$) would hold for any $\eps>0$. Indeed
from~(\ref{eq-dun-part-noncm}), it suffices to
estimate the sum over $X^{1/4}/(\log X)\leq d\leq X^{1/3-\eps}$.
This can be done using the Brun-Titchmarsh
inequality~(\ref{eq-bt-generale}) for $K_d$, since primes which are
totally split in $K(E[d])$ must also be so in $K_d$:
\begin{align*}
\sum_{X^{1/4}/(\log X)\leq d\leq X^{1/3-\eps}}
{\varphi(d)\pi_E(X;d,1)}&\ll_{\eps}
\frac{X}{\log X}\sum_d{\frac{\varphi(d)}{|PGL(2,\Zz/d\Zz)|}}\\
&\ll_{\eps} X^{3/4+\eps}.
\end{align*}
\par
One may ask similar questions with more general sets of conjugacy
classes replacing the identity element; this is left to the reader to
formulate, together with some potentially useful example for the
elliptic splitting problem.
\par
Besides the cyclotomic extension of $\Qq$, it seems few
Brun-Titchmarsh fields are known. We will see in
Section~\ref{ssec-cm-split} that the division fields of CM elliptic
curves provide further examples. But all those correspond to
(essentially) abelian Galois groups.

\begin{problem}
Find a Brun-Titchmarsh extension $K'/K$ corresponding to
an algebraic group $G/\Zz$ of dimension $>0$ with non-abelian
connected component.
\end{problem}

The known proofs of the classical Brun-Titchmarsh inequality and of
those for CM curves are based on sieve
methods: one can use almost any form of `additive' sieve
(see~\cite{halb-rich}) or a refined version of the large sieve
(see~\cite[\S 3 or \S 4]{bom-ls}). The latter may be generalized, to a
certain extent using techniques as in~\cite[Prop. 9]{km-density} to
handle Artin $L$-functions, but this requires to be useful that
all irreducible representations of the finite groups
$G(\Zz/d\Zz)$ be of degree $\leq \gamma'$ for some $\gamma'>0$
independent of $d$, which is equivalent to the connected component of
$G/\Zz$ being abelian. 
However, this fails to give useful information for the Brun-Titchmarsh
problem; this is because the required saving of the factor
$1/|G(\Zz/d\Zz)|$ comes, in the case of arithmetic progressions, from
summing over integers $n\equiv 1\mods{d}$ by writing $n=md+1$ and
summing over $m$. This underlying regularity is of course inexistent
in more complicated extensions.
\par
This suggests another problem: prove~(\ref{eq-bt-classic})
\emph{without} appealing to the regularity of arithmetic progressions.

\mysection{Elliptic twins}\label{sec-twins}

\subsection{Definition}

The first step in the direction of Problem~\ref{prob-outside}
introduces instead another interesting analytic problem.
Let $K=\Qq$ for simplicity.
Fix $X\geq 1$ and an integer $d$ such that $d^2>8X^{1/2}$. Let
$\{p_1,\ldots, p_k\}$ be the set of primes splitting completely
in $\Qq(E[d])$ (i.e. $d\mid \dun(E_{p_j})$) with
$p_j\leq X$, and assume they are indexed in increasing order, so that
$p_j< p_k$ if $j< k$.
\par
Consider $p=p_j$ and $q=p_{j+1}$ for some $j$. Since
$$
d^2\mid \dun(p)^2 \mid n_p(E)=p+1-a_p(E)\text{ and also }
d^2\mid \dun(q)^2 \mid q+1-a_q(E)
$$
we get by subtracting
$$
d^2\mid (q-p)+(a_p(E)-a_q(E)).
$$
Therefore, \emph{if} the right-hand side is non-zero, it follows that
$$
q\geq p+(a_p(E)-a_q(E))+d^2,
$$
but by the Riemann Hypothesis for $E_p$ and $E_q$, and the assumption
$d^2>8X^{1/2}>8q^{1/2}$, we have
$$
|a_p(E)-a_q(E)|\leq 2(\sqrt{p}+\sqrt{q})\leq \frac{d^2}{2},
$$
hence we get a gap between $p$ and $q$,
$$
q\geq p+\frac{d^2}{2},
$$
which is stronger than the ``trivial'' gap imposed by the congruence
$p\equiv q\equiv 1\mods{d}$.
\par
However, this is subject to the condition that
$$
(q-p)+(a_p(E)-a_q(E))\not=0
$$
which is equivalent with
$$
|E_p(\Fp_p)|\not=|E_q(\Fp_q)|.
$$
There is no reason this should not occur, and this prompts the
following general definition: 

\begin{definition}
Let $K$ be a number field and $E/K$ an elliptic curve. Two distinct
prime ideals $\ideal{p}$ and $\ideal{q}$ of $K$ are called \emph{elliptic
twins} for $E$ if 
$$
|E_{\ideal{p}}(\Fp_{\ideal{p}})|=|E_\ideal{q}(\Fp_{\ideal{q}})|
$$
i.e. $E$ has as many points reduced modulo $\ideal{p}$ and modulo
$\ideal{q}$. We say that $\ideal{p}$ has an $E$-twin, or simply a
twin.
\end{definition}

\begin{remark}
More generally, let $C/K$ be an algebraic curve (or even an arbitrary
algebraic variety) and fix $\mathcal{C}/\Oc_K[1/S]$ a model of $C$
defined over the integers of $K$ (minus a finite set $S$ of
primes). Two distinct prime ideals $\ideal{p}$ and $\ideal{q}$ of $K$
which are not in $S$ are called $C$-twins if 
$$
|C_{\ideal{p}}(\Fp_{\ideal{p}})|=|C_\ideal{q}(\Fp_{\ideal{q}})|
$$
We say that $\ideal{p}$ has a $C$-twin. 
Note that except for finitely many pairs, this is independent of the 
choice of the model $\mathcal{C}$, but for definiteness one may chose
one of the preferred models of $C$, or define twins for a variety defined
over an open subset of $\spec \Oc_K$.\footnote{ Especially since no
  variety is known to have infinitely many twin pairs...}
\par
N. Katz first suggested the following case, justifying the
\emph{rapprochement} with twin primes: instead of an
elliptic curve, consider the affine conic $C\,:\,x^2+y^2=1$ over
$\Qq$ (equivalently, to stay with algebraic groups, the restriction of
scalars from $\Zz[i]$ to $\Zz$ of the kernel of the norm map
$\mathbf{G}_{m/\Zz[i]}\ra \mathbf{G}_{m/\Zz}$). This ``is'' a model
over $\Zz$, and we have (remember $C$ is affine)
$$
|C(\Fp_p)|=
\begin{cases}
p+1&\text{ if $p\equiv 3\mods{4}$ ($p$ is inert in $\Qq(i)$)}\\
p-1&\text{ if $p\equiv 1\mods{4}$ ($p$ splits in $\Qq(i)$).}
\end{cases}
$$
Consequently, the condition
$|C_{p}|=|C_{q}|$ means either $p=q$, or
(1) $p\equiv 1\mods{4}$ and $p-2$ is
prime (it is then inert and
$|C_{p-2}|=(p-2)+1=p-1=|C_{p}|$), or (2) $p\equiv 3\mods{4}$
and $p+2$ is prime, which is (1) with $p$ and $p+2$
exchanged. Hence the $C$-twins are ``half'' the ordinary twin
primes, namely pairs $(p,p+2)$ with $p\equiv 3\mods{4}$.
\par
Note it doesn't seem to be possible to get the other half of all
twin primes\footnote{ 
 Numerical experiments (and standard conjectures) confirm that
those ``two-halves'' are equidistributed, in an obvious sense.
} in this manner: using a conic one would need a quadratic field $K$
with the property that $p$ is split in $K$ if and only if $p\equiv
3\mods{4}$
We ask:
\begin{question}
Is there an algebraic variety $\mathcal{X}/\Zz[1/2]$ with the property
that $p>2$ and $q>p$ are $\mathcal{X}$-twins if and only if $q=p+2$ ?
Is there one such that $p$ and $q$ are $\mathcal{X}$-twins if and only
if $q=p+2$ and $p\equiv 1\mods{4}$?
\end{question}
The author's guess is ``No''.
\end{remark}

The definition of elliptic twins certainly looks unnatural from a
geometric viewpoint: we compare the reduction of a curve modulo two
distinct primes. But in the absence of better ways of bounding the
number of \outside primes, and as analogues of the ordinary twin
primes, they are worth investigating.

\subsection{General facts}

We now introduce some more notation. Fix an elliptic curve $E/K$
defined over a number field. 
We define three arithmetic functions: 
\begin{align}
\label{eq-n}
n_{\ideal{p}}&=|E_{\ideal{p}}(\Fp_{\ideal{p}})|,\\
\label{eq-M}
\grandm(n) &= |\{ \ideal{p}\,|\, N\ideal{p}\leq X\,\text{ and }
n_{\ideal{p}}=n\}|,\\ 
\label{eq-M1}
\petitm(\ideal{p}) &= \grandm(n_{\ideal{p}}).
\end{align}
So $n_{\ideal{p}}$ and $\petitm(\ideal{p})$ are supported on primes of
$K$, and $\grandm(n)$ is defined for all $n\geq 1$.
\par
Of course $\ideal{p}$ has an $E$-twin 
if and only if $\petitm(\ideal{p})>1$. We'll say
that an $n\geq 1$ is a twin value if $\grandm(n)>1$, and call the
primes $\ideal{p}$ with $n_{\ideal{p}}=n$ the $E$-twins
associated to $n$.
\par
The main questions about elliptic twins concern the
behavior of those three functions. 
In particular:
\begin{question}\label{q-1}
What is the behavior of the function
\begin{equation}\label{eq-ejum}
\ejum(X)=|\{ n\leq X\,\mid\, n\text{ is a twin value }\}|
\end{equation}
counting the twin values up to $X$, or of
\begin{equation}\label{eq-ejumbis}
\ejumbis(X)=|\{\ideal{p}\,|\, N\ideal{p}\leq X\text{ and }\ideal{p}
\text{ has an $E$-twin}\}|.
\end{equation}
\end{question}
\begin{question}\label{q-2}
What is the behavior of the sum
\begin{equation}\label{eq-somp1}
\somp(X)=\sum_{N\ideal{p}\leq X}{\petitm(\ideal{p})}
\end{equation}
as $X\ra +\infty$ ?
\end{question}
\begin{question}\label{q-3}
More generally, for fixed $k\geq 0$, what is the behavior of the
moments of $\petitm(\ideal{p})$ and $\grandm(n)$
\begin{align}
\somg_k(X)&=
\sum_{n\leq X}{\grandm(n)^k}\\
\somp_k(X)&=\sum_{N\ideal{p}\leq X}{\petitm(\ideal{p})^k}.
\end{align}
\end{question}
\begin{question}\label{q-4}
Differently formulated: what can be said about $\grandm(n)$? How large
can it be compared to $n$, and how does it behave as $n\ra +\infty$?
\end{question}

Question~\ref{q-1} is the elliptic analogue of the classical
twin-prime problem. On the other hand, because the analogue of
the ``multiplicity'' $\grandm(n)$ is simply the constant $2$ for the
twin-prime problem, Questions~\ref{q-2},~\ref{q-3} and~\ref{q-4} do
not have a classical counterpart and are genuinely elliptic problems.
\par
Also of interest is the dependence on $E$ of all those quantities, in
particular the ``meta-question'' is: what global arithmetic invariants
of $E$ can be extracted from information about the functions
$\grandm(n)$ and $\petitm(n)$? (Recall that according to the Isogeny
Theorem, the curve $E/K$ is determined up to $K$-isogeny
by the function $\ideal{p}\mapsto n_{\ideal{p}}$). We will see that it
is likely that one can extract from the asymptotic of $\ejum(X)$
whether $E$ has CM or not.
Recall the notation $x^+$ and $x^-$ (\ref{eq-x-pm}).

\begin{lemma}\label{lm-simple}
Let $E/K$ be an elliptic curve over a number field. For any $n\geq 1$
we have
\begin{align}\label{eq-control-1}
n_{\ideal{p}}=n& \Rightarrow n^-\leq N\ideal{p}\leq n^+,\\
& \Rightarrow N\ideal{p}^-\leq n\leq N\ideal{p}^+.
\label{eq-control-2}
\end{align}
and
\begin{equation}\label{eq-triv-borne-M}
\grandm(n)\leq |\{\ideal{q}\,\mid\, \ideal{q} \text{ is prime
and } n^-\leq N\ideal{q}\leq n^+\}| \ll [K:\Qq]\frac{\sqrt{n}}{\log (n+1)},
\end{equation}
the implied constant being absolute.
\end{lemma}

\begin{proof}
The implications~(\ref{eq-control-1}) and~(\ref{eq-control-2}) are
just the Riemann Hypothesis~(\ref{eq-rh})
for $E_{\ideal{p}}$.
The bound on $\grandm(n)$ then follows trivially by definition; for
the last inequality, observe that if $\ideal{q}$ is prime and $n^-
\leq N\ideal{q}\leq n^+$, $N\ideal{q}$ is a prime power in that range,
of which the number is $\ll \sqrt{n}/(\log (n+1))$, with an absolute
implied constant. Each prime power $q^f$ can occur for at most $[K:\Qq]$
prime ideals since $\ideal{q}$ must be above $q$ in $K$
(compare~(\ref{eq-borne-ind-cm}) below).
\end{proof}

\begin{remark} The delicacy of the matter is indicated by the fact
that the size (about $\sqrt{n}$ ideals among $n$ with $N\ideal{a}\leq
n$) of this range is just such that even on the Generalized Riemann
Hypothesis it is not possible to ensure that it contains at least
one prime ideal $\ideal{p}$ for all $n$ large enough.
Indeed, on GRH we have
$$
\pi_K(X) = \li(X) + \oun(X^{1/2}(\log \Delta_K X)
$$
(where $\Delta_K$ is the absolute value of the discriminant of $K$;
the implied constant is absolute, see~\cite[]{serre-tch} for 
instance). This only implies 
$$
\pi_K(n^+)-\pi_K(n^-) \ll n^{1/2}(\log \Delta_K n)
$$
which is worse than the trivial bound obtained by counting all integral
ideals.
\end{remark}

The ``trivial'' bound~(\ref{eq-triv-borne-M}) is in a sense
best possible, because it is possible to find curves over a finite
prime field $\Zz/p\Zz$ with any value of $a_E$ satisfying $|a_E|\leq
2\sqrt{p}$. We state more formally this easy fact:

\begin{proposition}\label{pr-bad}
Let $n\geq 1$ be an integer. There exists an elliptic curve $E/\Qq$
with good reduction at all primes $p$ such that $n^-\leq p\leq n^+$,
and with $n_p=n$ for all such primes.
\end{proposition}

In contrast with the remark above, note that it is known
that for ``most'' integers $n$ the number of primes described is $\gg
\sqrt{n}/(\log n)$ (see e.g.~\cite{harman}, where this is shown to
hold for $n<p<n+n^{\delta}$ for any $\delta>1/10$; the case $\delta=1/2$
is much easier).

\begin{proof}
For $p$ with $n^-\leq p\leq n^+$, let $b_p=p+1-n$, so by
construction we have
$|b_p|\leq 2\sqrt{p}$. By work of Deuring~\cite{deuring}
(Honda-Tate theory for elliptic curves, see
Theorem~\ref{th-honda-tate} below), there
exists an elliptic curve $E/\Fp_p$ with $a_E(p)=b_p$.
Consider a Weierstrass equation
$$
E_p/\Fp_p\,:\,y^2+a_1(p)xy+a_3(p)y=x^3+a_2(p)x^2+a_4(p)x+a_6(p)
$$
for such a curve. By the Chinese Remainder
Theorem we can find $a_i\in \Zz$, $i=1,2,3,4,6$, reducing to $a_i(p)$
modulo $p$ for all $p$ with $n^-\leq p\leq n^+$. Then the curve
$$
y^2+a_1 xy+a_3 y=x^3+a_2 x^2+a_4 x+a_6
$$
is an elliptic curve (since it reduces to a non-singular curve modulo
those primes), and it has $n_p=n$ for all the primes $p$ in question.
\end{proof}

Of course, having constructed one $n$ with $\grandm(n)\gg
\sqrt{n}/(\log n)$ does not tell anything about the asymptotic growth of
$\grandm(n)$ as $n\ra +\infty$. The following trivial lemma shows (in
particular) that on average $\grandm(n)$ is much smaller.

\begin{lemma}\label{lm-same-thing}
Let $E/K$ be an elliptic curve over a number field, $k\geq 1$ an
integer. We have
$$
\somg_k(X)=
\sum_{\stacksum{N\ideal{p}\leq X^+}{n_{\ideal{p}}\leq X}}
{\petitm(\ideal{p})^{k-1}}.
$$
In particular
\begin{equation}\label{eq-approx}
\sum_{n\leq X}{\grandm(n)} = \pi_K(X) + \oun_{K}(\sqrt{X}),
\end{equation}
where $\pi_K(X)$ is the number of prime ideals of $K$ with
$N\ideal{p}\leq X$.
\end{lemma}

Note that $n_{\ideal{p}}\leq X$ \emph{implies} $N\ideal{p}\leq X^+$,
but this condition is included in the summation to recall how
the size of $N\ideal{p}$ is controlled.

\begin{proof}
We have
\begin{align*}
\sum_{n\leq X}{\grandm(n)^k}&=
\sum_{n\leq X}{\grandm(n)^{k-1}\Bigl(\sum_{n_{\ideal{p}}=n}{1}
\Bigr)}
=\sum_{n_{\ideal{p}}\leq X}{\grandm(n_{\ideal{p}})^{k-1}}\\
&=\sum_{\stacksum{N\ideal{p}\leq X^+}{n_{\ideal{p}}\leq X}}{
\petitm(\ideal{p})^{k-1}}.
\end{align*}
Then~(\ref{eq-approx}) follows by taking $k=1$ and noting that
$$
\Bigl|\sum_{\stacksum{N\ideal{p}\leq X^+}{n_{\ideal{p}}\leq X}}
{1}-\sum_{N\ideal{p}\leq X}{1}\Bigr|
\leq \sum_{X^-\leq N\ideal{p}\leq X^+}{1}
\ll_K \sqrt{X}.
$$
\end{proof}

\begin{question}\label{q-5}
Is it true that
\begin{equation}\label{eq-q-5}
\grandm(n)=\oun_{E,\eps}(n^{\eps})
\end{equation}
for all $\eps>0$?
\end{question}

We will see in Section~\ref{sec-cm-resultats} that this is true for CM
curves (and we will give a 
more precise result). Heuristic and numerical evidence point
to even stronger results, but note that because of
Proposition~\ref{pr-bad}, any progress requires using global
properties of the elliptic curve.
\par
If~(\ref{eq-q-5}) holds it follows that we have
\begin{equation}\label{eq-g-p-equiv}
\somg_k(X)
=
\somp_{k-1}(X)
+\oun_{E,\eps,k}(X^{1/2+\eps}),\text{ for any } \eps>0.
\end{equation}

Finally we remark that the two functions $\ejum(X)$ and $\ejumbis(X)$
are somewhat different, since $\ejumbis(X)$ counts the twins with
multiplicity. For this reason (see Section~\ref{sec-cm-resultats}), it is a
little bit easier to deal with.

\subsection{Heuristic}\label{ssec-heuristic}

Here we consider an elliptic curve $E/\Qq$ which doesn't have CM, and we
make some rough heuristics concerning elliptic twins.
It should be possible to give somewhat more convincing
arguments and more precise predictions using a probability model such
as that used by Lang-Trotter~\cite{lang-trotter}.
\par
For a prime number $p$, there are about $4\sqrt{p}$ possible values of
$a_p$, and according to the Sato-Tate conjecture, they should be
such that the angle $\theta_p\in [0,\pi]$ satisfying
$$
a_p=2\sqrt{p}\cos\theta_p
$$
is equidistributed with respect to the measure
$d\mu=\frac{2}{\pi}\sin^2 \theta d\theta$.
\par
Compared to the uniform measure, this measure is concentrated around
$0$, which should tend to limit the possibility of $E$-twins
occurring, since a twin $q$ must have $a_q=n_p-q-1$, so $q$ getting
relatively large sends $a_q$ towards the extreme, less probable, range
of possible values. In particular, for heuristic purpose, assuming
$a_p$ to be uniformly distributed should bias the result towards
\emph{more} twins.
\par
In a uniform situation, each possible prime $q$, $p^-\leq q\leq p^+$,
has probability about $1/\sqrt{q}$ of being a twin of $p$. Since $q$
must be prime, this makes a probability about
$$
\approx \frac{1}{\sqrt{p}} \times \frac{\sqrt{p}}{(\log p)} \approx
\frac{1}{\log p}
$$
for $p$ to have at least one twin. This is comparable to the situation
with classical twin primes $p$, the probability of $p+2$ being prime
being about $1/(\log p)$. In particular we can ask
\begin{question}
Let $E/\Qq$ be an elliptic curve over $\Qq$ without CM. Prove or
disprove that 
\begin{equation}\label{eq-ejum-conj}
\ejum(X)\sim c \frac{X}{(\log X)^2}
\end{equation}
for some $c>0$ as $X\ra +\infty$.
\end{question}
It is conjectured that the number $\pi_2(X)$ of twin
primes $\leq X$ satisfies
$$
\pi_2(X)\sim c_2 \frac{X}{(\log X)^2}\text{ with }
c_2=2\prod_{p\geq 3}{\Bigl(1-\frac{1}{(p-1)^2}\Bigr)}=1.3203\ldots
$$
In Section~\ref{sec-cm-resultats}, we'll see it seems more plausible that for
$E$ with CM, we have
$$
\ejumbis(X)\sim c \frac{X}{\log X},\quad
\ejum(X)\sim c \frac{X}{(\log X)^{3/2}}.
$$
Concerning the multiplicity question, the same vague heuristic
suggests that the probability that $p$ has $k$ twins is
about $1/(\log p)^k$, and this would seem to imply that the maximal
multiplicity is
\begin{equation}\label{eq-conj-multimax}
\petitm(p)=k\approx \frac{\log p}{\log\log p}.
\end{equation}
Again, in the CM case, Section~\ref{sec-cm-resultats} suggests that
$\petitm(p)$ can be much larger, almost as large as a divisor-like
function.
\par
For numerical experiments, see Section~\ref{sec-numerics} below.



\mysection{Curves with complex multiplication}\label{sec-cm-resultats}


The analytic problems we have raised can be analyzed much further for
CM curves. For elliptic twins, this will reveal some differences (so
that, for instance, the 
behavior of $\grandm(n)$, $\petitm(\ideal{p})$ should distinguish
between CM and non-CM curves) while highlighting in a different way
the connexion with the classical twin primes. We will prove upper
bounds for the moments of $\grandm(n)$. Those upper bounds are such
that general expectations about primes represented by polynomials lead
to believe that they are of the correct order of magnitude.

\subsection{Preliminaries}

We recall the basic facts of complex multiplication theory that
describe the reductions of a CM curve and their Frobenius
endomorphisms. The theory is basically due to Deuring; see for
instance~\cite[II]{silv-2} for a modern treatment.
\par
Let $E/H$ be an elliptic curve over a number field $H$ with CM by
an order $\Oc$ in the ring of integers $\Oc_K$ of a quadratic
imaginary field $K$. \emph{For simplicity, we will assume in this
  section that $K\subset H$, i.e. the defining field contains the CM
  field.} This excludes in particular the important case $H=\Qq$, but
the principle still applies in the general case, and we will
extend the results for one curve over $\Qq$ in
Section~\ref{ssec-num-twins}, so that a complete treatment could be easily
obtained (recall that in any case the composite field $HK$ is at most
a quadratic extension of $H$, so the case $H=\Qq$ is really
``complementary'' to the case $K\subset H$). For a given imaginary
quadratic order $\Oc$, it is known (\cite[]{silv-2}) that all elliptic
curves with CM by $\Oc$ can be defined over the ring-class field
associated to $\Oc$ (e.g., if $\Oc=\Oc_K$, over the Hilbert
class-field of $K$).
\par
The following notation will be used: for an imaginary quadratic field
$K/\Qq$, we let $\chi=\chi_K$ denote the Kronecker symbol for $K$,
i.e. the primitive quadratic Dirichlet character associated to $K$ by
class-field theory, and let $r(n)$ or $r_K(n)$ denote the arithmetic
function 
$$
r(n)=r_K(n)=|\{ \ideal{a}\subset \Oc_K\,\mid\, N\ideal{a}=n\}|
$$
so that the Dedekind zeta function of $K$ is given by
\begin{align*}
\zeta_K(s)&=\sum_{\ideal{a}}{(N\ideal{a})^{-s}}=
\sum_{n\geq 1}{r_K(n)n^{-s}}\\
&=\prod_{\ideal{p}}{(1-(N\ideal{p})^{-s})^{-1}}
=\zeta(s)L(s,\chi_K).
\end{align*}
In particular,
\begin{equation}\label{eq-simple-r}
r_K(n)=\sum_{d\mid n}{\chi_K(d)},\quad r_K(n)\leq d(n),
\end{equation}
where $d(n)$ is the ``number of divisors'' function.
It will be convenient to fix once and for all a basis
$(1,\omega)$ of $\Oc_K$ as a $\Zz$-module.
\par
The following result is that part of the Main Theorem of Complex
Multiplication that will be needed (see~\cite[]{silv-2}):
\begin{theorem}\label{th-cm}
With the above notation, there exists a map $\ideal{p}\mapsto
\psi(\ideal{p})\in\Oc$, from the set of prime ideals of $H$ where $E$ is
unramified to $\Oc$, with the property that $\psi(\ideal{p})$ is the
Frobenius automorphism for $E_{\ideal{p}}/\Fp_{\ideal{p}}$.
\end{theorem}

In fact, properly normalized, this map extends to the
Gr\"ossencharakter of $E$ (\cite[]{silv-2},~\cite[]{rubin}), but we do
not need this deeper fact. 
\par
We denote by $\frobs(E)$ the image of $\psi$, i.e. the set of all
Frobenius endomorphisms of $E$ at primes of $K$.
\par
By the properties of the Frobenius automorphism, if $\ideal{p}$ is an
unramified prime ideal of $H$, we have
\begin{equation}\label{eq-cm-1}
N^H_{\Qq} \ideal{p}=|\Fp_{\ideal{p}}|=N_{\Qq}^K(\psi(\ideal{p})),
\end{equation}
and
\begin{equation}
  \label{eq-cm-2}
n_{\ideal{p}}=N^K_{\Qq}(\psi(\ideal{p})-1).
\end{equation}

We will reduce the problems about prime ideals in $H$ to those of $K$
using the following simple lemma:

\begin{lemma}\label{lm-borne-cm}
With the same notation as before, for any prime ideal $\ideal{p}$ in
$K$, $\psi(\ideal{p})$ is divisible by a single prime $p$, and 
for any $z\in \Oc_K$ with this property
$$
|\{ \ideal{p} \,\mid\, \psi(\ideal{p})=z\}|
\leq [H:\Qq]
$$
\end{lemma}

\begin{proof}
Equation~(\ref{eq-cm-1}) proves the first statement. Then for any
$\ideal{p}$ in $H$ with $\psi(\ideal{p})=z$, the prime $p$ below
$\ideal{p}$ in $\Qq$ is independent of $\ideal{p}$: it is the unique
$p$ such that $N^K_{\Qq}z=p^{\nu}$ for some $\nu\geq 1$. Hence the
number of $\ideal{p}$ is $\leq [H:\Qq]$.
\end{proof}


\subsection{Elliptic twins}

We apply now the theory of complex multiplication to elliptic
twins. We keep the same notation and convention.
First we can answer Question~\ref{q-5} for a CM curve.


\begin{proposition}\label{pr-upper-ind-cm}
Let $E/H$ be a CM curve. We have for $n\geq 1$
\begin{equation}\label{eq-borne-ind-cm}
\grandm(n)\leq  [H:\Qq]r_K(n),
\end{equation}
and in particular for any $n\geq 1$ and any $\eps>0$,
\begin{equation}\label{eq-rep-q-5-cm}
\grandm(n)=\oun_{E,\eps}(n^{\eps}),
\end{equation}
the implied constant depending only on $E$ and $\eps$.
\end{proposition}

\begin{proof}
By~(\ref{eq-cm-2}), for any $\ideal{p}$ with
$n_{\ideal{p}}=n$, the integer
$z_{\ideal{p}}=\psi(\ideal{p})-1\in \Oc\subset \Oc_K$ is a solution to
the norm equation $N_{\Qq}^K z=n$ in $K$.
Moreover, if $z$ is any solution of this equation, all prime ideals
$\ideal{p}$ with $z_{\ideal{p}}=z$ satisfy $\psi(\ideal{p})=z+1$. 
Thus by Lemma~\ref{lm-borne-cm}, for each $z$ there are at most
$[H:\Qq]$ prime ideals $\ideal{p}$ with $z_{\ideal{p}}=z$,
hence~(\ref{eq-borne-ind-cm}) follows. 
\par
Now~(\ref{eq-rep-q-5-cm}) is immediate since $r(n)\ll_{\eps}
n^{\eps}$ (for instance, use~(\ref{eq-simple-r})).
\end{proof}

Our main result is the following theorem.

\begin{theorem}\label{th-upper-cm}
Let $E/H$ be a CM elliptic curve as above. For any $\eps>0$, we have
\begin{align}\label{eq-upper-cm-g}
\somg_k(X) \ll X(\log X)^{\beta(k-1)+\eps}
&\quad\text{ for $k\geq 1$}\\
\somp_k(X) \ll X(\log X)^{\beta(k)+\eps}
&\quad\text{ for $k\geq 0$},
\label{eq-upper-cm-p}
\end{align}
for $X\geq 2$, where
\begin{equation}
  \label{eq-beta-k}
\beta(k)=2^{k}-k-2.
\end{equation}
The implied constants depend on $k$, $K$, $H$ and $\eps$.
\end{theorem}

\begin{remark}\label{rem-annoying}
One can probably put $\eps=0$; indeed, this is the case for
$\somp_k(X)$ for $k=0$, $k=1$, and for all $k$ the proof yields a
stronger result with $(\log X)^{\eps}$ replaced by a power of
$\log\log X$; since I believe this is mistaken anyway (see the proof
of Proposition~\ref{pr-somme}), I prefer not to put this stronger
statement.
\end{remark}

For example,
$$
\sum_{N\ideal{p}\leq X}{\petitm(\ideal{p})}\ll_K \frac{X}{\log X},
\quad 
\sum_{N\ideal{p}\leq X}{\petitm(\ideal{p})^2}\ll_K
X.
$$
Moreover, we'll see in the course of proving the theorem that standard
conjecture about primes represented by polynomials imply that the
estimates~(\ref{eq-upper-cm-p}) and~(\ref{eq-upper-cm-g}) are of the
correct order of magnitude. For $\somp_0(X)$, this is just the Prime
Ideal Theorem in $K$ (and doesn't give any information about elliptic
twins). 
\par
Before starting the proof, we remark that by
Proposition~\ref{pr-upper-ind-cm} and~(\ref{eq-g-p-equiv}), the
bounds~(\ref{eq-upper-cm-g}) and~(\ref{eq-upper-cm-p}) are
equivalent. We will work with $\somp_k(X)$ for $k\geq 1$, the case
$k=0$ being obvious.


\subsection{Reduction to twin-prime-like counting}
\label{ssec-twin-prime}

The strategy of the proof is to reduce to some counting of (principal)
prime ideals in the ring $\Oc_K$, and to use~(\ref{eq-cm-2}) to put
the counting into the shape of ``parallel'' twin-prime-like
equations, for which upper bounds of the (conjectural) correct order
of magnitude can be efficiently and uniformly obtained by a sieve
method. In this case, we'll use Huxley's version of the large sieve in
number fields~\cite{huxley}.
\par
A \emph{prime element} in $\Oc_K$ is an integer $z$ such that $(z)$ is
a prime ideal.
We first reduce to those $\ideal{p}$ such that $\psi(\ideal{p})$ is
a prime element.
\par
\begin{lemma}\label{lm-reduc-fp-one}
Let $E/H$ be as above. We have for any $k\geq 0$ and any $\eps>0$
\begin{equation}
\label{eq-reduc-fp-one}
\sum_{\stacksum{N\ideal{p}\leq X}{f_{\ideal{p}}\geq 2}}
{\petitm(\ideal{p})^k}\ll_{\eps,K,k} X^{1/2+\eps}
\end{equation}
the implied constant depending only on $\eps$, $K$ and $k$. In the
sum, $f_{\ideal{p}}$ is the residual degree of $\ideal{p}$.
\end{lemma}
\begin{proof}
By~(\ref{eq-rep-q-5-cm}), we have
\begin{align*}
\sum_{\stacksum{N\ideal{p}\leq X}{f_{\ideal{p}}\geq 2}}
{\petitm(\ideal{p})^k} & \ll_{\eps,k} X^{\eps}
\sum_{\stacksum{N\ideal{p}\leq X}{f_{\ideal{p}}\geq 2}}
{1}\\
&\ll_{\eps,k} X^{\eps}
\sum_{\stacksum{p^k\leq X}{k\geq 2}}{r_K(p^k)}\\
&\ll_{\eps,k,K} X^{1/2+\eps}.
\end{align*}
\end{proof}
\par
Henceforth we only consider prime ideals $\ideal{p}$ of $H$ which are
of degree $1$. In particular, by~(\ref{eq-cm-1}), $\psi(\ideal{p})$
is then a prime element of $\Oc_K$.
\par
\medskip
Next we deal with the parameterization of elliptic twins. Recall that
an integer $z\in 
\Oc_K$ is \emph{primitive} if it is not divisible by any $d\in \Zz$,
$d\not= \pm 1$; in terms of the basis $(1,\omega)$ of $\Oc_K$, if
$z=a+b\omega$, this 
means that $a$ and $b$ are coprime. We let $\prim$ denote the set of
primitive elements in $\Oc_K$ \emph{modulo $\pm 1$}. Note that any
non-zero $z\in\Oc_K$ can be written $z=dv$ for some $d\in \Zz$ and
some $v\in\prim$: if $z=a+b\omega$, $d=(a,b)$, $v=z/d$. The pair
$(d,v)$ is unique, up to simultaneous sign-change.
\par
The norm $Nu$ of an element $u\in\prim$ is well-defined. So is the
complex-conjugation (i.e. the action of the Galois group of $K$). 
In addition, for a $k$-tuple
$\underline{u}=(u_1,\ldots,u_k)\in\prim^k$ we define the
\emph{discriminant} $\disc(\underline{u})$ to be
\begin{equation}\label{eq-disc-def}
\disc(\underline{u})=
\prod_{1\leq i<j\leq k}{(u_i \bar{u}_j-\bar{u}_i u_j)}.
\end{equation}
This is well-defined up to sign so it can be thought of as an integral
ideal in $K$. Note that, by primitivity, $\disc(\underline{u})=0$ if
and only if there exist $i\not= j$ such that $u_i=u_j$ (in $\prim$)
(see the proof of the next lemma).

\begin{lemma}\label{lm-bij}
Let $\nrmone$ be the set of $z\in K$ of norm $1$. There exists a
bijection 
$$
\bij\,:\,\prim\lra \nrmone
$$
given by $u\mapsto \bar{u}/u$ for $u\in \prim$.
\end{lemma}

\begin{proof}
Clearly $\bij$ maps $\prim$ into $\nrmone$. Moreover, $\bij$ is
injective: if $\bij(v)=\bij(w)$ with $v$, $w\in\prim$, we get $v/w\in
\Qq$ (because it is Galois-invariant), so we have $av=bw$ for some
$a$, $b\in \Qq$, $(a,b)=1$. Because $v$ and $w$ are primitive, this
implies that $|a|=|b|=1$, so $v=\pm w$.
\par
It remains to prove surjectivity. 
This amounts essentially to finding all pythagorean triples (when
$K=\Qq(i)$), but instead of doing it by hand, we can appeal to
Hilbert's Theorem 
90 for $K/\Qq$ (see e.g.~\cite[VIII-6]{lang}): for $z\in K^{\times}$,
$Nz=1$ is equivalent with $z=\bar{w}/w$ for some $w\in
K^{\times}$. Writing $w=vd/e$ for some $v\in \prim$ and $d$, $e\in
\Zz$, we have $z=\bar{v}/v=\bij(v)$.
\end{proof}

Note that one can write the discriminant $\disc(\underline{u})$ as a
Vandermonde determinant
$$
\disc(\underline{u})=(u_1\cdots u_k)^{k-1}
\prod_{1\leq i<j\leq k}{(\bij(u_i)-\bij(u_j))}
=(u_1\cdots u_k)^{k-1}\left|\bij(u_i)^{j-1}\right|_{i,j}.
$$

\begin{lemma}
Let $K$ be an imaginary quadratic field. For integers $w$, $z\in
\Oc_K$, we have
\begin{equation}\label{eq-nrm-1}
N(w-1)=N(z-1)
\end{equation}
if and only if there exists an $u\in\prim$ such that $w=f_u(z)$,
where $f_u$ is the linear form
\begin{equation}\label{eq-fu}
f_u(z)=\bij(u)(z-1)+1=\frac{\bar{u}}{u}(z-1)+1.
\end{equation}
Such an element $u\in\prim$ is unique.
\end{lemma}

\begin{proof}
This is an immediate consequence of the previous
lemma:~(\ref{eq-nrm-1}) holds if and only if
$N((w-1)/(z-1))=1$, therefore if and only if there exists a 
$u\in\prim$ (which is unique) with 
$$
\frac{w-1}{z-1}=\bij(u)=\frac{\bar{u}}{u},
$$
i.e. $w=\bij(u)(z-1)+1=f_u(z)$.
\end{proof}

Note that in this lemma we have $w=z$ if and only if $u=1$ and
$w=\bar{z}$ if and only if $u=z$.
\par
By~(\ref{eq-cm-2}), it follows that if $n_{\ideal{p}}=n_{\ideal{q}}$,
there exists $u\in\prim$ such that
$\psi(\ideal{p})=f_u(\psi(\ideal{q}))$. For a given $u$, since
$\psi(\ideal{p})$ is a prime element, this is similar to the
classical twin-prime problem: the question is to find prime elements
$\pi\in\Oc_K$ such that $f_u(\pi)$ 
is also prime (note that $f_u$ can not be properly defined for prime
\emph{ideals}). 
\par
There are infinitely many $u\in \prim$, but there is a
(congruence) condition for $f_u(z)$ to be an integer when $z\in \Oc_K$,
and this will restrict the values of $u$ occurring in a sum like
$\somp_k(X)$.

\begin{lemma}\label{lm-congruence-cond}
Let $\underline{u}=(u_1,\ldots,u_k)\in \prim^k$. For $z\in \Oc_K$, we
have 
$$
f_{u_i}(z)\in \Oc_K\text{ for all } i,\ 1\leq i\leq k,
$$
if and only if $z\equiv 1\mods{[\underline{u}]}$, where 
$[\underline{u}]$ is the (ideal) l.c.m of the elements $u_1$,\ldots,
$u_k$. 
\end{lemma}

\begin{proof}
It suffices to treat the case $k=1$, by definition of the l.c.m. Since
$$
f_u(z)=\frac{\bar{u}}{u}(z-1)+1
$$
we have $f_u(z)\in \Oc_K$ if and only if $\bij(u)(z-1)\in
\Oc_K$. Since $u$ is primitive, $u$ and $\bar{u}$ are coprime, so this
is equivalent with $z-1\in (u)$, i.e. $z\equiv 1\mods{u}$.
\end{proof}


In other words, the ``twin-prime problem'' for $f_u$ concerns only
prime elements $\pi\in \Oc_K$ with $\pi\equiv 1\mods{u}$.

\begin{corollary}\label{cor-mp}
Let $\ideal{p}$ be a prime ideal of $H$ with
$N\ideal{p}\leq X$. We have
\begin{equation}\label{eq-mp}
\petitm(\ideal{p})=
\sum_{\stacksum{Nu\leq X^+}{f_u(\psi(\ideal{p}))\in
\frobs(E)}}{1}.
\end{equation}
\end{corollary}

\begin{proof}
By the above we get directly
\begin{equation}\label{eq-pmp1}
\petitm(\ideal{p})=
\sum_{\stacksum{u\in\prim}{f_u(\psi(\ideal{p}))\in
\frobs(E)}}{1}.
\end{equation}
Let $z\in\Oc_K$ be an integer with $Nz\leq X$ and $f_u(z)\in\Oc_K$. By
Lemma~\ref{lm-congruence-cond}, we can write
$$
z=uv+1\text{ for some } v\in\Oc_K,
$$
which implies $Nu\leq N(z-1)\leq X^+$. Hence the result.
\end{proof}

\begin{corollary}\label{cor-mp1}
Let $\ideal{p}$ be a prime ideal of degree $1$ of $H$ with no twin of
degree $\geq 2$. We have
\begin{equation}\label{eq-mp1}
\petitm(\ideal{p})\leq 
2[H:\Qq]\sum_{\stacksum{Nu\leq \sqrt{X}+1}{f_u(\psi(\ideal{p}))
\text{ is prime}}}
{1}.
\end{equation}
\end{corollary}

\begin{proof}
As in the proof of Proposition~\ref{pr-upper-ind-cm}, to each
prime element $\pi$ of $\Oc_K$, there correspond at most $[H:\Qq]$
prime ideals $\ideal{p}$ of $H$ with $\psi(\ideal{p})=\pi$. Hence the
previous corollary implies
\begin{equation}\label{eq-mp2}
\petitm(\ideal{p})\leq [H:\Qq] 
\sum_{\stacksum{Nu\leq X^+}{f_u(\psi(\ideal{p}))\text{ is prime}}}{1}.
\end{equation}
Write $\pi=\psi(\ideal{p})$ for simplicity. For $z\in\Oc_K$ such that
$Nz\leq X$ and $f_u(z)\in\Oc_K$ we have by
Lemma~\ref{lm-congruence-cond}
\begin{equation}\label{eq-zuv}
z=uv+1\text{ for some } v\in\Oc_K,
\end{equation}
and
$$
f_u(z)=\bar{u}v+1.
$$
We can use the classical trick of Dirichlet of switching divisors:
remark that taking $v$ instead of $u$ in~(\ref{eq-zuv}) leads to
$$
f_v(z)=\bar{v}u+1=\overline{f_u(z)}.
$$
In particular, if $f_u(z)$ is prime, so is $f_v(z)$, hence both $u$
and $v$ occur together in~(\ref{eq-mp2}). Since one of them has norm
$\leq \sqrt{X^+}=\sqrt{X}+1$, the corollary follows.
\end{proof}

We now rewrite the sum $\somp_k(X)$.
\begin{lemma}
Let $k\geq 1$. We have
$$
\somp_k(X)=\sum_{\stacksum{\underline{u}\in\prim^k}
{Nu_i\leq X^+}}{\somp_{\underline{u}}(X)}
+\oun_{\eps,k,K}(X^{1/2+\eps})
$$
for any $\eps>0$, where
$$
\somp_{\underline{u}}(X)=|\{ \ideal{p}\text{ degree $1$ in } H
\,\mid\, N\ideal{p}\leq X, 
f_{u_i}(\psi(\ideal{p}))\in \frobs(E)
\text{ for } 1\leq i\leq k\}|
$$
for $\underline{u}=(u_1,\ldots, u_k)$.
\end{lemma}

\begin{proof}
By Corollary~\ref{cor-mp},
$$
\petitm(\ideal{p})=
\sum_{\stacksum{u\in\prim}
{f_u(\psi(\ideal{p}))\in\frobs(E)}}{1}.
$$
By Lemma~\ref{lm-reduc-fp-one} we can reduce to prime ideals of degree
$1$,
$$
\somp_k(X)=\sum_{\stacksum{N\ideal{p}\leq X}{f_{\ideal{p}}=1}}
{\petitm(\ideal{p})^k}+\oun(X^{1/2+\eps}).
$$
Expanding the $k$-th power and changing the order of
summation, the result follows.
\end{proof}

\begin{corollary}\label{lm-reduc-final}
Let $k\geq 1$. We have
$$
\somp_k(X)\ll \sum_{\stacksum{\underline{u}\in\prim^k}
{Nu_i\leq \sqrt{X}+1}}{\somp^+_{(1,\underline{u})}(X)}
+\oun_{\eps,k,K}(X^{1/2+\eps})
$$
for any $\eps>0$, 
where $(1,\underline{u})$ is a $(k+1)$-tuple, and for any $k$-tuple
$\underline{v}=(v_1,\ldots, v_k)$ we let
\begin{equation}\label{eq-reduc-final}
\somp^+_{\underline{v}}(X)=
|\{ z\in \Oc_K
\,\mid\, Nz\leq X \text{ and }
f_{v_i}(z)\text{ is prime for }
1\leq i\leq k\}|.
\end{equation}
The implied constant depends on $\eps$, $k$ and $H$.
\end{corollary}

\begin{proof}
Instead of Corollary~\ref{cor-mp}, we use Corollary~\ref{cor-mp1};
note that it may well
happen that $\ideal{p}$ is of degree $1$ but has a twin of 
degree $\geq 2$ and such twins are not counted
in~(\ref{eq-reduc-final}), but the contribution of such twins is
trivially $\ll X^{1/2+\eps}$, 
by the same argument used in Lemma~\ref{lm-reduc-fp-one}.
\end{proof}

Note that $[\underline{v}]=[\underline{u}]$ if
$\underline{v}=(1,\underline{u})$. 
\par
\medskip
Theorem~\ref{th-upper-cm} is a consequence of the following two
propositions:

\begin{proposition}\label{pr-upper-indiv}
Let $K/\Qq$ be an imaginary quadratic field, $k\geq 1$ an integer and
let $\underline{u}\in \prim^k$ with $N[\underline{u}]\leq X$. 
Assume that $u_i\not= u_j$ for $i\not= j$ (in $\prim$).
Then we have 
$$
\somp^+_{\underline{u}}(X)\ll_K
\frac{\phi_k(\underline{u})}{N[\underline{u}]}
\frac{X}{(\log (X/N[\underline{u}]))^{k}}
$$
for $X\geq 2$, where
\begin{equation}\label{eq-petit-phi}
\phi_k(\underline{u})=
\prod_{\ideal{p}\mid N[\underline{u}]\disc(\underline{u})}
{\Bigl(1+\frac{k}{N\ideal{p}}\Bigr)}.
\end{equation}
The implied constant depends only on $K$ and $k$.
\end{proposition}

\begin{proposition}\label{pr-somme}
Let $K/\Qq$ be an imaginary quadratic field, 
$k\geq 0$ an integer. For any $\eps>0$ we have
$$
\sum_{\stacksum{\underline{u}\in\prim^k}
{Nu_i\leq X}}{\frac{\phi_k(1,\underline{u})}{N[\underline{u}]}}
\ll_K (\log X)^{\gamma(k)+\eps}
$$
for $X\geq 2$, where $\gamma(k)=2^k-1$. The implied constant depends
only on $K$, $k$ and $\eps$. For $k=0$, we put, by convention,
$\prim^k=\{1\}$.
\end{proposition}

These will be proved in Section~\ref{ssec-twin-one}
and~\ref{ssec-twin-somme} respectively.
\par
\medskip
To finish the proof of Theorem~\ref{th-upper-cm}, let
$$
\somp_k^+(X)=\sum_{\stacksum{\underline{u}\in\prim^k}
{Nu_i\leq \sqrt{X}+1}}{\somp^+_{(1,\underline{u})}(X)}
$$
split the sum
$\somp_k^+(X)$ into $k$ subsums $\somp_{k,j}^+(X)$, $0\leq j\leq k$,
where $\somp_{k,j}^+(X)$ is the sum of the $\somp^+_{(1,\underline{u})}(X)$
for those $\underline{u}\in\prim^k$ where there are $j+1$ 
values among the components of $(1,\underline{u})$ 
i.e. the set $\{ 1, u_i \}$ has $j+1$ elements.
\par
By Lemma~\ref{lm-reduc-final} and Proposition~\ref{pr-upper-indiv}
(applied to the corresponding tuples $(1,\underline{i})$) for
$(j+1)$-tuples, we have 
\begin{align*}
\somp_{k,j}^+(X)&\ll 
\sum_{\underline{u}}{\somp_{(1,\underline{u})}^+(X)} + X^{1/2+\eps}\\
&\ll \sum_{\underline{u}}
{\frac{\phi_{j+1}(1,\underline{u})}{N[\underline{u}]}
\frac{X}{(\log X)^{j+1}}} + X^{1/2+\eps}\\
&\ll \frac{X}{(\log X)^{j+1}}
\sum_{\stacksum{\underline{v}\in\prim^j}
{Nv_i\leq X}}{\frac{\phi_{j+1}(1,\underline{v})}{N[\underline{v}]}}
\\
&\ll X(\log X)^{\gamma(j)-j-1+\eps}\text{ for any $\eps>0$ by
Proposition~\ref{pr-somme}.} 
\end{align*}
In the next-to-last inequality, we used the fact that if the set $\{
1, u_i\}$ has $j+1$ elements, $[\underline{u}]=[\underline{v}]$ where $v$
is any $j$-tuple whose components are the $j$ elements of $\{ u_i \}$,
and applied Proposition~\ref{pr-somme} for $j$ (there is a
multiplicity for each $\underline{v}$, but it is a combinatorial
function of $j$ and $k$ only).
\par
Summing over $j\leq k$, the theorem follows, since $j\mapsto
\gamma(j)-j-1=2^j-j-2$ is increasing for $j\geq 0$ ($0\mapsto -1$,
$1\mapsto -1$, $2\mapsto 0$, $3\mapsto 3\ldots$). The implied constant
depends on  $k$, $K$, and $H$.
\par
\begin{remark}
We conclude by justifying the assertion that Theorem~\ref{th-upper-cm}
should provide the correct order of magnitude for $\somg_k(X)$ and
$\somp_k(X)$ as $X\ra +\infty$ (up to the $(\log X)^{\eps}$ factor,
see Remark~\ref{rem-annoying}). First, for
$\somp^+_{\underline{u}}(X)$, we are
counting integers $z$ congruent to $1$ modulo $[\underline{u}]$ such that
the $k+1$ linear forms $f_{u_i}(z)$ take simultaneously prime
values. For any $u\in\prim$, we have $f_u(1)=1$, hence there
is no non-trivial common divisor to the values $f_u(z)$ for
$z\in\Oc_K$, $z\equiv 1\mods{u}$. Also, if no two $u_i$ coincide in
$\prim$, the condition that $f_{u_i}(z)$ be prime are
``independent''. Thus the usual heuristic predict that there should be
infinitely many $z\equiv 1\mods{\underline{u}}$ for which the
$f_{u_i}$ take prime values, and moreover, each of those $k+1$
conditions should be satisfied with ``probability'' $1/(\log X)$ for
$Nz\leq X$.
\par
Since the congruence condition limits the values of $z$ allowed, this
justifies that Proposition~\ref{pr-upper-indiv} gives the asymptotic
behavior, up to the  arithmetic factor $\phi_k(\underline{u})$, which
is very small; the asymptotic behavior should be
\begin{equation}\label{eq-asym-conj}
\somp_{\underline{u}}^+(X)\sim
c(\underline{u})\frac{X}{N[\underline{u}](\log X)^{k+1}}
\end{equation}
as $X\ra +\infty$,
for some (more complicated) arithmetic function $c(\underline{u})\leq
\phi_k(\underline{u})$.
Any other heuristic confirms this, of course; that based
on cancellation in long averages involving the M\"obius
function could in theory provide a prediction for the value of
$c(\underline{u})$ as an Euler product.
\par
If it seems reasonable to expect that~(\ref{eq-asym-conj}) holds, one
may also expect that it does uniformly at least in a range
$N[\underline{u}]\leq X^{\delta}$ for some $\delta>0$, and this would
provide a lower bound for $\somp^+(X)$ of the same order of
magnitude.
\par
The reader will easily convince herself that all other overcounting
done in deriving Theorem~\ref{th-upper-cm} should have at most the
effect of introducing a multiplicative constant: this includes the
step from Frobenius elements $\frobs(E)$ to all prime elements in
$K$ and the overcounting used in the proof of
Proposition~\ref{pr-somme} (because of the logarithmic scaling of that
sum).
\end{remark}



\subsection{Twin-primes in quadratic fields}
\label{ssec-twin-one}

In this section we prove Proposition~\ref{pr-upper-indiv}.
The argument is cleaner when the ideal $[\underline{u}]$ is principal:
the reader may assume that it is so in a first reading.
\par
Apart from the fact that we work over a quadratic field, the problem
is quite standard, and the proof will be close to, for instance, the
arguments in~\cite[\S 3]{bom-ls}.
\par
We will use the large sieve for $K$, in the version given by
Huxley~\cite[Th. 2]{huxley}. First some notation: for an integral ideal
$\ideal{n}$ of $K$, we denote by $(\Oc_K/\ideal{n})\spcheck$ the group of
additive characters of $\Oc_K/\ideal{n}$. We write 
$$
\sums_{\psi}{\alpha_{\psi}}
$$
for a sum over the \emph{primitive} characters of $\Oc_K/\ideal{n}$,
i.e. those which are not induced by a character modulo $\ideal{m}$ for
some $\ideal{m}$ dividing $\ideal{n}$. Also we denote by
$$
\sumb_{\ideal{n}}{\alpha_{\ideal{n}}}
$$
a sum over \emph{squarefree} ideals $\ideal{n}$.
\par
We define the height of $z\in\Oc_K$ by
$$
h(z)=\max (|a|,|b|)
$$
for $z=a+b\omega$. There exists a constant $\kappa>0$ such that 
\begin{equation}
  \label{eq-incl-n-h}
Nz\leq X \text{ implies } h(z)\leq \kappa \sqrt{X}
\text{ for }X>0\text{ and } z\in\Oc_K
\end{equation}
(one can take $\kappa=2$ for all $K$ if the basis $(1,\omega)$ is the
``canonical'' one).

\begin{theorem}\label{th-huxley}(Huxley)
Let $K/\Qq$ be an imaginary quadratic field. We have
\begin{equation}\label{eq-ls}
\sum_{N\ideal{n}\leq Q}{
\sums_{\psi\in (\Oc_K/\ideal{n})\spcheck}
{\Bigl|
\sum_{h(z)\leq X}{a(z)\psi(z\,\mathrm{mod}\ \ideal{n})}
\Bigr|^2}}
\ll (X^2+Q^2) \sum_{h(z)\leq X}{|a(z)|^2},
\end{equation}
where $(a(z))$ is any sequence of complex numbers, $Q$ and $X$
are any real numbers $\geq 1$.
The implied constant is absolute.
\end{theorem}

From this, proceeding as in~\cite[Th. 6]{bom-ls}, we derive an
arithmetic sieve result: a \emph{sieve} here is a pair $(M,\Omega)$
where 
$$
M=\{z\in\Oc_K\,\mid\, h(z)\leq X\}
$$
for some $X\geq 1$ and $\Omega$ is a map which associates 
a subset $\Omega(\ideal{p})\subset \Oc/\ideal{p}$ to prime ideals
$\ideal{p}$ with norm $N\ideal{p}\leq Q$. We denote
$\omega(\ideal{p})=|\Omega(\ideal{p})|$. The corresponding
\emph{sifted set} is
\begin{equation}\label{eq-sifted-set}
\mathcal{M}=
\{
z\,\mid\, h(z)\leq X\text{ and }
z\mods{\ideal{p}}\notin \Omega(\ideal{p})
\text{ for all }\ideal{p}
\}.
\end{equation}

\begin{corollary}\label{cor-sieve}
Let $K/\Qq$ be an imaginary quadratic field
and $(M,\Omega)$ a sieve. We have
$$
|\mathcal{M}|\ll \frac{X^2+Q^2}{J}
$$
where 
$$
J = \sumb_{N\ideal{n}\leq Q}{J(\ideal{n})},
\text{ with }
J(\ideal{n})=\prod_{\ideal{p}\mid \ideal{n}}
{\frac{\omega(\ideal{p})}{N\ideal{p}-\omega(\ideal{p})}}
$$
for $\ideal{n}$ squarefree.
The implied constant is absolute.
\end{corollary}

We will apply Corollary~\ref{cor-sieve} to the situation of
Proposition~\ref{pr-upper-indiv}. 
\par
\medskip
To setup the situation, we use the ideal-class group of
$K$. Let $\ideal{a}=[\underline{u}]$, let $\ideal{b}_0$ be an integral
ideal of $K$ with minimal norm in the ideal class 
inverse to that of $\ideal{a}$, say
$$
\ideal{a}\ideal{b}_0=(a_0)\text{ for some } a_0\in\Oc_K.
$$
If $z\in\Oc_K$ satisfies $z\equiv
1\mods{\ideal{a}}$, 
there exists an integral ideal $\ideal{b}$ such that
$$
(z-1)=\ideal{a}\ideal{b}
$$
and since $(z-1)$ is principal, $\ideal{b}$ and $\ideal{b}_0$ are in the
same ideal class, i.e. there exists $b\in K^{\times}$ such that
\begin{equation}
\label{eq-zmu}
\ideal{b}=b\ideal{b}_0, \text{ hence }
(z-1)=(a_0 b).  
\end{equation}
Since $\ideal{b}$ is integral, the denominator of $b$ is bounded (by
that of $N\ideal{b}_0$), i.e. there exists $d_0\in\Zz$, independent of
$z$ and with $d_0\mid N\ideal{b}_0$, such that
$b=c/d_0$ for some $c\in \Oc_K$.
\par
Hence, using~(\ref{eq-zmu}), there exists a unit
$\eps\in\Oc_K^{\times}$ (a finite group of order $\leq 6$) such that
\begin{equation}\label{eq-z-c}
z=\eps\frac{ a_0 c}{d_0}+1.
\end{equation}
Therefore
\begin{equation}\label{eq-1-bis}
\somp^+_{\underline{u}}(X)
\leq \sum_{\eps\in\Oc_K^{\times}}
{|\{ c\in \Oc_K\,\mid\, 
z=\eps \frac{a_0c}{d_0}+1\text{ satisfies }
N(z)\leq X \text{ and }
f_{u_i}(z)\text{ is prime}
\}|}.
\end{equation}
By~(\ref{eq-z-c}) and the definition of $a_0$, $d_0$, if $Nz\leq X$ we
have 
\begin{equation}\label{eq-n-c}
Nc\leq \frac{d_0^2}{Na_0}X^+\leq \frac{X^+}{N\ideal{a}}.
\end{equation}
For $\eps\in\Oc_K^{\times}$, consider the sieving problem
$(M,\Omega_{\eps})$ consisting in sieving the set
$$
M=\{z\in\Oc_K\,\mid\, h(z)\leq \kappa(\sqrt{(X+1)/N\ideal{a}})\}
$$
(where $\kappa$ is as in~(\ref{eq-incl-n-h})),
by prime ideals $\ideal{p}$ with $N\ideal{p}\leq \sqrt{X/N\ideal{a}}$,
with $\Omega_{\eps}(\ideal{p})$ defined as follows: let
\begin{equation}
  \label{eq-omega-1}
\Omega^+(\ideal{p})=\{ -\frac{u_i d_0}{\eps \bar{u}_i a_0}\,\mid\,
1\leq i\leq k\} \subset (\Oc_K/\ideal{p})
\end{equation}
(with the convention that any ratio where the denominator is $0$
modulo $\ideal{p}$ is omitted), and define
$$
\Omega_{\eps}(\ideal{p})=
\begin{cases}
\Omega^+(\ideal{p})&\text{ if } |\Omega^+(\ideal{p})|=k\\
\emptyset&\text{ otherwise.}
\end{cases}
$$

\begin{lemma}
Let $\mathcal{M}_{\eps}$ denote the sifted set for the sieving problem
above.
We have
$$
\somp^+_{\underline{u}}(X)
\leq \sum_{\eps}{|\mathcal{M}_{\eps}|}.
$$
\end{lemma}

This is an immediate consequence of the previous
inequality~(\ref{eq-1-bis}),~(\ref{eq-incl-n-h}) and the definition of the
sieve (one could of course be more precise and not disregard the
primes $\ideal{p}$ with $|\Omega^+(\ideal{p})|<k$).

\begin{lemma}
We have $\omega_{\eps}(\ideal{p})=0$ if and only if
$$
\ideal{p} \mid N[\underline{u}]\disc(\underline{u})
$$
where the discriminant is defined in~(\ref{eq-disc-def}).
\end{lemma}

\begin{proof}
This is clear: the factor $N[\underline{u}]=N\ideal{a}$ arises from
the possibility that the denominators in~(\ref{eq-omega-1}) are
divisible by $\ideal{p}$, whereas the discriminant occurs from the
possibility that
$$
\frac{u_i}{\bar{u}_i}
=\frac{u_j}{\bar{u}_j}\mods{\ideal{p}}
\text{ i.e. }
u_i \bar{u}_j - \bar{u}_i u_j=0\mods{\ideal{p}}
$$
for some $i\not= j$.
\end{proof}

Note that $\disc(\underline{u})\not=0$ in the application to
Proposition~\ref{pr-upper-indiv} since no two $u_i$ coincide.
\par
\medskip
By Corollary~\ref{cor-sieve}, we deduce that
\begin{equation}
  \label{eq-sieve-res}
\somp^+_{\underline{u}}(X)
\ll \frac{1}{J}\frac{X}{N\ideal{a}}
\end{equation}
for $X\geq 2$ with an absolute implied constant, where
$$
J=\sumb_{\stacksum{N\ideal{n}\leq \sqrt{X/N\ideal{a}}+1}
{(\ideal{n},\disc(\underline{u})N\ideal{a})=1}}
{J(\ideal{n})},\text{ with }
J(\ideal{n})=
\prod_{\ideal{p}\mid \ideal{n}}
{\frac{k}{N\ideal{p}-k}}.
$$
Note that $\ideal{n}\mapsto J(\ideal{n})$ is an arithmetic function
that depends only on $k$, not on $\underline{u}$.
\par
It only remains to find a lower bound for $J$ to get an upper bound
for $\somp^+_{\underline{u}}(X)$; the only issue of note is the
uniformity in $\underline{u}$. All the arguments below are standard
(see e.g.~\cite[Th. 2.4]{halb-rich}), but by lack of a convenient
reference, especially in the context of a number field, we give all
details.
\par
For $\ideal{n}$ squarefree we have
$J(\ideal{n})\geq J^{\flat}(\ideal{n})$, 
where $J^{\flat}(\ideal{n})$ is the totally multiplicative arithmetic
function on integral ideals of $K$ such that
$$
J^{\flat}(\ideal{p})=
\begin{cases}
0&\text{if }N\ideal{p}\leq k^2\\
k/N\ideal{p}&\text{otherwise.}
\end{cases}
$$
Therefore
\begin{equation}\label{eq-res-1}
J\geq \sumb_{\stacksum{N\ideal{n}\leq \sqrt{X/N\ideal{a}}+1}
{(\ideal{n},\disc(\underline{u})N\ideal{a})=1}}
{J^{\flat}(\ideal{n})}.
\end{equation}
We consider the generating series
$$
z(s)=\sumb_{\ideal{n}}{J^{\flat}(\ideal{n})
(N\ideal{n})^{-s}}=
\prod_{N\ideal{p}>k^2}{(1+k(N\ideal{p})^{-s-1})}
$$
which converges absolutely\footnote{ In particular, has no zero.} for
$\Reel(s)>0$, and the closely related 
$$
w(s)=\sum_{\ideal{n}}{J^{\flat}(\ideal{n})
(N\ideal{n})^{-s}}=
\prod_{N\ideal{p}>k^2}{(1-k(N\ideal{p})^{-s-1})^{-1}}
$$
which also converges absolutely in the same region.

\begin{lemma}\label{lm-compar-zw}
There exists a Dirichlet series
$$
y(s)=\sum_{\ideal{n}}{Y(\ideal{n})
(N\ideal{n})^{-s}}
$$
such that
\begin{equation}\label{eq-factor}
z(s)=y(s)w(s),
\end{equation}
and $y(s)$ converges absolutely for $\Reel(s)>-1/2$. 
\end{lemma}

\begin{proof}
This is clear by comparing the Euler factors of $z(s)$ and $w(s)$,
using the fact that the zeros of $1-k(N\ideal{p})^{-s-1}$ have
$\Reel(s)<-1/2$ for $N\ideal{p}>k^2$.
\end{proof}

\begin{lemma}\label{lm-k}
There exists a constant $c>0$ such that
$$
\sumb_{N\ideal{n}\leq Y}{J^{\flat}(\ideal{n})}
=c (\log Y)^k +\oun((\log Y)^{k-1})
$$
for $Y\geq 2$.
\end{lemma}

\begin{proof}
This is obvious by comparison of $w(s)$ with $\zeta_K(s+1)^k$, which
has a pole of order $k$ at $s=0$, and contour integration: we have (as in
Lemma~\ref{lm-compar-zw}) 
$$
w(s)=\zeta_K(s+1)w_1(s)
$$
for some Dirichlet series $w_1(s)$ which converges absolutely in the
region $\Reel(s)>-1/2$. 
\end{proof}

\begin{lemma}\label{lm-arith-factor}
Let $T(\ideal{n})$ be a completely multiplicative arithmetic function
of integral ideals of $K$ such that:
\par\indent
\emph{(i)} There exists $A>0$ such that
$$
T(\ideal{p})\leq \frac{A}{N\ideal{p}}
$$
for all prime ideals $\ideal{p}$.
\par\indent
\emph{(ii)} There exists $c>0$ and $\gamma>0$ such that
$$
\sum_{N\ideal{n}\leq Y}{T(\ideal{n})}
= c (\log Y)^{\gamma}+\oun((\log Y)^{\gamma-1})
$$
for $Y\geq 2$.
\par
Fix an integer $B\geq 1$. Then for all $Y\geq 2$ and all non-zero
integral ideals $\ideal{q}$ such that $N\ideal{q}\leq Y^B$, we have
$$
\sum_{\stacksum{N\ideal{n}\leq Y}{(\ideal{n},\ideal{q})=1}}
{T(\ideal{n})}= (\log Y)^{\gamma}
\Bigl(
\sum_{\ideal{d}\mid \ideal{q}}{\mu(\ideal{d})T(\ideal{d})}
\Bigr)
+\oun((\log Y)^{\gamma-1}(\log \log Y)^{A}),
$$
the implied constant depending on $A$, $B$ and $\gamma$.
\end{lemma}

In this statement and in the proof, we use $d(\ideal{n})$ 
(resp. $\mu(\ideal{n})$) to denote the divisor function (resp. M\"obius)
function for integral ideals. The latter is defined as usual
(i.e. $\mu(\ideal{p}^k)=(-1)^k$ for any prime ideal $\ideal{p}$ and
$k\geq 0$, and $\mu$ multiplicative). The M\"obius inversion formula
holds:
$$
\sum_{\ideal{d}\mid\ideal{n}}{\mu(\ideal{d})}=
\begin{cases}
1&\text{ if }\ideal{n}=1\\
0&\text{ otherwise}.
\end{cases}
$$
We will use the following easy estimate
\begin{equation}\label{eq-est}
\prod_{\ideal{p}\mid\ideal{n}}
{\Bigl(1+\frac{A}{N\ideal{p}}\Bigr)}
\ll (\log \log N\ideal{n})^A
\end{equation}
for all non-zero integral ideals $\ideal{n}$. The implied constant
depends only on $A$.

\begin{proof}
We have by M\"obius inversion
\begin{align*}
\sum_{\stacksum{N\ideal{n}\leq Y}{(\ideal{n},\ideal{q})=1}}
{T(\ideal{n})}&=
\sum_{\ideal{d}\mid\ideal{q}}{\mu(\ideal{d})
\sum_{N\ideal{n}\leq Y/N\ideal{d}}{
T(\ideal{n}\ideal{d})
}}\\
&=\sum_{\stacksum{\ideal{d}\mid\ideal{q}}{N\ideal{d}<Y^{\delta}}}
{\mu(\ideal{d})T(\ideal{d})
\sum_{N\ideal{n}\leq Y/N\ideal{d}}{
T(\ideal{n})
}}+
\oun(Y^{-\delta+\eps})
\end{align*}
for any (fixed) $\delta>0$ and $\eps<\delta$, having used the complete
multiplicativity, and~(i) and~(ii) to estimate the remaining sum over
large divisors of $N\ideal{q}$:
\begin{align*}
\sum_{\stacksum{\ideal{d}\mid\ideal{q}}{N\ideal{d}\geq Y^{\delta}}}
{\mu(\ideal{d})T(\ideal{d})
\sum_{N\ideal{n}\leq Y/N\ideal{d}}{
T(\ideal{n})
}}&\ll
\frac{(\log Y)^{\gamma}}{Y^{\delta}}
\sum_{\ideal{d}\mid\ideal{q}}{d(\ideal{d})^{A}}\\
&\ll \frac{(\log Y)^{\gamma}}{Y^{\delta}}
d(\ideal{q})^{A+1}
\\
&\ll_{\eps,B,A} Y^{-\delta+\eps}
\end{align*}
(by the assumption $N\ideal{q}\leq Y^B$). The implied constant depends
on $\eps$, $A$, $B$ and $\gamma$.
\par
Using again~(ii) we have
\begin{align*}
\sum_{\stacksum{\ideal{d}\mid\ideal{q}}{N\ideal{d}<Y^{\delta}}}
{\mu(\ideal{d})T(\ideal{d})
\sum_{N\ideal{n}\leq Y/N\ideal{d}}{
T(\ideal{n})
}}&=
c\sum_{\stacksum{\ideal{d}\mid\ideal{q}}{N\ideal{d}<Y^{\delta}}}
{\mu(\ideal{d})T(\ideal{d})
\Bigl(\Bigl(\log \frac{Y}{N\ideal{d}}\Bigr)^{\gamma}
+\oun( (\log Y)^{\gamma-1})\Bigr)}\\
&=c(\log Y)^{\gamma}
\sum_{\stacksum{\ideal{d}\mid\ideal{q}}{N\ideal{d}<Y^{\delta}}}
{\mu(\ideal{d})T(\ideal{d})}
+\oun((\log Y)^{\gamma-1}(\log \log Y)^A)
\end{align*}
by expanding the logarithm and estimating
\begin{align*}
(\log Y)^{\gamma-1}
\Bigl|\sumb_{\stacksum{\ideal{d}\mid\ideal{q}}
{N\ideal{d}< Y^{\delta}}}
{\mu(\ideal{d})T(\ideal{d})}\Bigr|
&\leq (\log Y)^{\gamma-1}
\prod_{\ideal{p}\mid\ideal{q}}
{(1+T(\ideal{p}))}\\
&\leq (\log Y)^{\gamma-1}
\prod_{\ideal{p}\mid\ideal{q}}
{\Bigl(1+\frac{A}{N\ideal{p}}\Bigr)}
\ll_A (\log \log N\ideal{q})^A\\
&\ll_{A,B} (\log \log Y)^A
\end{align*}
by~(\ref{eq-est}).
\par
It remains to get rid of $\delta$, which is possible since
$$
\sum_{\stacksum{\ideal{d}\mid\ideal{q}}
{N\ideal{d}\geq Y^{\delta}}}
{\mu(\ideal{d})T(\ideal{d})}\ll
Y^{-\delta}\sum_{\ideal{d}\mid\ideal{q}}
{d(\ideal{d})^A}\ll Y^{-\delta+\eps}.
$$
Choosing $\delta$ small enough and $\eps<\delta$, the lemma follows.
\end{proof}

We come back to~(\ref{eq-res-1}) and write, using~(\ref{eq-factor})
$$
\sumb_{\stacksum{N\ideal{n}\leq Y}
{(\ideal{n},\disc(\underline{u})N\ideal{a})=1}}
{J^{\flat}(\ideal{n})}=
\sum_{\stacksum{N\ideal{m}\leq Y}
{(\ideal{m},\disc(\underline{u})N\ideal{a})=1}}
{Y(\ideal{m})\sum_{\stacksum{N\ideal{n}\leq Y/N\ideal{m}}
{(\ideal{n},\disc(\underline{u})N\ideal{a})=1}}{
J^{\flat}(\ideal{n})
}}.
$$
To the inner sum we can apply Lemma~\ref{lm-arith-factor} with
$\ideal{q}=\disc(\underline{u})N\ideal{a}$ and $\gamma=k$:
the assumptions hold for some $A$, $B$ and $\gamma=k$ by
Lemma~\ref{lm-k}. Therefore
\begin{align}
\sumb_{\stacksum{N\ideal{n}\leq Y}
{(\ideal{n},\disc(\underline{u})N\ideal{a})=1}}
{J^{\flat}(\ideal{n})}&=
c\sum_{\stacksum{N\ideal{m}\leq Y}
{(\ideal{m},\disc(\underline{u})N\ideal{a})=1}}
{Y(\ideal{m})
\bigl(\log (Y/N\ideal{a})^{k}
+\oun( (\log Y)^{k-1}(\log \log Y)^{k})\bigr)}
\nonumber\\
&=c(\log Y)^k\Bigl(
\sum_{\ideal{d}\mid \disc(\underline{u})N\ideal{a}}
{\mu(\ideal{d})J^{\flat}(\ideal{d})}
\Bigr)
\sum_{\stacksum{N\ideal{m}\leq Y}
{(\ideal{m},\disc(\underline{u})N\ideal{a})=1}}
{Y(\ideal{m})}\label{eq-res-2}\\
&\quad\quad+\oun( (\log Y)^{k-1}(\log \log Y)^{k}),\nonumber
\end{align}
by again expanding the logarithm, and using the fact that for any
$B\geq 0$ the series
$$
\sum_{\ideal{m}}{Y(\ideal{m})(\log
N\ideal{m})^B}
$$
is absolutely convergent. Now apply the following lemma to
$Y(\ideal{m})$ and $\ideal{q}=\disc(\underline{u})N\ideal{a}$:

\begin{lemma}
\label{lm-acv}
Let $Y(\ideal{n})$ be a multiplicative arithmetic function,
$y(s)$ its generating Dirichlet series. Assume that the Euler product
for $y(s)$ converges absolutely for $\Reel(s)>-1/2$. Then for any
non-zero integral ideal $\ideal{q}$ we have
$$
\sum_{\stacksum{N\ideal{m}\leq Y}{(\ideal{m},\ideal{q})=1}}
{Y(\ideal{m})}\gg 1
$$
for $X\geq 2$, the implied constant depending only on the function
$Y$.
\end{lemma}

\begin{proof}
By a standard application of contour integration and Perron's
formula. The size of $\ideal{q}$ does not matter here because the sum
always involves $\ideal{m}=1$, with a contribution $=1$. In slightly
more detail: it is well-known (see e.g.~\cite[]{tit-zeta}) that
$$
\frac{1}{2i\pi}\int_{1-iT}^{1+iT}
{y^s \frac{ds}{s}}
=h(y)
+\oun\Bigl(\frac{y}{T|\log y|}\Bigr)
$$
for all $y>0$ and $T>0$, where $h(y)=1$ for $y>1$, $h(y)=0$ for $y<1$
and $h(1)=1/2$.
\par
Let $y_{\ideal{q}}(s)$ be the generating Dirichlet series of
$Y(\ideal{m})$ restricted to those $\ideal{m}$ coprime to
$\ideal{q}$. Choosing $X$ of the form $1/2+m$ for some integer $m$, as
we may without loss of generality, we have
\begin{equation}\label{eq-int-1}
\frac{1}{2i\pi}\int_{1-iT}^{1+iT}
{y_{\ideal{q}}(s) X^s \frac{ds}{s}}
=
\sum_{\stacksum{N\ideal{m}\leq X}{(\ideal{m},\ideal{q})=1}}
{Y(\ideal{m})}
+\oun(XT^{-1})
\end{equation}
since
$$
\sum_{\ideal{m}}{\frac{Y(\ideal{m})}
{N\ideal{m}|\log (X/N\ideal{m})|}}< + \infty
$$
(use the absolute convergence of $\sum{Y(\ideal{m})}$
and $|N\ideal{m}(\log X/N\ideal{m})|\gg 1$).
\par
On the other hand, by Cauchy's theorem we have
\begin{equation}\label{eq-int-2}
\frac{1}{2i\pi}\int_{\mathcal{C}}{y_{\ideal{q}}(s)X^s\frac{ds}{s}}
=y_{\ideal{q}}(0)
\end{equation}
where $\mathcal{C}$ is the boundary of the rectangle $[-1/4, 1] \times
[-T, T]$. By absolute convergence, the integral on the horizontal
pieces and on the vertical line $\Reel(s)=-1/4$ are
\begin{align*}
\frac{1}{2i\pi}\Bigl\{
\int_{1/4-iT}^{1-it}{}
+\int_{1+iT}^{-1/4+iT}{} 
y_{\ideal{q}}(s) X^s\frac{ds}{s}
\Bigr\}
&\ll XT^{-1}\\
\frac{1}{2i\pi}\int_{-1/4-iT}^{-1/4+it}{y_{\ideal{q}}(s)
X^s\frac{ds}{s}}
&\ll X^{-1/4}
\end{align*}
the implied constant depending only on
$Y$. Hence~(\ref{eq-int-1}) and~(\ref{eq-int-2}) show that 
$$
\sum_{\stacksum{N\ideal{m}\leq X}{(\ideal{m},\ideal{q})=1}}
{Y(\ideal{m})}
=y_{\ideal{q}}(0)+\oun(XT^{-1})+\oun(X^{-1/4}).
$$
Taking $T=X^2$ for instance gives
$$
\sum_{\stacksum{N\ideal{m}\leq X}{(\ideal{m},\ideal{q})=1}}
{Y(\ideal{m})}\gg y_{\ideal{q}}(0)
$$
the implied constant depending only on $Y$.
\par
Since $y_{\ideal{q}}(0)$ is the same absolutely
convergent Euler product as $y(0)$, except that primes dividing
$\ideal{q}$ are omitted, and
any partial product of an absolutely
convergent infinite product has a uniform lower bound, it follows that
$$
y_{\ideal{q}}(0)\gg 1,
$$
thereby proving the lemma.
\end{proof}

Since moreover
$$
\sum_{\ideal{d}\mid \disc(\underline{u})N\ideal{a}}
{\mu(\ideal{d})J^{\flat}(\ideal{d})}
=\prod_{\ideal{p}\mid  \disc(\underline{u})N\ideal{a}}
{\bigl(1-J^{\flat}(\ideal{p})\bigr)}>0,
$$
because $J^{\flat}(\ideal{p})<1$ for all $\ideal{p}$ (this is why small
primes had to be excluded), the inequality~(\ref{eq-res-2}) proves
that
\begin{equation}\label{eq-res-3}
J\gg \prod_{\ideal{p}\mid  \disc(\underline{u})N\ideal{a}}
{\bigl(1-J^{\flat}(\ideal{p})\bigr)}
\end{equation}
the implied constant depending on $k$ and $K$ only.

\begin{lemma}
For all $\ideal{p}$ we have
$$
1-J^{\flat}(\ideal{p})\geq (1-k^{-2})\Bigl(1+\frac{k}{N\ideal{p}}\Bigr).
$$
\end{lemma}
\begin{proof}
This is obvious from the definition.
\end{proof}

Proposition~\ref{pr-upper-indiv} follows
from~(\ref{eq-sieve-res}),~(\ref{eq-res-3}) and this lemma.  

\subsection{Proof of Proposition~\ref{pr-somme}}
\label{ssec-twin-somme}

In this section we prove Proposition~\ref{pr-somme}. For $k=0$, the
result is obvious with no need of the factor $(\log X)^{\eps}$, since
the sum is reduced to $u=1$. So we assume $k\geq 1$.
\par
We have by~(\ref{eq-est})
$$
\phi_k(1,\underline{u})\ll (\log \log X)^k
$$
with an absolute implied constant, hence by positivity
\begin{equation}\label{eq-positivity}
\sum_{\stacksum{\underline{u}\in\prim^k}
{Nu_i\leq X}}{\frac{\phi_k(1,\underline{u})}{N[\underline{u}]}}
\ll (\log\log X)^k\sum_{n\leq X^k}{\frac{\rho(n)}{n}}
\end{equation}
for $X\geq 2$ (the constant depending only on $k$), where $\rho(n)$ is
the arithmetic function defined by 
\begin{equation}\label{eq-rho}
\rho(n) = |\{ (u_1,\ldots,u_n)\text{ ideals in } \Oc_K\,\mid\,
N[u_1,\ldots,u_n]=n\}|.
\end{equation}
Thus we drop the condition that the $u_i$ be integers or primitive,
and drop the size condition $Nu_i\leq X$ on the solutions 
of $N[\underline{u}]=n$, and this shouldn't change the order of
magnitude because of the logarithmic weight.
\par
The arithmetic function $\rho(n)$ is multiplicative.

\begin{lemma}\label{lm-rho-borne}
Let $n\geq 1$ be an integer. We have
$$
\rho(n)\leq d(n)^{2k}
$$
where $d(n)$ is the function ``number of divisors''.
\end{lemma}

\begin{proof}
In~(\ref{eq-rho}), $Nu_i\mid n$ for all $i$, so there are at
most $d(n)^k$ choices of $(Nu_1,\ldots,Nu_k)$, and for each
of those there are
$$
r(Nu_1)\cdots r(Nu_k)\leq r(n)^k\leq d(n)^k
$$
choices of $(u_1,\ldots u_k)$.
\end{proof}

\begin{lemma}\label{lm-rho-p}
Let $p$ be a prime number.
We have
$$
\rho(p)=(1+\chi(p))(2^k-1).
$$
\end{lemma}

\begin{proof}
We have $N[u_1,\ldots,u_k]=p$ if and only if 
\begin{equation}\label{eq-interim}
[u_1,\ldots,u_k]=\pi
\end{equation}
where $\pi$ is an ideal such that $N\pi=p$.
\par
For a given $\pi$, the solutions
$\underline{u}$ of $N\underline{u}=p$ correspond bijectively to 
$k$-tuples of integers $(\nu_1,\ldots \nu_k)$ such that
$$
u_i=\pi^{\nu_i},
$$
with $0\leq \nu_i\leq 1$ and at least one of the $\nu_i$ is $=1$.
The number of such tuples is equal to $2^k-1$ (all tuples except
$(0,\ldots,0)$).
\par
The number of solutions of $N\pi=p$ is $1+\chi(p)$ for all
primes $p$, and the lemma follows.
\end{proof}

Proposition~\ref{pr-somme} is a consequence of~(\ref{eq-positivity}) and
Lemmas~\ref{lm-rho-borne} and~\ref{lm-rho-p}, applying to $\rho$
the following very standard result (compare
Section~\ref{ssec-twin-one}) applied with $\gamma=2^k-1$. 
\begin{lemma}
Let $\rho(n)$ be a multiplicative arithmetic function satisfying: 
\par\noindent
\emph{(i)} There exists $A>0$ such that 
\begin{equation}\label{eq-i}
\rho(n)\leq d(n)^A\text{ for all } n\geq 1,
\end{equation}
\emph{(ii)} There exists an integer $\gamma$ such that for all primes
$p$ we have $\rho(p)=\gamma(1+\chi(p))$.
\par
Then there exists $c>0$ such that
$$
\sum_{n\leq X}{\frac{\rho(n)}{n}}\sim c(\log X)^{\gamma}
$$
as $X\ra +\infty$.
\end{lemma}

\begin{proof}
Let
$$
z(s)=\sum_{n\geq 1}{\rho(n)n^{-s}}
$$
be the Dirichlet generating series of $\rho$. By~(i), the series
converges and defines a holomorphic function for $\Reel(s)>1$.
By multiplicativity, $z(s)$ has an absolutely convergent Euler product
expansion
$$
z(s)=\prod_{\chi(p)=1}{(1+2 \gamma p^{-s} +
\rho(p^2) p^{-2s} + \ldots )}
\prod_{\chi(p)=0}{(1+\gamma p^{-s}+\ldots)}
\prod_{\chi(p)=-1}{(1+\rho(p^2)p^{-2s}+\ldots)}.
$$
Hypothesis (ii) implies that one can factorize
$$
z(s)=\zeta_K(s)^{\gamma}z_1(s)
$$
where $z_1(s)$, first defined by this equation for $\Reel(s)>1$,
admits analytic continuation to a holomorphic function on
$\Reel(s)>1/2$. Indeed one has
$$
\zeta_K(s)=\prod_{\chi(p)=1}{(1-2 p^{-s}+p^{-2s})^{-1}}
\prod_{\chi(p)=-1}{(1-p^{-2s})^{-1}},
$$
so the products over split and inert primes already converge for
$\Reel(s)>1/2$, while the coefficient of $p^{-s}$ in the $p$-Euler
factor for $z_1(s)$ vanishes.
\par
Since $\zeta_K(s)$ has a simple pole at $s=1$, it follows that $z(s)$
has a pole of order $\gamma$ at $s=1$, so a standard contour
integration proves the lemma.
\end{proof}

For $k=1$, we can easily get rid of the annoying factor $\log \log X$, as
mentioned in Remark~\ref{rem-annoying}.

\begin{proposition}
We have
$$
\sum_{\stacksum{u\in\prim}{Nu\leq X}}{\frac{\phi_1(1,u)}{Nu}}
\ll (\log X)
$$
for $X\geq 2$, the implied constant depending only on $K$.
\end{proposition}

\begin{proof}
We allow ourself to be a little sketchy: we have
$$
\phi_1(1,u)=\prod_{\ideal{p}\mid Nu (u-\bar{u})}{(1+(N\ideal{p})^{-1})}.
$$
Assume $K=\Qq(\sqrt{-4D})$ with $4D$ a fundamental discriminant
$4D\equiv 0\mods{4}$ so that 
$(1,\sqrt{-D})$ is a $\Zz$-basis of $\Oc_K$ (the remaining case being
similarly treated) and $N(a+b\sqrt{-D})=a^2+Db^2$.
\par
By trivial estimate, we have for $u=a+b\sqrt{-D}$
$$
\phi_1(1,u)\leq \frac{\psi(2D)}{2D}\frac{\psi(a^2+Db^2)}{a^2+Db^2}
\frac{\psi(b)}{b},
$$
(recall $\psi$ is defined in~(\ref{eq-psi})).
Hence
$$
\sum_{\stacksum{u\in\prim}{Nu\leq X}}{\frac{\phi_1(1,u)}{Nu}}
\leq 2
\sum_{d\leq X}{\frac{\mu(d)^2}{d}\sum_{0\leq |b|\leq (X/D)^{1/2}}
{\frac{\psi(b)}{b}\sum_{\stacksum{0\leq a\leq \sqrt{X-Db^2}}{d\mid
a^2+Db^2}}{\frac{1}{a^2+Db^2}}}}.
$$
The contribution of $b=0$ is $\ll 1$ (since $d\mid a^2$ and $d$
squarefree imply $d\mid a$). For $|b|\geq 1$, in the inner sum we
write $a=da_1+\alpha$ for some $\alpha$, $0\leq \alpha<d$, such that
$\alpha^2=-Db^2\mods{d}$. For given $\alpha$, by partial summation, the
inner sum over $a_1$ is easily seen to be $\ll (bd\sqrt{D})^{-1}$, uniformly in
$\alpha$. The result then follows since the number of $\alpha$ for a
given squarefree $d$ is at most the number of divisors of $d$, and
$$
\sum_{b\leq (X/D)^{1/2}}{\frac{\psi(b)}{b^2}}\ll \log X,\text{ and }
\sum_{n\geq 1}{\frac{d(n)\mu(n)^2}{n^2}}<+\infty.
$$
\end{proof}

Extending this kind of argument for $k\geq 2$ might be possible
although certainly cumbersome since the various $u_i$ would become
mixed up together. The issue is whether $\disc(\underline{u})$ can
have too often too small prime factors, and doesn't seem completely
trivial.


\subsection{The elliptic splitting problem}
\label{ssec-cm-split}

Because the condition $d\mid \dun(\ideal{p})$ is equivalent to the
congruence $\frob_{\ideal{p}}\equiv 1\mods{d}$ in the endomorphism
ring of $E$, we can again apply sieve to obtain a Brun-Titchmarsh
inequality for totally split primes in $K(E[d])$ for a CM curve.
In particular, the extension $K(E[\infty])/K$ is a Brun-Titchmarsh
field for a CM curve.



\begin{theorem}\label{th-bt-cm}
Let $E/H$ be a CM curve with complex multiplication by an order $\Oc$ of a
quadratic field $K/\Qq$, and $H'=H(E[\infty])$ its division
field. Assume that $H$ contains $K$. Then $H'/H$ is a Brun-Titchmarsh
field corresponding to the restriction of scalars
$G=\Res_{\Oc/\Zz}(\Gg_m)$.
\end{theorem}

First remark that the extension $H'/H$ enters in the setup described
in Section~\ref{ssec-bt} for the general Brun-Titchmarsh problem,
because of part 1. of Theorem~\ref{th-galois-nf} and the general
ramification properties of $E[d]$. 

\begin{proposition}
Let $H$ be a number field, $E/H$ an elliptic curve with CM by an order
$\Oc\subset K\subset H$ and let $d\geq 1$ be an integer. We have
$$
\pi_E(X;d,1)\ll [H:\Qq] \frac{X}{\varphi_{\Oc}(d)(\log X/d^2)}
$$
for $d\leq X$, where the implied constant is absolute and
$\varphi_{\Oc}(d)= |(\Oc/d\Oc)^{\times}|$.
\end{proposition}

This proposition clearly implies the theorem since
$\Gal(H(E[d])/H)$ is of bounded index in
$G(\Zz/d\Zz)=(\Oc/d\Oc)^{\times}$. In turn,
since $\ideal{p}$ is split in $H(E[d])/H$ if and only if the Frobenius
$\psi(\ideal{p})$ satisfies $\psi(\ideal{p})\equiv 1\mods{d}$, 
it follows immediately
from Lemma~\ref{lm-borne-cm} and the next proposition:
\begin{proposition}\label{pr-cm-split-bis}
Let $K/\Qq$ be an imaginary quadratic field. Then 
$$
\pi_K(X;d,1)\ll_K \frac{X}{\varphi_{K}(d)(\log X/d^2)}
$$
the implied constant depending only on $K$.
\end{proposition}

\begin{proof}
This is almost a (simpler) special case of
Proposition~\ref{pr-upper-indiv} (for $k=1$ with $d$ instead of $u$;
it is not included in that Proposition since the latter assumes
$u\not\in \Zz$), so we can be very 
sketchy. One applies the large sieve, as in
Section~\ref{ssec-twin-one}, to sieve
$$
M=\{z\in\Oc_K\,\mid\, h(z)\leq \sqrt{X}/d\}
$$
by primes $\ideal{p}$ with $N\ideal{p}\leq \sqrt{X}/d=Q$, with
$\Omega(\ideal{p})=\{-1/d\mods {\ideal{p}}\}$, if $\ideal{p}$ does not
divide $d$ and $\Omega(\ideal{p})=\emptyset$ otherwise. By
Corollary~\ref{cor-sieve} we derive
$$
\pi_K(X;d,1)\ll \frac{X}{d^2}\frac{1}{J}
$$
with 
$$
J=\sumb_{N\ideal{n}\leq Q}{\prod_{\stacksum{\ideal{p}\mid\ideal{n}}{
(d,\ideal{p})=1}}
{\frac{1}{N\ideal{p}-1}}}.
$$
Evaluating this sum in the usual manner, the result follows.
\end{proof}

Note the following simple corollary of Theorem~\ref{th-bt-cm} for the
elliptic splitting problem, which is still not very strong however
(recall the expected order of magnitude is $X$).

\begin{corollary}
Let $E/H$ and $K$ be as in the proposition. We have
$$
S_E(X;\dun)\ll_E X(\log X)^{1/2}
$$
for $X\geq 2$.
\end{corollary}

\begin{proof}
We split the sum
$$
S_E(X;\dun)=\sum_{d\leq \sqrt{X}+1}{\varphi(d)\pi_E(X;d,1)}
$$
in two ranges $d\leq B$ and $B<d\leq \sqrt{X}+1$ 
where $B=(\sqrt{X}+1)/A$ 
for some $A\geq 1$ to be chosen later. In the first range,
applying the Brun-Titchmarsh inequality yields
\begin{align*}
\sum_{d\leq B}{\varphi(d)\pi_E(X;d,1)}
&\ll_E \frac{X}{\log X/B^2}\sum_{d\leq B}
{\frac{\varphi(d)}{\varphi_K(d)}}\\
&\ll_E X\frac{\log (\sqrt{X}/A)}{\log A}.
\end{align*}
In the other range, we use instead the trivial bound coming from
Lemma~\ref{lm-borne-cm} and overcounting all integers $z\in\Oc_K$
instead of only prime elements, which gives
$$
\pi_E(X;d,1)\ll_K [H:\Qq] |\{ z\in\Oc_K\,\mid\,
Nz\leq X\text{ and } z\equiv 1\mods{d}\}|
\ll_K [H:\Qq] \Bigl(\frac{X}{d^2}+1\Bigr).
$$
Hence
$$
\sum_{B<d\leq \sqrt{X}+1}{\varphi(d)\pi_E(X;d,1)}
\ll_K [H:\Qq]X\sum_{B<d\leq \sqrt{X}+1}
{\frac{\varphi(d)}{d^2}}
\ll_K [H:\Qq] X \log A.
$$
We now choose $A=\exp(\sqrt{\log X})$ and it follows that
$$
S_E(X;\dun)\ll_E  X(\log X)^{1/2},
$$
as desired.
\end{proof}

\begin{remark}
The Brun-Titchmarsh property and the Bombieri-Vinogradov Theorem in
$K$ can be used to prove a (weak) lower bound
$$
S_E(X;\dun)\gg_E X\frac{\log \log X}{\log X}
$$
(better than the trivial lower bound $X/\log X$ arising by taking the
single term $d=1$ in~(\ref{eq-dun-som}) only by $\log\log X$).
The factor $\varphi(d)$ is the reason of the difficulties in the
direction of lower bounds.
\end{remark}

\mysection{Local study of totally split primes}
\label{sec-local}

We now change the point of view, motivated by the
considerations of the previous sections. We wish to understand, given
$d\geq 1$, for which finite fields $\Fp_q$ there does exist \emph{some}
elliptic curve $E/\Fp_q$ with $\dun(E)=d$, or more generally with
its $d$-torsion points rational over $\Fp_q$. In the cyclotomic case
the answer is simple: $\Fp_q$ contains all the $d$-th roots of unity
if and only if $q\equiv 1\mods{d}$. And the analogue of $\dun$ is the
largest $d$ for which all $d$-th roots of unity are in $\Fp_q$,
therefore it is simply $q-1$.
\par
We will first study this question using the methods introduced by
Deuring~\cite{deuring}. The results can also be extracted from papers
of Schoof~\cite{schoof-2}, Howe~\cite{howe}, Tsfasman-Vladut (and
maybe others I have not seen). But those are written with a slightly 
different emphasis. Then we recover similar results using
modular curves and the trace formula, before giving some applications.

\subsection{Results using endomorphism rings}

We first deal quickly with the case of supersingular elliptic
curves.

\begin{proposition}\label{pr-supersing}
Let $E/\Fp_q$ be a supersingular elliptic curve over a finite
field with characteristic $p$. We have
\begin{equation}\label{eq-supersing-small}
\dun(E)\leq 2
\end{equation}
unless $E$ satisfies $a(E)^2=4q$, in which case
\begin{equation}\label{eq-supersing-big}
\dun(E)=
\begin{cases}
\sqrt{q}+1&\text{ if $a(E)=-2\sqrt{q}$}\\
\sqrt{q}-1&\text{ if $a(E)=2\sqrt{q}$}.
\end{cases}
\end{equation}
\end{proposition}

\begin{proof}
All this is contained in~\cite[Lemma 4.8]{schoof-2} for instance, but
most of it is easy to see. For instance, if $a(E)^2=4q$ (so $q$ is a
square) the Frobenius
$\frob$ is a solution of the quadratic equation
$X^2-a(E)X+q=0$, which has the double root $\pm\sqrt{q}\in \Zz$ (with sign
chosen has in the statement of the proposition). So $\frob-1\in
\Zz\subset \End(E)$, and Lemma~\ref{lm-cong-end}
implies~(\ref{eq-supersing-big}). 
\par
For the other cases, it is known that $a(E)^2=q$, $2q$ or $3q$, or
$a(E)=0$. If $a(E)=0$ (the only possibility over $\Fp_p$), for
instance, the congruence $a(E)\equiv 2\mods{\dun(E)}$
proves~(\ref{eq-supersing-small}). Similarly in the other cases the
congruences of Lemma~\ref{lm-cong-end} either
prove~(\ref{eq-supersing-small}), or a weaker bound like
$\dun(E)\leq 4$, which will suffice here (see~\cite[Lemma
4.8]{schoof-2} for complete details).
\end{proof}

This has the following global corollary which shows that supersingular
primes have a small contribution to~(\ref{eq-sum-dun}).

\begin{corollary}
Let $E/\Qq$ be an elliptic curve. We have
\begin{align}
\sum_{\stacksum{p\leq X}{a_p(E)=0}}{\dun(E)}&\ll_E 
\frac{X}{\log X}\quad\text{ if $E$ has CM}\\
\sum_{\stacksum{p\leq X}{a_p(E)=0}}{\dun(E)}&\ll_E X^{3/4}
\quad\text{otherwise}.
\end{align}
for all $X\geq 2$, the implied constant depending on $E$ only.
\end{corollary}

\begin{proof}
If $E$ has CM, the number of supersingular primes $p\leq X$ is well
known to be (see e.g.~\cite{lang-trotter}) $\sim X/(2\log X)$, while
if $E$ doesn't have CM, Elkies~\cite{elkies} has shown that the number
of supersingular primes 
$p\leq X$ of $E$ is $\ll_E X^{3/4}$. Since $\dun(p)\leq 2$ by
Proposition~\ref{pr-supersing}, the result follows.
\end{proof}

\begin{remark}
In the non-CM case, Serre's proof~\cite{serre-mg} that the number of
supersingular primes $\leq X$ is $o(X/\log X)$ suffices to show that 
$$
\sum_{\stacksum{p\leq X}{a_p(E)=0}}{\dun(E)}
=o\Bigl(\frac{X}{\log X}\Bigr)
$$
as $X\ra +\infty$.
\end{remark}

From now on we assume that $E/\Fp_q$ is an ordinary elliptic curve
over a finite field with $q$ elements. We let $\Oc=\End(E)$, $K$ the
field of fraction of $\Oc$, $\Oc_K$ the ring of integers of $K$. Let
$\frob\in \Oc$ be the Frobenius 
endomorphism of $E$. The main tool to find $\dun(E)$ is
Lemma~\ref{lm-cong-end}.

\begin{lemma}\label{lm-div-ok}
Let $d\geq 1$ be an integer. We have $d\mid \frob-1$ in $\Oc_K$ if and
only if $a(E)\equiv 2\mods{d}$ and $n(E)=q+1-a(E)=0\mods{d^2}$.
\end{lemma}

\begin{proof}
Let $\frob'=(\frob-1)/d\in K$, so $d\mid \frob-1$ in $\Oc_K$ if and
only if $\frob'\in \Oc_K$. But since $\frob'\not\in \Zz$, since $E$ is
ordinary, its minimal polynomial over $\Zz$ is
$$
(X-\frob')(X-\overline{\frob'})=X^2-\frac{a(E)-2}{d}X+\frac{n(E)}{d^2}.
$$
Hence the result since $\Oc_K$ is the integral closure
of $\Zz$ in $K$.
\end{proof}

We can check that this gives back the other congruences.

\begin{lemma}
Let $a$, $q\geq 2$, $d\geq 1$ be integers such that
$$
\begin{cases}
a\equiv 2\mods{d}\\
q+1-a\equiv 0 \mods{d^2}.
\end{cases}
$$
Then $q\equiv 1\mods{d}$ and $a^2-4q\equiv 0\mods{d^2}$.
\end{lemma}

\begin{proof}
We have modulo $d$
$$
0=q+1-a=q+1-2=q-1
$$
and modulo $d^2$
$$
a^2-4q=(q+1)^2-4q=(q-1)^2=0.
$$
\end{proof}

\begin{lemma}
Let $E/\Fp_q$ as before. We have $d\mid \frob-1$ in $\Oc_K$ if and
only if $a^2-4q=0\mods{d^2}$ and $n(E)\equiv 0\mods{d^2}$.
\end{lemma}

\begin{proof}
Let again $\frob'=(\frob-1)/d\in K$. In terms of $\frob'$, the two
assumptions are
\begin{align*}
N(\frob-1)&=d^2N\frob'=0\mods{d^2}\\
(\frob-\overline{\frob})^2&=d^2(\frob'-\overline{\frob'})^2
=0\mods{d^2},
\end{align*}
hence we see that $N\frob'\in \Zz$ and
$(\frob'-\overline{\frob'})^2\in \Zz$.
\par
The latter is also $(\frob'+\overline{\frob'})^2-4N\frob'$, hence we
deduce that $(\Tr \frob')^2\in \Zz$. Since $\Tr(\sigma')\in \Qq$,
it must be an integer, hence the result.
\end{proof}

Those easy results give a good handle on the condition $d\mid \frob-1$
in $\Oc_K$. The problem is that $\Oc$ is in general a proper order in
$\Oc_K$. However, the necessary congruence conditions are also
sufficient ``up to isogeny''.

\begin{proposition}\label{pr-dun-iso}
Let $\Fp_q$ be a finite field with $q$ elements, $d\geq 1$ an integer
coprime with $q$.
\par
There exists an ordinary elliptic curve $E/\Fp_q$ with $E[d]\subset
E(\Fp_q)$, i.e. $d\mid \dun(E)$, if and only if there exists $a\in\Zz$
such that
$$
\begin{cases}
|a|<2\sqrt{q}\\
(a,q)=1\\
a\equiv 2\mods{d}\\
q+1-a\equiv 0\mods{d^2}.
\end{cases}
$$
\end{proposition}

For the proof we need some results which are part of Honda-Tate theory
for elliptic curves (which goes back to Deuring), and others due to
Waterhouse~\cite{water} concerning the endomorphism rings of elliptic
curves over finite fields.

\begin{theorem}\label{th-honda-tate}(Deuring, Honda, Tate)
Let $\Fp_q$ be a finite field with $q$ elements. Given an integer $a$
such that $|a|<2\sqrt{q}$ and $(a,q)=1$, there exists an ordinary
elliptic curve $E$ over $\Fp_q$ with $a(E)=a$.
\end{theorem}

See for instance~\cite[Th. 4.1]{water}.

\begin{theorem}\label{th-endos}(Deuring, Waterhouse)
Let $\Fp_q$ be a finite field with $q$ elements, $a$ an integer with
$|a|<2\sqrt{q}$ and $(a,q)=1$. Let $K=\Qq(\sqrt{a^2-4q})$ and let
$\Oc$ be an order of $K$. There exists an ordinary elliptic curve
$E/\Fp_q$ with $a(E)=a$ and $\End(E)=\Oc$ if and only if $\Oc$
contains the roots of 
$$
X^2-aX+q=0.
$$
\end{theorem}

See~\cite[Th. 4.2 (2)]{water}. Note that this second result requires
Tate's Theorem identifying relating isogenies between elliptic curves 
with Galois-invariant maps between their $\ell$-adic Tate modules,
$(\ell,q)=1$.

\begin{proof}[Proof of Proposition~\ref{pr-dun-iso}]
The condition is necessary. Conversely, if $a$ exists as
described, Theorem~\ref{th-endos} shows that there exists $E/\Fp_q$
with $a(E)=a$ and $\End(E)=\Oc_K$, where $K$ is the imaginary
quadratic field $K=\Qq(\sqrt{a^2-4q})$.
\par
The congruence conditions on $a(E)$
and $n(E)$ then mean (Lemma~\ref{lm-div-ok}) that $d\mid \frob-1$ in
$\Oc_K=\End(E)$, hence $d\mid \dun(E)$.
\end{proof}

\begin{remark}
If $q=p\geq 5$ is prime, one can remove the condition $(a,p)=1$ on $a$
from the statement of the proposition. Indeed, if $p\mid a$, we have
$a=0$, and since $a\equiv 
2\mods{d}$, the only values of $d$ occurring are $d=1$ and $d=2$. But
those can be obtained from ordinary elliptic curves: $d=1$ by any $E$,
and $d=2$ by a Legendre curve
$$
E_{\lambda}\,:\,y^2=x(x-1)(x-\lambda)
$$
(which always has $2\mid \dun(E_{\lambda})$)
for some $\lambda\in \Fp_p-\{0,1\}$. Indeed, the condition that
$E_{\lambda}$ be 
ordinary is equivalent (see e.g.~\cite[V-4]{silv}) to
$H_p(\lambda)\not=0$, where $H_p$ is the Hasse-Deuring polynomial
$$
H_p=\sum_{j=0}^{(p-1)/2}{\binom{(p-1)/2}{j}^2X^j}\in \Fp_p[X].
$$
Since $0\leq \deg H_p=(p-1)/2<p-2$, there is a $\lambda\in
\Fp_p-\{0,1\}$ which is not a root of $H_p$, hence a corresponding
ordinary $E_{\lambda}$ with $2\mid \dun(E_{\lambda})$.
\par
On the other hand, if $q$ is a square, let $d=\sqrt{q}+1$. Then $d$
satisfies all the assumptions of Proposition~\ref{pr-dun-iso} with
$a=-2\sqrt{q}$, except $(a,q)=1$. But this is the only value of $a(E)$
for which one could have $d\mid \dun(E)$, and it corresponds to
supersingular curves, so that in general $(a,q)=1$ is a necessary
assumption.
\end{remark}

In applications, we are interested in the invariant $\dun(E)$, and
$d=\dun(E)$ means not only $E[d]\subset E(\Fp_q)$, but also that no
larger $d$ (coprime with $q$) satisfies this. However,
Proposition~\ref{pr-dun-iso} remains true with $d=\dun(E)$ instead of
$d\mid \dun(E)$ in the conclusion.

\begin{proposition}\label{pr-dun-ind}
Let $E/\Fp_q$ be an elliptic curve over a finite field with
$\dun(E)=d$. For every $\delta\mid d$, there exists an elliptic curve
$E'/\Fp_q$ which is $\Fp_q$-isogenous to $E$ and satisfies
$\dun(E')=\delta$.
\end{proposition}

\begin{corollary}\label{cor-dun-iso}
Let $\Fp_q$ and $d\geq 1$, $(d,q)=1$, be as above.
There exists an ordinary elliptic curve $E/\Fp_q$ with $\dun(E)=d$ if
and only if there exists $a\in\Zz$ 
with $|a|<2\sqrt{q}$, $(a,q)=1$, and such that 
$$
\begin{cases}
a\equiv 2\mods{d}\\
q+1-a\equiv 0\mods{d^2}.
\end{cases}
$$
\end{corollary}

\begin{proof}[Proof of the proposition]
Write $d=\delta\delta'$. We have, with the same notation as usual,
$\frob'=(\frob-1)/d\in \Oc$. It suffices to find a smaller order
$\Oc'\subset \Oc$ with $\delta'\frob'=(\frob-1)/\delta\in \Oc'$ but
for which there is no $e>1$ with $\delta'\frob'/e\in\Oc'$. Then, since
$\frob\in\Oc'$, Theorem~\ref{th-endos} shows that there exists
$E'/\Fp_q$, isogenous to $E$ (hence ordinary), with
$\End(E')=\Oc'$. Then $\dun(E')=\delta$ by construction
(Lemma~\ref{lm-cong-end}).
\par
To construct $\Oc'$, we write $\Oc=\Zz\oplus \omega\Zz$
(see~\cite[7-A]{cox}), and 
correspondingly $\frob'=m+n\omega$, for some $m$, $n\in \Zz$. So 
$\delta'\frob'=\delta'm+\delta'n\omega$. Let $\Oc'$ be the order $\Zz\oplus
cn\omega\Zz$ of $K$. Then $\delta'\frob'\in \Oc'$, but for any $e\geq
1$, we have
$$
\frac{\delta'\frob'}{e}=\frac{m\delta'}{e}+\frac{\delta' n\omega}{e}
$$
and for this to be in $\Oc'$ we must have $e=1$, showing that
$\Oc'$ satisfies the conditions required.
\end{proof}

\begin{remark}\label{rm-critere}
Over the base field $\Fp_p$, it is again possible to remove the
condition $(a,p)=1$. Putting back supersingular curves, the following
statement holds:
\par
Let $d\geq 1$ be an integer. There exists an elliptic curve $E$ over
$\Fp_p$ with $\dun(E)=d$ if and only if there exists $a$,
$|a|<2\sqrt{p}$, such that
$$
\begin{cases}
a\equiv 2\mods{d}\\
p+1-a\equiv 0\mods{d^2}.
\end{cases}
$$
In particular, this is always true for $d=1$ and $d=2$ (the latter for
$p\geq 5$).
\end{remark}

\begin{remark}
As a side remark and pretext to mention another interesting problem of
analytic number theory, the case $d=1$ can be studied purely
analytically from Theorem~\ref{th-honda-tate} and the distribution of
squarefree numbers in short intervals. Indeed, if $a\not=0$ is such that
$p+1-a$ is squarefree, any elliptic curve $E/\Fp_p$ with $a(E)=a$ must
have $\dun(E)=1$. Hence the existence of such an $E$ (for $p$ large
enough only, however) follows from any
``non-trivial'' estimate for error term in the asymptotic formula for
the number of squarefree numbers $n\leq X$
$$
\sum_{n\leq X}{|\mu(n)|}=\frac{1}{\zeta(2)}X+O(X^{\theta})
$$
as $X\ra +\infty$, with $\theta<1/2$, since this implies in particular
$$
\sum_{|a|<2\sqrt{p}}{|\mu(p+1-a)|}>0.
$$
The value $\theta=1/2$ is easily obtained, any improvement requiring
non-trivial cancellation in some exponential sums. See for
instance~\cite[p. 46]{graham-kolesnik} where it is shown that
$\theta=4/9+\eps$ is possible, for any $\eps>0$.
\end{remark}

\subsection{Results using the trace formula}

The criterion obtained in Corollary~\ref{cor-dun-iso} is quite
convenient. However, from our point of view, it is more natural to fix
a prime $p$ (or prime power) and look for which $d\mid p-1$ there exists
$E/\Fp_p$ with $\dun(E)=d$.
\par
A criterion of that type arises naturally if we use, instead of
endomorphism rings, the theory of modular curves and the
Eichler-Selberg trace formula. Although Corollary~\ref{cor-dun-iso}
and Remark~\ref{rm-critere} would suffice for the applications in the
next section, this approach is sufficiently independent and
instructive to be included here.


\begin{theorem}\label{th-criterion}
Let $p$ be a prime number, $d\mid p-1$ an
integer. Write $d=ef$ with $(e,2)=1$, $f\mid 2^{\infty}$. There exists
$E/\Fp_p$ with $\dun(E)=d$ if and only if there 
exists $a$ with $|a|<2\sqrt{p}$ such that
\par
\emph{(1)} We have $e^2\mid a^2-4p$;
\par
\emph{(2)} If $f\not=1$, there exists $\eps=\pm 1\mods{f}$
such that $\eps^2-a\eps+p=0\mods{f^2}$.
\end{theorem}

We need some geometric preliminaries.
For any integer $d\geq 1$, there exists a smooth
affine curve $Y(d)$ naturally defined over $\Qq(\mmu_d)$, with good
reduction at all primes $p\nmid d$, which is a coarse moduli scheme
for ``elliptic curves with a $d$-level structure'' (see~\cite{katz-mazur}
or~\cite{de-ra}). Over $\Cc$, $Y(d)(\Cc)$ is the ``usual'' quotient 
$$
\Gamma(d)\backslash \mathbf{H}
$$
of the upper half-plane by the principal congruence subgroup
$$
\Gamma(d)=\{ g\in SL(2,\Zz)\,\mid\, g\equiv 1\mods{d}\}.
$$
\par
Moreover, $Y(d)$ has an integral model over the ring of integers
$\Zz[\mmu_d]$ of the cyclotomic field $\Qq(\mmu_d)$. Notice that
$p\equiv 1\mods{d}$ means that the $p$ is totally split in this field,
hence $\Zz[\mmu_d]/(p)\simeq (\Fp_p)^{\varphi(d)}$. The above ``moduli
scheme'' sentence implies in particular (see~\cite[VI-3]{de-ra}) that
for $p\equiv 1\mods{d}$,
it is the same to give a point in $Y(d)(\Fp_p)$ as to give a pair
$(E,(e_1,e_2))$ of an elliptic curve $E/\Fp_p$ together with
two $\Fp_p$-rational points of order $d$, $e_1$ and $e_2$, such that the
Weil pairing $\gena e_1,e_2 \genb$ is equal to a fixed primitive $d$-th
root of unity (these pairs taken up to isomorphism). In other words we
have (see also~\cite{howe} for a description of other modular curves
over finite fields): 

\begin{lemma}
Let $p$ be a prime number, $d\geq 1$ an integer such that $d\mid
p-1$. Then there exists $E/\Fp_p$ with $E[d]\subset E(\Fp_p)$ if and
only if $Y(d)(\Fp_p)\not= \emptyset$.
\end{lemma}

We are thus reduced to finding points on the curve $Y(d)$ over the
finite field $\Fp_p$.
\par
The curve $Y(d)$ has a natural compactification $X(d)$, which over
$\Cc$ amounts to adding the cusps to $\mathbf{H}$ before taking the
quotient by $\Gamma(d)$. The projective curve $X(d)$ has also good
reduction at all $p$ not dividing $d$ (and a moduli description in
terms of ``generalized elliptic curves''). 
\par
For $p$ a prime of good reduction, the local zeta function of $X(d)$
$$
Z(X(d),p)=\exp\Bigl(\sum_{n\geq 1}
{\frac{|X(d)(\Fp_{p^n})|}{n}T^n}\Bigr)
$$
is, by general results (due to F.K. Schmidt in this case of curves over
finite fields), a rational function of the 
form
$$
Z(X(d),p)=\frac{P_d(T)}{(1-T)(1-pT)},
$$
where $P_d$ is a polynomial of degree $2g(d)$, $g(d)$ being the genus
of $X(d)$. From this and the definition of $Z(X(d),p)$, one can deduce
immediately that 
$$
|X(d)(\Fp_p)|=p+1-\sum_{i}{\alpha_i}
$$
where
$$
P_d=\prod_{i=1}^{2g(d)}{(1-\alpha_i T)}.
$$
The point of using the compactified curve $X(d)$ is that we have the
following consequence of the computation of the zeta functions of
modular curves by Shimura (\cite[\S 7.5]{shimura-1}).

\begin{theorem}\label{th-shimura}
Let $d\geq 1$ be an integer, $p\equiv 1\mods{d}$ a prime number. We
have
$$
|X(d)(\Fp_p)|=p+1-\Tr (T_p|S_2(\Gamma(d))),
$$
where the last term is the trace of the Hecke operator $T_p$ acting on
the space $S_2(\Gamma(d))$ 
of weight $2$ holomorphic cusp forms for the congruence subgroup
$\Gamma(d)$.
\end{theorem}

More
precisely, Shimura's result gives the zeta function for models of
$X(d)$ over $\Qq$, of which there exist several; but all give the
same $X(d)$ over $\Qq(\mmu_d)$, hence the result since we consider $p$
totally split in $\Qq(\mmu_d)$.
\par
The Eichler-Selberg trace formula gives an expression for the
trace, which one may use to find when $X(d)(\Fp_p)\not=\emptyset$;
this idea is used by Jordan~\cite{jordan}. However, he works
with Shimura curves, which are compact, and his main interest
is at primes of bad reduction.
\par
Here we have to take the cusps into account, since they do not
correspond to elliptic curves. Over $\Cc$, the cusps of $X(d)$ are
described in~\cite[Lemma 1.42]{shimura-1}. We need to know which are
rational over $\Fp_p$.
\par
Let $\varphi^+(d)$ denote the number of even Dirichlet characters
modulo $d$ (i.e. $\chi(-1)=1$). This is given by
\begin{equation}\label{eq-phiplus}
\varphi^+(d)=
\begin{cases}
{\displaystyle{\frac{\varphi(d)}{2}}}&\text{if $d>2$}\\
\varphi(d)=1&\text{if $d=2$}.
\end{cases}
\end{equation}
By orthogonality of characters, we have for any $x\in \Zz$
\begin{equation}\label{eq-ortho}
\sum_{\chi\text{ even}}{\chi(x)}=
\begin{cases}
\varphi^+(d)&\text{if $x\equiv \pm 1\mods{d}$}\\
0&\text{otherwise.}
\end{cases}
\end{equation}
(This will be needed later on).

\begin{lemma}\label{lm-count-cusps}
Let $d\geq 1$ be an integer, $p\equiv 1\mods{d}$ a prime number. 
All the cusps of $X(d)$ are $\Fp_p$-rational, and in particular
$$
|(X(d)-Y(d))(\Fp_p)|=\varphi^+(d)\psi(d)=
\begin{cases}
\frac{1}{2}\varphi(d)\psi(d)&\text{ if $d>2$}\\
\varphi(d)\psi(d)&\text{ if $d=2$}.
\end{cases}
$$
\end{lemma}

\begin{proof}
This follows from Theorem 10.9.1 (3) of~\cite{katz-mazur} which says
(in particular) that the cusps of $Y(d)/\Zz[\mmu_d]$ are rational over
$\Zz[\mmu_d]$, and ``do not vary'' by base change to any
$\Zz[\mmu_d]$-algebra; heuristically, 
cusps rational over $\Fp_p$ ``correspond'' to level $d$
structures on the Tate curve $Tate(q)/\Fp_p((q))$ rational over
$\Fp_p((q))$. Since the 
$d$-torsion of the latter is isomorphic as a Galois module
(\cite[V-3]{silv-2}) to 
$$
\Zz/d\Zz \times \boldsymbol{\mu}_d
$$
and $d\mid p-1$ so $\mmu_d\subset \Fp_p^{\times}$, it is visible that 
all level $d$ structures on $Tate(q)$ are $\Fp_p((q))$-rational.
\par
The number of cusps over $\Cc$ is found in~\cite[p. 22]{shimura-1},
or can be recomputed directly from the result in~\cite{katz-mazur}
quoted above.
\end{proof}

We will now state the trace formula in the form needed. A 
paper by Fomenko~\cite{fomenko} should include it, but I have
not been able to see it. On the other hand, the trace formula for
$\Gamma(d)$ is not easily derived from general accounts: for instance,
it does not correspond to an ``Eichler order'', so the arguments
in~\cite[Ch. 6]{miyake}, for instance, can not be adapted
straightforwardly. We can 
circumvent these difficulties by reducing to the much better known
case of Hecke congruence subgroups $\Gamma_0(N)$, for which we can
quote for instance~\cite{miyake},~\cite{serre-tp} or~\cite{hamer}
(among many other non-conflicting sources).

\begin{lemma}
Let $d\geq 1$ an integer and $p\equiv 1\mods{d}$ a prime number. There
exists an isomorphism of vector spaces
$$
u\,:\,
S_2(\Gamma(d))\lra \bigoplus_{\chi\text{ even}}
{S_2(\Gamma_0(d^2),\chi)} 
$$
where the direct sum is over all even Dirichlet characters modulo $d$,
$S_2(\Gamma_0(d^2),\chi)$ is the space of weight $2$ cusp forms for
$\Gamma_0(d^2)$ with nebentypus $\chi$, which satisfies
$$
u\circ T_p=T_p \circ u,
$$
where on the right $T_p$ is the direct sum of Hecke operators acting
on $S_2(\Gamma_0(d^2),\chi)$. 
\end{lemma}

\begin{proof}
We first introduce the congruence subgroups
$$
\Gamma_0(d,d)=
\{
g=
\begin{pmatrix}
a&b\\
c&d
\end{pmatrix}
\in SL(2,\Zz)\,\mid\, b=c=0\mods{d}
\}.
$$
We have $\Gamma(d)\lhd \Gamma_0(d,d)$ with quotient
$(\Zz/d\Zz)^{\times}$. 
\par
The even Dirichlet characters modulo $d$ are extended to characters of
$\Gamma_0(d,d)$ by
$$
\chi(g)=\chi(d).
$$
Then the natural action of $\Gamma_0(d,d)/\Gamma(d)$ on
$S_2(\Gamma(d))$ gives
the direct sum decomposition
$$
S_2(\Gamma(d))=\bigoplus_{\chi\text{ even}}{S_2(\Gamma_0(d,d),\chi)}
$$
(for odd $\chi$, $S_2(\Gamma_0(d,d),\chi)=0$).
\par
Since $d\mid p-1$, we have $\chi(p)=1$ for any character modulo $d$,
and this implies that $T_p$ acting on $S_2(\Gamma(d))$ is the direct
sum of the $T_p$ acting on $S_2(\Gamma_0(d,d),\chi)$
(see~\cite[3.5.6]{shimura-1}; it amounts to the fact that a $\chi(p)$
appears in the explicit formula for $T_p$ acting on
$S_2(\Gamma_0(d,d),\chi)$ but not for $T_p$ on $S_2(\Gamma(d))$).
\par
Moreover $\Gamma_0(d,d)$ is conjugate to $\Gamma_0(d^2)$ in
$SL(2,\Rr)$ by
$$
g\mapsto 
\begin{pmatrix}
d^{-1/2}&0\\
0&d^{1/2}
\end{pmatrix}
g
\begin{pmatrix}
d^{1/2}&0\\
0&d^{-1/2}
\end{pmatrix}
$$
This induces an isomorphism
\begin{equation}\label{eq-iso-chi}
S_2(\Gamma_0(d,d),\chi)\ra S_2(\Gamma_0(d^2),\chi)
\end{equation}
given by
$$
f\mapsto f\mid_2
\begin{pmatrix}
d^{1/2}&0\\
0&d^{-1/2}
\end{pmatrix},
$$
(where $\cdot\mid_2\cdot $ denotes the usual weight $2$ action of
$SL(2,\Zz)$ on functions).
Hence we have an isomorphism
$$
u\,:\,S_2(\Gamma(d))\ra
\bigoplus_{\chi\text{ even}}{S_2(\Gamma_0(d^2),\chi)}.
$$
Since $T_p$ commutes with the isomorphism~(\ref{eq-iso-chi}), as is
well known (compare~\cite[4.6.1]{miyake}), $u$ is also compatible.
\end{proof}

\begin{corollary}
Let $d\geq 1$ be an integer, $p\equiv 1\mods{d}$ a prime number. 
We have
$$
\Tr T_p|S_2(\Gamma(d))=\sum_{\chi\text{ even}}{\Tr
T_p|S_2(\Gamma_0(d^2),\chi)}, 
$$
where the sum is over even Dirichlet characters modulo $d$.
\end{corollary}

To state the trace formula for $S_2(\Gamma_0(d^2),\chi)$,
we require some further notation. Recall $\chi$ is an even character.
\par
If $\Oc$ is an order in an imaginary quadratic field, we let
$H(\Oc)$ denote its class number, divided by half the number of
units (i.e. $1$ unless $\Oc=\Zz[i]$, where it's $2$, or
$\Oc=\Zz[\mmu_3]$, where it's $3$). We denote by $\Oc(\delta)$ the
order with discriminant $\delta<0$, and let
$H(\delta)=H(\Oc(\delta))$. 
\par
If $\Oc\subset \Oc(a^2-4p)$ is a sub-order with index $f$, and $N\geq
1$, we denote
\begin{equation}\label{eq-muc}
\mu_{\chi}(\Oc,a,p,N)=
\frac{\psi(N)}{\psi(N/(N,f))}
\sum_{\stacksum{x\mods{N}}
{x^2-ax+p=0\mods{N(N,f)}}}{\chi(x)}
\end{equation}
(it makes sense).

\begin{theorem}
Let $d\geq 1$ be an integer, $\chi$ an even Dirichlet character modulo
$d$ and $p\equiv 1\mods{d}$ a prime number. We
have
$$
\Tr T_p|S_2(\Gamma_0(d^2),\chi)=
t_d(\chi)-t_e(\chi) -t_h(\chi)
$$
where
\begin{align}
t_d(\chi)&=
\begin{cases}
p+1&\text{ if $\chi=1$}\\
0&\text{ otherwise}
\end{cases}\nonumber\\
t_e(\chi)&=\frac{1}{2}
\sum_{\stacksum{a\in\Zz}{a^2<4p}}{
\sum_{\Oc\subset \Oc(a^2-4p)}{
H(\Oc)\mu_{\chi}(\Oc,a,p,d^2)}}\label{eq-ell-term}\\
t_h(\chi)&=
\frac{1}{2}\sum_{b\mid p}
{\sum_{c\mid d^2}
{\varphi\Bigl(\Bigl(
\frac{d^2}{c},c\Bigr)\Bigr)\chi(y_c)}},
\label{eq-hyp-term}
\end{align}
where $y_c$ is an integer modulo $d^2/((d^2/c,c))$ such that
\begin{align*}
y_c&\equiv b\mods{c}\\
y_c&\equiv p/b\mods{d^2/c}.
\end{align*}
\end{theorem}

\begin{remark}
The notation follows the genesis of these terms, for example in
Shimura's formulation~\cite{shimura-2} of the trace formula as a
kind of Lefschetz formula for correspondences: $t_d$ refers to the
``dual term'', as it should be understood as coming from an $H^2$,
which is non-zero only for weight $2$ and trivial character; $t_e$ refers
to the contribution of elliptic elements, and $t_h$ to the
contribution of hyperbolic elements. There is no parabolic
contribution here because we are working with $T_p$ and $p$ is not a
square.
\end{remark}

\begin{proof}
Serre~\cite[4.1]{serre-tp} quotes a general formula for all
levels and characters. To deduce the form claimed, notice that the term
denoted $A_1$ vanishes since $p$ is not a square and the term denoted
$A_4$ gives directly $t_d$. The term $-A_3$, we claim, is the same as
$t_h$. Indeed, we have from loc. cit.
$$
-A_3=\frac{1}{2}\sum_{b\mid p}{\mathrm{Inf}(b,p/b)
\sum_{c}{\varphi\Bigl(\Bigl(
\frac{d^2}{c},c\Bigr)\Bigr)\chi(y_c)}}
$$
where the sum over $c$ is restricted to divisors of $d^2$ such that
\begin{align}
\Bigl(\frac{d^2}{c},c\Bigr)&\mid \frac{p}{b}-b\label{eq-1}\\
\Bigl(\frac{d^2}{c},c\Bigr)&\mid \frac{d^2}{d_{\chi}}\label{eq-2}
\end{align}
($d_{\chi}$ is the conductor of $\chi$). Now first for $b\mid p$ we
have $\mathrm{Inf}(b,p/b)=1$, and also $p/b-b=\pm(p-1)$.
Also, for all $c\mid d^2$, we have
\begin{equation}\label{eq-utile}
\Bigl(\frac{d^2}{c},c\Bigr)\mid d.
\end{equation}
Indeed, proceeding locally at each prime $\ell$, if $d=\ell^{\nu}$,
and $c=\ell^{\mu}$ with $\mu\leq 2\nu$, the exponent of $(d^2/c,c)$ is
$\mathrm{Inf}(\mu,2\nu-\mu)\leq \nu$.
\par
Since $d_{\chi}\mid d$, and $d\mid p-1$, this shows that the two
restrictions~(\ref{eq-1}) and~(\ref{eq-2}) on $c$ are satisfied for
all $c\mid d^2$.
\par
Similarly, the term $-A_2$ in loc. cit. is the same as $t_e$ (recall
the weight is $2$). 
\end{proof}

\begin{corollary}
Let $d\geq 1$ an integer, $p\equiv 1\mods{d}$ a prime number. We have
$$
\Tr T_p|S_2(\Gamma(d))= p+1-t_e-t_h
$$
where
\begin{align*}
t_e&=\sum_{\chi\text{ even}}{t_e(\chi)}\\
t_h&=\sum_{\chi\text{ even}}{t_h(\chi)}.
\end{align*}
\end{corollary}

The next observation is elementary but crucial.

\begin{proposition}\label{pr-cusps}
Let $d\geq 1$ an integer, $p\equiv 1\mods{d}$ a prime number. Then
$t_h$ is equal to the number of $\Fp_p$-rational cusps of $X(d)$,
i.e. $\varphi^+(d)\psi(d)$.
\end{proposition}

\begin{proof}
The point is that the integer $y_c$ in~(\ref{eq-hyp-term}) can be
chosen, for $b\mid p$ and any $c\mid d^2$, to satisfy
\begin{equation}\label{eq-yc}
y_c\equiv 1\mods{d},
\end{equation}
and since the character $\chi$ is modulo $d$ (not $d^2$), we have
$\chi(y_c)=1$ for any $\chi$.
\par
To see~(\ref{eq-yc}), we work locally at all primes $\ell$ as before.
We have $b=1$ or $b=p$: both situations are similar, so assume $b=1$. 
Then writing $\ell^{\nu}$ for the $\ell$-component of $d$,
$\ell^{\mu}$ for that of $c$, $\mu\leq 2\nu$, the conditions on $y_c$
are
$$
\begin{cases}
y_c\equiv 1&\mods{\ell^{\mu}}\\
y_c\equiv p&\mods{\ell^{2\nu-\mu}}.
\end{cases}
$$
We have either $\mu\geq \nu$, in which case the first equation implies
$y_c\equiv 1\mods{\ell^{\nu}}$, or $2\nu-\mu>\nu$, in which case the
second implies $y_c\equiv p\equiv 1\mods{\ell^{\nu}}$, since $p\equiv
1\mods{d}$. Those local congruences patch, proving the claim for
$b=1$, and $b=p$ is symmetric.
\par
Using~(\ref{eq-ortho}), and the fact that $b=1$ and $b=p$ have the
same contribution, we can now write $t_h$ as
$$
t_h=\varphi^+(d)\sum_{c\mid d^2}{\varphi\Bigl(
\Bigl(\frac{d^2}{c},c\Bigr)\Bigr)}.
$$
We now use $\delta=(d^2/c,c)$ as new summation variable. Recall that
$\delta\mid d$ (\ref{eq-utile}). We get
\begin{equation}\label{eq-interm}
t_h=\varphi^+(d)\sum_{\delta\mid d}{\varphi(\delta)
M(\delta)}
\end{equation}
where
$$
M(\delta)=|\{ c\mid d^2\,\mid\, (d^2/c,c)=\delta\}|.
$$
We work again at each prime $\ell$ separately, with $\ell^{\nu}$ the
component of $d$, $\ell^{\rho}$ that of $\delta$. The
$\ell$-component $\ell^{\mu}$ of $c$ must therefore satisfy
$$
\mathrm{Inf}(\mu,2\nu-\mu)=\rho.
$$
Given $\rho$, there are two choices of $\mu$, namely $\mu=\rho$ or
$\mu=2\nu-\rho$ (since $\rho\leq \nu$), \emph{unless} $\rho=\nu$,
since in this case they coincide.
\par
It is clear that
$$
f(d)=\sum_{\delta\mid d}{\varphi(\delta)
M(\delta)}
$$
is multiplicative. Now we compute the value at $\ell^{\nu}$ using the
above: we have
$$
\begin{cases}
M(\ell^{\rho})=2&\text{ if $\rho<\nu$}\\
M(\ell^{\nu})=1,
\end{cases}
$$
so
\begin{align*}
f(\ell^{\nu})&=\sum_{\rho=0}^{\nu-1}
{2\varphi(\ell^{\rho})}+\varphi(\ell^{\nu})\\
&=2\ell^{\nu-1}+\ell^{\nu}-\ell^{\nu-1}\\
&=\psi(\ell^{\nu}).
\end{align*}
Comparing this and~(\ref{eq-interm}) with Lemma~\ref{lm-count-cusps},
the proposition is proved.
\end{proof}

\begin{remark}
I did not find mention in the literature of this fact that the
hyperbolic terms in the trace formula ``count the cusps'', although
that must be well-known. This applies obviously to more general
subgroups, with corresponding applications to elliptic curves over
finite fields using their moduli interpretation. It would be
interesting to see if there are higher-rank analogues, and their
consequences. 
\end{remark}

\begin{corollary}
Let $d\geq 1$ an integer, $p\equiv 1\mods{d}$ a prime number. We have
$$
|Y(d)(\Fp_p)|=t_e.
$$
In particular, there exists an elliptic curve $E/\Fp_p$ with
$E[d]\subset E(\Fp_p)$ if and only if $t_e>0$.
\end{corollary}

\begin{proof}
We have by Proposition~\ref{pr-cusps}
\begin{align*}
|Y(d)(\Fp_p)|&=|X(d)(\Fp_p)|-t_h\\
&=p+1-\Tr T_p-t_h\text{ (by Theorem~\ref{th-shimura})}\\
&=p+1-(p+1-t_e-t_h)-t_h\\
&=t_e.
\end{align*}
\end{proof}

Because of the average over $\chi$, $t_e$ is a sum of
terms each of which is obviously $\geq 0$. This makes it possible to
find a criterion to have $Y(d)(\Fp_p)\not= 
\emptyset$ (compare~\cite{jordan}). If the formula for $t_e$ involved
any oscillatory sum, it would be much harder to exploit it.
\par
For a quadratic imaginary order $\Oc\subset \Oc(a^2-4p)$
with index $f$ we let
$$
\mu(\Oc,a,p,d^2)=\sum_{\chi\text{ even}}{\mu_{\chi}(\Oc,a,p,d^2)},
$$
and
$$
\mu(a,p,d^2)=\mu(\Oc(a^2-4p),a,p,d^2).
$$
From~(\ref{eq-ell-term}) we have
\begin{equation}\label{eq-te-final}
t_e=\frac{1}{2}\sum_{\stacksum{a\in\Zz}{a^2<4p}}
{\sum_{\Oc\subset \Oc(a^2-4p)}{
H(\Oc)\mu(\Oc,a,p,d^2)}}.
\end{equation}

\begin{lemma}\label{lm-mu0}
We have for an order $\Oc\subset \Oc(a^2-4p)$ of index $f$
$$
\mu(\Oc,a,p,d^2)=\frac{\varphi^+(d)\psi(d^2)}{\psi(d^2/(d^2,f))}
\mu_0(a,p,d,f)
$$
where
$$\mu_0(a,p,d,f)=
|\{ x\mods{d^2}\,\mid\,
x=\pm 1\mods{d}\text{ and }
x^2-ax+p=0\mods{d^2(d^2,f)}.\}|
$$
\end{lemma}

This is simply the orthogonality relation~(\ref{eq-ortho}).
We let $\mu_0(a,p,d)=\mu_0(a,p,d,1)$.

\begin{corollary}\label{cor-critere}
Let $d\geq 1$ be an integer and $p\equiv 1\mods{d}$ a prime number. We
have
$Y(d)(\Fp_p)\not=\emptyset$ if and only if there exists an integer $a$
with $|a|<2\sqrt{p}$ such that $\mu(a,p,d)>0$, if and only if there
exists $a$ with $|a|<2\sqrt{p}$ such that the equation
$x^2-ax+p=0\mods{d^2}$ has a solution $x$ with $x=\pm 1\mods{d}$.
\end{corollary}

\begin{proof}
From~(\ref{eq-te-final}), we have $t_e>0$ if and only if there exists
$a$ and $\Oc\subset \Oc(a^2-4p)$ with $\mu(\Oc,a,p,d^2)>0$. 
But if this condition holds, seeing from the definition that 
$$
\mu(a,p,d^2)\geq \mu(\Oc,a,p,d^2),
$$
we have $\mu(a,p,d^2)>0$ also.
\par
The last statement is a rephrasing of this condition using
Lemma~\ref{lm-mu0}.
\end{proof}

We thus need to find a condition on $a$ for the existence of a
solution to the system
\begin{align}
&x=\pm 1\mods{d}\label{eq-x1}\\
&x^2-ax+p=0\mods{d^2}\label{eq-quad}.
\end{align}

By the chinese remainder theorem, this admits a solution if and only
if it does locally at every prime $\ell$. So we find equivalent
conditions for $d=\ell^{\nu}$. First we consider $\ell$ odd.

\begin{lemma}\label{lm-local-odd}
Let $\ell$ be an odd prime, $d=\ell^{\nu}$. The system above admits a
solution if and only $d^2=\ell^{2\nu}\mid a^2-4p$.

\end{lemma}

\begin{proof}
Let $\Delta=a^2-4p$ denote the discriminant of the quadratic
equation~(\ref{eq-quad}). Completing the square to rewrite it as
\begin{equation}
  \label{eq-complete-sq}
\Bigl(x-\frac{a}{2}\Bigr)^2-\Delta=0\mods{d^2}
\end{equation}
(since $\ell$ is odd) shows that there is a solution
to~(\ref{eq-quad}) if and only if $\Delta$ is a square modulo $d^2$.
\par
First assume that $d^2\mid \Delta$. Then reducing
modulo $d$ and using $p\equiv 1\mods{d}$ we see that
\begin{equation}
  \label{eq-disc-d}
a^2=4\mods{\ell^{\nu}}
\end{equation}
which implies $a=\pm 2\mods{\ell^{\nu}}$ since $\ell$ is
odd. By~(\ref{eq-complete-sq}), $x=a/2$ is thus a root
of~(\ref{eq-quad}) satisfying $x\equiv\pm 1\mods{\ell^{\nu}}$. 
\par
Conversely, assume that the system has a solution
$x$. Reducing~(\ref{eq-quad}) modulo $d$ leads to $2-ax\equiv
0\mods{d}$, i.e. $x=a/2\mods{d}$ (since $x=\pm 1\mods{d}$). Let
$x=a/2+dy$. Using~(\ref{eq-complete-sq}), we have
$$
\Delta=\Bigl(x-\frac{a}{2}\Bigr)^2=d^2y^2=0\mods{d^2}.
$$
\end{proof}


We now do the same with $\ell=2$.

\begin{lemma}\label{lm-local-deux}
Let $\ell=2$, $d=\ell^{\nu}$. The
system~(\ref{eq-quad}),~(\ref{eq-x1}) above admits a 
solution if and only if $a=2b$ is even and for some $\eps=\pm 1$ we
have 
$$
\begin{cases}
b^2-p=2^{2\nu-2}y\mods{2^{2\nu}}
\text{ and }b-\eps\equiv 2^{\nu-1}y\mods{d},
\text{ with }y\equiv 0,1\mods{4}
&\text{ if $\nu\geq 2$}\\
p+1=2b\mods{4}&\text{ if $\nu=1$}.
\end{cases}
$$
\end{lemma}

\begin{proof}
This is similar to the previous one, although more tedious, and
we leave it as an exercise, as it will not be used in the sequel.
\end{proof}

\begin{remark}
If we write $x=1+dy$ in~(\ref{eq-quad}), we obtain the corresponding
equation for $y$
$$
p+1-a+dy(2-a)=0,
$$
so if $d^2\mid p+1-a$, $d\mid a-2$, any $y$ (in particular $x=1$) is a
solution (compare Proposition~\ref{pr-dun-iso}). However, for
composite $d$ other cases
are possible. In other words, the $a$ of Theorem~\ref{th-criterion} is
not necessarily the same as the $a$ of Proposition~\ref{pr-dun-iso}:
for instance take $p=241$, $d=15$. Here $a=8$ satisfies
$d^2\mid a^2-4p$, but $p+1-a=234\not\equiv 0\mods{225}$. On the other
hand, $a=17$ satisfies $a\equiv 2\mods{d}$ and $p+1-a\equiv
0\mods{d^2}$.
\end{remark}

Theorem~\ref{th-criterion} is a consequence of
Corollary~\ref{cor-critere} and Lemma~\ref{lm-local-odd}, and also 
Proposition~\ref{pr-dun-ind}.
One could incorporate
Lemma~\ref{lm-local-deux} to the statement, instead of rephrasing the
system of equations~(\ref{eq-x1}), (\ref{eq-quad}) at $2$, but it
would be more complicated. 
\par
For $d$ odd, one can further rederive,
using~(\ref{eq-te-final}), Theorem 4.9 (i) 
of~\cite{schoof-2}, namely:
\begin{proposition}\label{pr-schoof}
Let $p$ be a prime number, $d\mid p-1$ an odd integer. The number of
isomorphism classes of elliptic curves $E/\Fp_p$ with $\dun(E)\geq d$
is equal to 
$$
\sum_{\stacksum{|a|<2\sqrt{p}}{a\equiv p+1\mods{d^2}}}
{H((a^2-4p)/d^2)}.
$$
\end{proposition}

\begin{remark}
One can also tackle the question of finding points on $Y(d)$ over
finite fields by using the Riemann Hypothesis for the curve $X(d)$,
namely the inequality
$$
|N_n-(p^n+1)|\leq 2g(d)p^{n/2}
$$
for $n\geq 1$,
where $N_n=X(d)(\Fp_{p^n})$ and $g(d)$ is the genus of $X(d)$. This
implies
$$
|X(d)(\Fp_{p^n})|\geq p^n+1-2g(d)p^{n/2},
$$
and if $p^n$ is large enough compared to $d$ so that this lower bound
exceeds the number of cusps, it follows that
$Y(d)(\Fp_{p^n})\not=\emptyset$. 
\par
This approach is developed, in greater generality, by
Howe~\cite{howe}. For our purpose, we are very interested in values of
$d$ large compared to $p$ (and in the base field, $n=1$). The
inequality above is then not precise enough.
\par
Indeed we have
$$
g(d)=1+d\varphi^+(d)\psi(d)\frac{d-6}{12d}
$$
(see e.g.~\cite[(1.6.4)]{shimura-1}), of size about $d^3$, while
(Lemma~\ref{lm-count-cusps}) the 
number of cusps is $\varphi^+(d)\psi(d)$,
of size about $d^2$, so the condition to ensure
$Y(d)(\Fp_p)\not=\emptyset$, namely
$$
p+1-2g(d)\sqrt{p}>\varphi^+(d)\psi(d)
$$
is true roughly speaking for $p$ of size at least $d^6$. This is
weaker than Lemma~\ref{lm-always} below gives from
Remark~\ref{rm-critere} or Theorem~\ref{th-criterion}. 
\end{remark}

\subsection{Applications}

The previous sections give some rather simple criteria for the
existence of an elliptic curve over a finite field with a given value
of $\dun(E)$. We will deduce here some results about the possible
values of $\dun(E)$ for all elliptic curves defined over a given
finite field. Let
$$
\duns(p)=\{d \geq 1\,\mid\, d=\dun(E)\text{ for some } E/\Fp_p\}.
$$
What can be said about $\duns(p)$?
\par
We list some properties previously established:
\begin{itemize}
\item $\duns(p)$ is a subset of the set of divisors of $p-1$, indeed
  (\ref{eq-a-priori}) a subset of the set
$$
\{d\mid p-1\,\mid\,  d\leq \sqrt{p}+1\}.
$$
\item $\duns(p)$ contains $1$ and $2$.
\item $\duns(p)$ is \inductive (i.e. if $d\in \duns(p)$ and $e\mid d$,
  we have $e\in \duns(p)$, by Proposition~\ref{pr-dun-ind}).
\end{itemize}

We now consider $\duns(p)$ on average over primes $p$, and will
describe, in a certain sense, which divisors of $p-1$ belong to
$\duns(p)$. It is of particular interest to consider
primes $p$ such that $p-1$ has some divisor $d>p^{1/4}$, and see 
which of those $d$ are in $\duns(p)$.
\par
First we count on average the divisors of $p-1$ which are of a certain
size. Let
\begin{equation}\label{eq-dalpha}
\dalpha{n}=|\{
d\mid n\,\mid\, d<n^{\alpha}\}|
\end{equation}
for $n\geq 1$ and $\alpha>0$.
\par
We recall the Bombieri-Vinogradov theorem, already mentioned before.

\begin{theorem}\label{th-bom-vin}
For any $A>0$ there exists $B>0$ such that
$$
\sum_{d\leq \sqrt{X}/(\log X)^B}
{\max_{(a,d)=1} \bigl|\pi(X;d,a)-\frac{\li(X)}{\varphi(d)}
\bigr|}\ll_A \frac{X}{(\log X)^A},
$$
the implied constant depending only on $A$.
\end{theorem}

For a proof, see e.g.~\cite[\S 7]{bom-ls}.

\begin{lemma}\label{lm-count-divs}
Let $\alpha>0$ be a real number. We have
$$
\sum_{p\leq X}{\dalpha{p-1}}=f(\alpha)cX+
\oun_{\alpha}\Bigl(\frac{X}{\log X}\Bigr)
$$
where
$$
f(\alpha)=
\begin{cases}
\alpha  &
\text{if $0<\alpha\leq 1/2$}\\
(1-\alpha) &
\text{if $1/2\leq \alpha\leq 1$}\\
 1 &
\text{if $\alpha\geq 1$}
\end{cases}
$$
and
$$
c=\frac{\zeta(2)\zeta(3)}{\zeta(6)}.
$$
The implied constant depends on $\alpha$ only.
In particular,
$$
\sum_{p\leq X}{\sum_{\stacksum{d\mid p-1}{d\leq \sqrt{X}+1}}{1}}
\sim \frac{c}{2}X
$$
as $X\ra +\infty$.
\end{lemma}

\begin{proof}
This is a (simpler) variant of the proof of~(\ref{eq-tit-asy})
using the Bombieri-Vinogradov theorem and the Brun-Titchmarsh
inequality.
Indeed, if $\alpha=\frac{1}{2}$, this is
a stronger form of~(\ref{eq-tit-asy}) with explicit error term (see
for instance~\cite{fouvry}; the proof in~\cite[3.5]{halb-rich} gives
a slightly worse error term $X(\log \log X)/(\log X)$).
\par
If $\alpha<\frac{1}{2}$, we let $\beta=1/\alpha>2$ and write
$$
\sum_{p\leq X}{\dalpha{p-1}}=\sum_{d<(X-1)^{\alpha}}
{\bigl(\pi(X;d,1)-\pi(d^{\beta}+1;d,1)\bigr)}.
$$
Since $\alpha<1$, the Brun-Titchmarsh inequality yields
$$
\sum_{d<(X-1)^{\alpha}}
{\pi(d^{\beta}+1;d,1)}\ll
\sum_{d<(X-1)^{\alpha}}{\frac{d^{\beta}}{\varphi(d)(\log (d^{\beta}+1))}}
\ll \frac{X}{\log X}.
$$
Moreover, by the Bombieri-Vinogradov Theorem we have
$$
\sum_{d<(X-1)^{\alpha}}
{\bigl|\pi(X;d,1)-\frac{\li(X)}{\varphi(d)}\bigr|
}\ll \frac{X}{(\log X)^A}
$$
for any $A>0$. Since
$$
\sum_{d<(X-1)^{\alpha}}
{\frac{\li(X)}{\varphi(d)}}\sim \alpha cX\text{ as $X\ra +\infty$,}
$$
this proves the first part for $\alpha\leq \frac{1}{2}$.
\par
If $\frac{1}{2}<\alpha<1$, we use Dirichlet's trick to switch divisors
$$
\dalpha{n}=\ddalpha{1-\alpha}{n}
$$
to reduce to $1-\alpha$. Finally, for $\alpha\geq 1$,
$\dalpha{n}=d(n)$, and this is Linnik's theorem~(\ref{eq-tit-asy})
again, with error term.
\par
The last statement follows from the case $\alpha=\frac{1}{2}$, noting
that
$$
\sum_{p\leq X}{|\{d\mid p-1\,\mid\, \sqrt{X}\leq d\leq \sqrt{X}+1\}|}
=O(\sqrt{X}).
$$
\end{proof}

\begin{lemma}\label{lm-always}
Let $d\geq 1$ be an integer and $p\equiv 1\mods{d}$ a prime number. If 
$$
d< 2p^{1/4}
$$
we have $d\in \duns(p)$.
\end{lemma}

\begin{proof}
This follows from the criterion of Remark~\ref{rm-critere}, for
instance. The assumption means
that $4\sqrt{p}> d^2$, hence all $a\in \Zz/d^2\Zz$ have a lift to $\Zz$
with $|a|<2\sqrt{p}$. In particular, there is an $a$, $|a|<2\sqrt{p}$,
with $a\equiv p+1\mods{d^2}$. Since $p\equiv 1\mods{d}$, we have $a\equiv
2\mods{d}$, and by Remark~\ref{rm-critere}, $d\in\duns(p)$.
\par
For odd $d$, one can also appeal to Theorem~\ref{th-criterion} in the
same way: $a^2\mods{d^2}$ runs over all squares modulo $d^2$, and $4p$
is a square modulo $d^2$ (since $p\equiv 1\mods{d}$; indeed, if
$p=1+md$, $4p\equiv (2+md)^2\mods{d^2}$). So there exists $a$ with
$4p=a^2\mods{d^2}$, i.e. $d^2\mid a^2-4p$.
\end{proof}

\begin{remark}
One can see from the proof that this lemma is essentially best
possible, in the sense (for instance) that for any $\theta>1/4$, there
exist $p$ and $d$ with $d<p^{\theta}$ and $d\not\in\duns(p)$. This
confirms again that the condition that $\dun(E)$ be of size larger
than $p^{1/4}$ reflects a critical threshold in this subject.
\end{remark}

\begin{proposition}\label{pr-all}
We have
$$
\sum_{p\leq X}{|\duns(p)|}= \frac{cX}{4}+\oun\Bigl(\frac{X}{\log X}\Bigr)
$$
for $X\geq 2$, with an absolute implied constant.
\end{proposition}

Actually, we will prove a more precise result. As suggested by
Lemma~\ref{lm-always}, we partition $\duns(p)$ in two subsets according
to whether $d<2p^{1/4}$ or $d>2p^{1/4}$ (there can not be
equality); call those subsets $D_s(p)$ and $D_{\ell}(p)$,
respectively. 
\par
We then have:
\begin{proposition}\label{pr-small}
We have
$$
\sum_{p\leq X}{|D_s(p)|}= \frac{cX}{4}+\oun\Bigl(\frac{X}{\log X}\Bigr)
$$
for $X\geq 2$.
\end{proposition}

\begin{proposition}\label{pr-large}
We have
$$
\sum_{p\leq X}{|D_{\ell}(p)|}\ll \frac{X}{\log X}
$$
for $X\geq 2$.
\end{proposition}

Proposition~\ref{pr-all} follows immediately.

\begin{proof}[Proof of Proposition~\ref{pr-small}]
By Lemma~\ref{lm-always}, we have
$$
D_s(p)=\{d\mid p-1\,\mid\,d<2p^{1/4}\}
$$
so $|D_s(p)|=d_{1/4}(p-1)+\delta(p)$, where
$$
\delta(p)=|\{d\mid p-1\,\mid\, (p-1)^{1/4}\leq d<2p^{1/4}\}|.
$$
By Lemma~\ref{lm-count-divs}, it suffices to show that
$$
\sum_{p\leq X}
{\sum_{\stacksum{d\mid p-1}{(p-1)^{1/4}\leq d<2p^{1/4}}}
{1}}\ll \frac{X}{\log X}.
$$
This follows as before from the Brun-Titchmarsh inequality, writing
$$
\sum_{p\leq X}
{\sum_{\stacksum{d\mid p-1}{(p-1)^{1/4}\leq d<2p^{1/4}}}
{1}}
=\sum_{d\leq 2X^{1/4}}
{(\pi(d^4+1;d,1)-\pi(d^4/16;d,1))}.
$$
Equivalently, one may simply adapt the proof of
Lemma~\ref{lm-count-divs} for $\alpha=1/4$.
\end{proof}

\begin{proof}[Proof of Proposition~\ref{pr-large}]
By Remark~\ref{rm-critere}, we have $d\in D_1(p)$ if and only if
there exists $a\in \Zz$ with $|a|<2\sqrt{p}$ such that
$$
\begin{cases}
a\equiv 2\mods{d}\\
a\equiv p+1\mods{d^2}
\end{cases}
$$
Notice that $p\equiv 1\mods{d}$ is equivalent with $a\equiv 2\mods{d}$
if the last congruence holds.
\par
Now we remark that if $d\in D_{\ell}(p)$, then such an $a$ is unique:
indeed, if $a_1$ and $a_2$ satisfy the above conditions, we have
$a_1\equiv a_2\mods{d^2}$. Since $d^2>4\sqrt{p}$ and
$|a_i|<2\sqrt{p}$, this is possible only if $a_1=a_2$.
\par
Therefore we can write
$$
\sum_{p\leq X}{|D_{\ell}(p)|}=
\sum_{p\leq X}{\sum_{|a|<2\sqrt{p}}
{\sum_{\stacksum{d\mid p-1}
{\stacksum{d^2\mid p+1-a}{d^2>4\sqrt{p}}}}
{1}}}.
$$
We exchange the order of summation, getting
$$
\sum_{p\leq X}{|D_{\ell}(p)|}=
\sum_{d\leq \sqrt{X}+1}{
\sum_{\stacksum{|a|<2\sqrt{X}}{a\equiv 2\mods{d}}}{
\sum_p{1}}}
$$
where the inner sum is over primes $p$ satisfying the size conditions:
$$
\begin{cases}
p\leq X\\
p<d^4/16\\
a^2/4<p
\end{cases}
$$
and the congruence
$$
p\equiv a-1\mods{d^2},
$$
in other words
$$
\sum_{p\leq X}{|D_{\ell}(p)|}=
\sum_{d\leq \sqrt{X}+1}{
\sum_{\stacksum{|a|<2\sqrt{X}}{\stacksum{a\equiv 2\mods{d}}{2|a|<d^2}}}
{
\Bigl(
\pi\bigl(\inf({\textstyle{\frac{d^4}{16}}},X);d^2,a-1\bigr)
-\pi(a^2/4;d^2,a-1)\Bigr)
}}.
$$
We drop the second term by positivity, and write
\begin{align*}
\sum_{p\leq X}{|D_{\ell}(p)|}&\leq
\sum_{d<2X^{1/4}}{\sum_{\stacksum{2|a|<d^2}{a\equiv 2\mods{d}}}
{\pi(d^4/16;d^2,a-1)}}\\
&+\sum_{2X^{1/4}\leq d\leq \sqrt{X}+1}
{\sum_{\stacksum{|a|<2\sqrt{X}}{a\equiv 2\mods{d}}}{
\pi(X;d^2,a-1)
}}.
\end{align*}
By the Brun-Titchmarsh inequality~(\ref{eq-bt-classic}), the first
term is
\begin{align*}
\sum_{d<2X^{1/4}}{\sum_{\stacksum{2|a|<d^2}{a\equiv 2\mods{d}}}
{\pi(d^4/16;d^2,a-1)}}
&\ll \sum_{d<2X^{1/4}}{\sum_{\stacksum{2|a|<d^2}{a\equiv 2\mods{d}}}
{\frac{d^4}{\varphi(d^2)\log d}}}\\
&\ll \frac{X}{\log X}.
\end{align*}
For the second term, we further split the range of $d$ into
$2X^{1/4}\leq d\leq X^{1/2-\delta}$ and $X^{1/2-\delta}<d\leq
\sqrt{X}+1$, where $0<\delta<1/2$. For the second range, where $d$ is
very large, we simply
overcount all integers $n\equiv a-1\mods{d^2}$ instead of primes,
getting
\begin{align*}
\sum_{X^{1/2-\delta}< d\leq \sqrt{X}+1}
{\sum_{\stacksum{|a|<2\sqrt{X}}{a\equiv 2\mods{d}}}{
\pi(X;d^2,a-1)}}&\ll
\sum_{X^{1/2-\delta}< d\leq \sqrt{X}+1}
{\frac{\sqrt{X}}{d} \times \frac{X}{d^2}}\\
&\ll X^{1/2+3\delta},
\end{align*}
so if $\delta<1/6$, this saves a power of $X$ instead of merely $\log
X$.
\par
Finally, we have again by~(\ref{eq-bt-classic})
\begin{align*}
\sum_{2X^{1/4}\leq d\leq X^{1/2-\delta}}
{\sum_{\stacksum{|a|<2\sqrt{X}}{a\equiv 2\mods{d}}}{
\pi(X;d^2,a-1)}}&\ll
\sum_{2X^{1/4}\leq d\leq X^{1/2-\delta}}
{\sum_{\stacksum{|a|<2\sqrt{X}}{a\equiv 2\mods{d}}}{
\frac{X}{\varphi(d^2)\log X}}}\\
&\ll
\frac{X^{3/2}}{\log X}\sum_{2X^{1/4}\leq d\leq X^{1/2-\delta}}
{\frac{1}{d^2\varphi(d)}}\\
&\ll \frac{X}{\log X}.
\end{align*}
\end{proof}

\begin{remark}
Here the criterion given by the trace formula could also have been
used, but it would be slightly more complicated, mainly because of the
possible multiplicity of $a$ occurring for the same $d$.
\end{remark}

As a variant, we mention, and leave as an exercise, what happens for
elements of $D_1(p)$ larger than $p^{1/4+\theta}$ for some fixed
$\theta>0$. 
\begin{proposition}
Let $\theta>0$ be a real number. We have
$$
\sum_{p\leq X}{|\{d\in D_1(p)\,\mid\, d>p^{1/4+\theta}\}|}
\ll_{\theta} X^{1-2\theta}
$$
for $X\geq 2$, the implied constant depending only on $\theta$.
\end{proposition}

We also leave as an exercise the following estimate on the
average number of isomorphism classes of $E/\Fp_p$ with $\dun(E)\geq
2p^{1/4}$ (use Proposition~\ref{pr-schoof} and the trivial estimate
$H(\Delta)\ll \Delta^{1/2}\log \Delta$, see
e.g.~\cite[Th. 7-24]{cox}).

\begin{proposition}\label{pr-large-poids}
We have
$$
\sum_{p\leq X}
{\sum_{\stacksum{d\mid p-1}{\stacksum{2\nmid d}{d\geq 2p^{1/4}}}} 
{|\{E/\Fp_p\,\mid\, \dun(E)\geq d\}|}}
\ll X^{5/4}.
$$
\end{proposition}

For comparison, the total number of isomorphism classes of $E/\Fp_p$
with $p\leq X$ is $\sim X^2$ (there are $p$ possible $j$-invariants
and, except for cubic and biquadratic twists for $j=0$, $1728$, two
isomorphism classes for each $j$-invariant, see
e.g.~\cite[X-5]{silv}). 

\begin{remark}
For heuristic purposes in trying to make guesses about the
distribution of \outside primes for elliptic curves, it is really a
lower-bound for $|D_{\ell}(p)|$ that one would like to have on
average, or more precisely for the quantity in
Proposition~\ref{pr-large-poids}. This looks like a fairly hard
problem: one can see in the proof of Proposition~\ref{pr-large} that
it boils down to assertions about the equidistribution of primes $\leq
Y$ to moduli which are $\gg Y^{1/2}$, and moreover with ``initial
term'' $a-1$ which vary. The latter constraint, in particular, seems
currently incompatible with the methods developed by Bombieri,
Friedlander and Iwaniec~\cite{bom-fried-iwa}.
\end{remark}

\mysection{Numerical examples}
\label{sec-numerics}

The various problems we have considered lend themselves easily to
numerical experimentation using computer packages for elliptic curves
computations. We have used the PARI/GP system and written scripts to
perform the following computations, for an elliptic curve $E/\Qq$
given by a Weierstrass equation:
\begin{itemize}
\item{
Compute the invariants $\dun(p)$, $\dde(p)$ at a prime $p$, and the
sum $S_E(X;\dun)$. Also, find the weak \outside primes of $E$ which
are $\leq X$, and if the order of the Galois groups $G_d$ can be
computed, the \outside primes $\leq X$.
}
\item{
Compute the multiplicity functions $\grandm(n)$ or $\petitm(p)$, the
number of $E$-twins $\leq X$ and more generally
the various moments $\somg_k(X)$, $\somp_k(X)$.
}
\end{itemize}

The numerical results can be compared to the predictions, when we have
some. Especially if $E$ is a Serre curve (Section~\ref{sec-comp-ce}),
one can compare $S_E(X;\dun)$ with the conjectural asymptotic
$$
S_E(X;\dun)\sim c(E)\li(X).
$$
\par
The PARI system does not implement (yet) the computation of $\dun(p)$
as a primitive function although, based on Cohen's description of the
Shanks-Mestre algorithm to compute $a_p$ (\cite[7.4.3]{cohen-1}),
this should be almost as fast as computing $a_p$.
However one can write a simple enough algorithm by computing 
the exponent (i.e. $\dun\dde$) of $E_p(\Fp_p)$ by looking for
an element of maximal order, either by ``exhaustion'' or more
efficiently (as suggested by K. Belabas) by picking up a few
``random'' points on $E_p(\Fp_p)$ and 
taking the l.c.m of their orders.\footnote{~In the
computations below, this was done with $20$ random points, so
in theory the results might be off by a small amount. However, it is
easy to repeat the computations for the primes yielding ``large''
values of $\dun(p)$, thus ensuring their correctness.
} Moreover, for primes $p$ with
$|E_p(\Fp_p)|$ squarefree, one has $\dun(p)=1$ without further
computations, and this happens quite often if the curve has no
non-trivial rational $2$-torsion points.
\par
Computing elliptic twins is even simpler, and the computation of the
sums 
$$
\somg_k(X)=\sum_{n\leq X}{\grandm(n)^k}
$$
can be performed using very little memory by operating by blocks of
$n$. Numerically, $\grandm(n)$ is always very small so $\somg_k(X)$ is
very close to $\somp_{k-1}(X)$ (compare~(\ref{eq-g-p-equiv})). Also we
computed the modified first moment
$$
\somg'(X)=\sum_{\stacksum{n\leq X}{n\text{ twin value}}}
{\grandm(n)}.
$$
Note that we have obviously
$$
\somg'(X)=\ejumbis(X)+\oun(\sqrt{X}).
$$
(see~(\ref{eq-ejumbis}) for $\ejumbis(X)$).

\subsection{The test curves}

We used two non-CM curves, which are Serre curves, and one CM
curve. Here are their id-sheets:

\begin{example}\label{example-1}
Consider the curve (see~\cite[5.9.2]{serre-div},~\cite[I \S
7]{lang-trotter})
$$
E\,:\, y^2=x^3+6x-2
$$
with $j(E)=2^93$, discriminant $-2^6 3^5$, conductor
$1728$. It has rank $0$.
By~\cite[Th. 7.1]{lang-trotter}, this curve is a Serre curve and
$m=3$ in this case.
\par
Using Corollary~\ref{cor-const-serre},~(\ref{eq-c0}), we have
\begin{align}
c'(E)&=\frac{5461}{5425}=1.0066\ldots\nonumber\\
c(E)&=c'(E)c_0=1.2668\ldots
\label{eq-cste-ex1}
\end{align}
\end{example}

\begin{example}\label{example-2}
Consider the curve (see~\cite[5.5.6]{serre-div})
$$
F\,:\, y^2+y=x^3-x
$$
with $j(F)=2^{12}3^3/37$, discriminant $37$, conductor
$37$. It has rank $1$, the point $(0,0)$ being of infinite order.
It is also a Serre curve and $m=37$. (It is also studied by Mazur and
Swinnerton-Dyer in~\cite{mazur-sd}).
\par
Using Corollary~\ref{cor-const-serre}, we have
\begin{align}
c'(F)&=\frac{1732338101}{1732332625}=1.000003\ldots\nonumber\\
c(F)&=c'(F)c_0=1.2584\ldots
\label{eq-cste-ex2}
\end{align}
(the value of $c(F)$ differs from $c_0$ by less than $10^{-5}$).
\end{example}

\begin{example}\label{example-3}
The last curve is the CM curve~(\ref{eq-cm}) of
Example~\ref{ex-cm-curve}, namely
$$
A\,:\,y^2=x^3-x,
$$
(with CM by $\Zz[i]$). The expected behavior is now
$$
S_A(X;\dun)\sim c(A)X
$$
with $c(A)$ given by~(\ref{eq-const-cm}).
\end{example}

\subsection{Numerical examples: the elliptic splitting problem}

We now give a few examples of computations of averages of $\dun$.
Here are some experimental data for $p\leq \pmax$, for the curves $E$
and $F$ of Examples~\ref{example-1} and~\ref{example-2}.
\par
\begin{center}
\begin{supertabular}{c|c|c|c|c|c}
$X$ & $\pi(X)$ & $S_E(X;\dun)$ & Ratio & $S_F(X;\dun)$ & Ratio \\
\hline
100,000 &    9592     & 11945   & 1.24530 & 11944    & 1.24520 \\
500,000 &    41538    & 52418   & 1.26192 & 51969    & 1.25111 \\
1,000,000 &  78498    & 99144   & 1.26301 & 98465    & 1.25436 \\
5,000,000 &  348513   & 440751  & 1.26466 & 438079   & 1.25699 \\
10,000,000 & 664579   & 841232  & 1.26581 & 835662   & 1.25743 \\
15,000,000 & 970704   & 1229075 & 1.26616 & 1220393  & 1.25722 \\
20,000,000 & 1270607  & 1608929 & 1.26626 & 1597802  & 1.25751 \\
30,000,000 & 1857859  & 2352704 & 1.26635 & 2336778  & 1.25778 \\
40,000,000 & 2433654  & 3081940 & 1.26638 & 3061994  & 1.25818 \\
50,000,000 & 3001134  & 3800076 & 1.26621 & 3775641  & 1.25807 \\
60,000,000 & 3562115  & 4510928 & 1.26636 & 4480730  & 1.25788 \\
\end{supertabular}
\end{center}

The agreement with the expected behavior seems quite good, but it
should be noticed that only values of $d$ (in the sense
of~(\ref{eq-dun-som})) which are fairly small actually occur in this
range. In accordance with~(\ref{eq-cste-ex1})
and~(\ref{eq-cste-ex2}), the sum for $E$ tends to be slightly larger
than that for $F$.
\par
All outside primes $\leq 300,000,000$ were computed. It turns out that
there are very few of them. Here is the complete list, indicating the
prime $p$, the value of $\dun(p)$ and the order of the Galois group
$G_d$

\begin{center}
\begin{supertabular}{c|c|c}
$p$ & $\dun(E,p)$ & $|G_d|$ \\
\hline
196561 & 140 &  92897280 \\
4095037 & 162 &  76527504 \\
13403893 & 114 &  17729280 \\
30626899 & 106 &  46433088 \\
53629561 & 184 &  410370048 \\
54460963 & 258 &  480598272 \\
76391737 & 172 &  320398848 \\
132576571 & 127 &  258080256 \\
138085949 & 143 &  345945600 \\
145030393 & 312 &  966131712 \\
\end{supertabular}
\end{center}

There are $20$ additional weak \outside primes, for instance
$p=779761$ with $\dun(p)=36=p^{\alpha}$ with $\alpha=0.26\ldots$
\par
The impact of the single very large value of $\dun$ at $p=196561$ is
quite noticeable: we have
\begin{center}
\begin{supertabular}{c|c|c|c}
$X$ & $S_E(X;\dun)$ & $\pi(X)$ & Ratio \\
\hline
196560 & 22218 & 17700 & 1.2552 \\
196561 & 22358 & 17701 & 1.2630 \\
\end{supertabular}
\end{center}

\medskip
In another direction, here is a table listing, for those $d\leq 140$
for which at least one 
$p\leq 3,000,000$ splits completely in $\Qq(E[d])$, how many do:
$\pi_X(E;d,1)$ is in the second row, the third is the ratio
$\pi(X)/\pi_E(X;d,1)$, for comparison with $|G_d|$. 

\begin{center}
\begin{supertabular}{c|c|c|c|c|c|c}
$d$ & 2 & 3 & 4 & 5 & 6 & 7\\
\hline
Number & 13032 & 1624 & 783 & 164 & 502 & 28 \\
Ratio &  6.0223 & 48.335 & 100.25 & 478.63 & 156.36 & 2803.4\\
$|G_d|$ & 6 & 48 & 96 & 480 & 144 & 2016 \\
\hline\hline

$d$ & 8 & 9 & 10 & 11 & 12 & 13 \\
\hline
Number & 40 & 17 & 33 & 7 & 28 & 4 \\
Ratio &  1962.4 & 4617.4 & 2378.6 & 11213. & 2803.4 & 19624.\\
$|G_d|$ & 1536 & 3888 & 2880 & 13200 & 2304 & 26208 \\
\hline\hline

$d$ &  14 & 15 & 16 & 17 & 18 & 19 \\
\hline
Number & 6 & 2 & 1 & 1 & 8 & 1 \\
Ratio & 13082. & 39248. & 78496. & 78496. & 9812.0 & 78496 \\
$|G_d|$ &  12096 & 23040 & 24576 & 78336 & 11664 & 123120 \\
\hline\hline

$d$ &   20 & 21 & 23 & 24 & 28 & 30 \\
\hline
Number & 1 & 1 & 2 & 2 & 1 & 1 \\
Ratio & 78496. & 78496. & 39248. & 39248. & 78496. & 78496.\\
$|G_d|$ & 46080 & 96768 & 267168 & 36864 & 193536 & 69120 \\
\hline\hline

$d$ &   35 & 36 & 70 & 140 \\
\cline{1-5}
Number & 1 & 1 & 1 & 1 \\
Ratio &78496. & 78496. & 78496. & 78496. \\
$|G_d|$ & 967680 & 186624 & 5806080 & 92897280\\
\end{supertabular}
\end{center}

\medskip
As for $F$, here is the table listing the \outside primes $\leq
300,000,000$.
\begin{center}
\begin{supertabular}{c|c|c}
$p$ & $\dun(F,p)$ & $|G_d|$ \\
\hline
8317 & 11 &  13200 \\
63317 & 22 &  79200 \\
657493 & 44 &  1267200 \\
1258667 & 37 &  1822176 \\
11019023 & 98 &  29042496 \\
\end{supertabular}
\end{center}

One can see again that those $p$ for which $\dun(p)$ is large have an
important effect; here we have
\begin{center}
\begin{supertabular}{c|c|c|c}
$X$ & $S_F(X;\dun)$ & $\pi(X)$ & Ratio \\
\hline
63313 & 7849 & 6343 & 1.2374 \\
63317 & 7871 & 6344 & 1.2407\\
657491 & 66953 & 53378 & 1.2543 \\
657493 & 66997 & 53379 & 1.2551 \\
\end{supertabular}
\end{center}

Here is the table of the number of primes $p\leq
3,000,000$ which split in $\Qq(F[d])$ for $2\leq d\leq 44$ (those $d$
for which no $p$ splits are omitted):
\begin{center}
\begin{supertabular}{c|c|c|c|c|c}
$d$ & 2 & 3 & 4 & 5 & 6 \\
\hline
Number & 13034 & 1645 & 790 & 152 & 268 \\
Ratio & 6.0224 & 47.718 & 99.363 & 516.42 & 292.89 \\
$|G_d|$ & 6 & 48 & 96 & 480 & 288 \\
\hline\hline

$d$ & 7 & 8 & 9 & 10 & 11\\
\hline
Number & 30 & 56 & 15 & 22 & 10\\
Ratio &  2616.5 & 1401.7 & 5233.1 & 3568.0 & 7849.7 \\
$|G_d|$ & 2016 & 1536 & 3888 & 2880 & 13200 \\
\hline\hline

$d$ &  12 & 13 & 14 & 15 & 16 \\
\hline
Number & 16 & 2 & 4 & 2 & 4 \\
Ratio & 4906.0 & 39248. & 19624. & 39248. & 19624. \\
$|G_d|$ & 4608 & 26208 & 12096 & 23040 & 24576 \\
\hline\hline

$d$ &  21 & 22 & 24 & 44 \\
\cline{1-5}
Number & 1 & 3 & 2 & 1 \\
Ratio & 78497. & 26165. & 39248. & 78497. \\
$|G_d|$ & 96768 & 79200 & 73728 & 1267200 \\

\end{supertabular}
\end{center}

For the CM curve $A$ of Example~\ref{example-3}, we get the following
for $p\leq 30,000,000$, where we compare $S_A(X;\dun)$ with $X$ in the
last column: 
\begin{center}
\begin{supertabular}{c|c|c}
$X$ & $S_A(X;\dun)$ & Ratio \\
\hline
10000 & 5410 & 0.5410 \\
100000 & 55578 & 0.5558 \\
500000 & 267450 & 0.5349 \\
1000000 & 529742 & 0.5297 \\
5000000 & 2633630 & 0.5267 \\
10000000 & 5274876 & 0.5275 \\
15000000 & 7839124 & 0.5226 \\
20000000 & 10386178 & 0.5193 \\
25000000 & 13027268 & 0.5211 \\
30000000 & 15665348 & 0.5222 \\
\end{supertabular}
\end{center}

The expected linear growth of $S_G(X;\dun)$ seems also apparent.

\subsection{Numerical examples: elliptic twins}
\label{ssec-num-twins}

Motivated by the rough heuristic of Section~\ref{ssec-heuristic}, for
non-CM curves we compare $\somg'(X)$ with\footnote{ As usual, this  
gives a much better approximation than $X/(\log X)^2$.
}
$$
\li_2(x)=\int_2^x{\frac{dt}{(\log t)^2}}
=\li(x)-\li(2)-\frac{x}{\log x}+\frac{2}{\log 2}.
$$
\par
The first table lists some values of $X$, $\somg'(X)$ and
$\somg'(X)/\li_2(X)$ for the curves $E$ and $F$, for $X\leq 10^8$. 
\par
\medskip
\begin{center}
\begin{supertabular}{c||c|c|c|c}
$X$ & $\somg'_E(X)$ &  $\somg'_E(X)/\li_2(X)$ & $\somg'_F(X)$ &
$\somg'_F(X)/\li_2(X)$ \\
\hline
1000 & 32   & 0.9226 & 29   & 0.8361   \\
10000 & 133& 0.8198 & 154  & 0.9492  \\
100000 & 1110  & 1.1736 & 1062  & 1.1229\\
1000000  & 7364  & 1.1788 & 7349  & 1.1764\\
5000000  & 29583  & 1.2079 & 29045 & 1.1860 \\
10000000  & 54036  & 1.2143 & 52734 &  1.1850\\
20000000  & 98582  & 1.2136 & 97226 &  1.1969  \\
40000000  & 181587 & 1.2197 & 178934 &  1.2018\\
60000000  & 259489 & 1.2206 & 255478 &  1.2018  \\
80000000  & 333974 & 1.2193 & 329150 &  1.2017\\
99980000  & 407033 & 1.2205 & 401293 &  1.2033    \\
\end{supertabular}
\end{center}

\medskip
Next we list the multiplicities $\grandm(n)$ occurring for twin values
$n$: in this range, $\grandm(n)\leq 5$, and the number of integers
with a given $\grandm(n)=k>1$ is as follows:
\par
\medskip
\begin{center}
\begin{supertabular}{c||c|c|c|c}
$k$ & 2 & 3 & 4 & 5 \\
\hline
$E$ & 194197 & 5982 & 167 & 5 \\
$F$ & 191817 & 5685 & 146 & 4 \\
\end{supertabular}
\end{center}
\medskip
The values of $n\leq 10^8$ with $\grandm_E(n)=5$ are
$$
n\in \{13269240,14469576,20024896,52472068,64703760\}
$$
and those with $\grandm_F(n)=5$ are
$$
n\in \{5597128,64220836,85004608, 86998320\}.
$$
To compare with~(\ref{eq-conj-multimax}), note that
$$
\frac{\log x}{\log \log x}=
\begin{cases}
5.7980&\text{ for } x=10^7\\
6.3225&\text{ for } x=10^8.
\end{cases}
$$
Because of the very small number of $n$ with $\grandm(n)>2$,
$\ejum(X)$ (see~(\ref{eq-ejum})) is almost equal to $\frac{1}{2}
\somg'(X)$. In particular, the numerical data seems to
confirm~(\ref{eq-ejum-conj}) for $E$ and $F$.
\par
\medskip
We now consider the CM curve $A/\Qq$. Of course, the field of
definition does not contain the CM field, as assumed in
Section~\ref{sec-cm-resultats}. However, it is very simple to adapt the
arguments there to this case.
\par
For supersingular $p$, i.e. $p\equiv 3\mods{4}$, we have $n_p=p+1$; in
particular if we write 
$$
\grandm(n)=
\grandm_o(n)+\grandm_s(n),
$$
where $\grandm_o(n)$ (resp. $\grandm_s(n)$) is the number of ordinary
primes $p$ with $n_{p}=n$ (resp. supersingular primes), it
follows that $\grandm_s(n)=0$ or $1$ according to whether $n-1$ is
prime $\equiv 3\mods{4}$ or not (note that $n_p\equiv 0\mods{4}$ for
all $p$ since $A[2]\subset A(\Qq)$, so $n_p-1\equiv 3\mods{4}$ for all
$p$). 
\par
We thus get the bound
\begin{equation}\label{eq-borne-ind-cm-q}
\grandm(n)\leq 1+\frac{1}{2}r(n)
\end{equation}
instead of~(\ref{eq-borne-ind-cm}).
\par
To estimate $\somg_k(X)$, write
\begin{align*}
\somg_k(X)&=\sum_{n\leq X}{(\grandm_o(n)+\grandm_s(n))^k}\\
&=\sum_{j=0}^k{\binom{k}{j}
\sum_{n\leq X}{\grandm_s(n)^{k-j}\grandm_o(n)^{j}}}\\
&\leq \sum_{j=0}^k{\binom{k}{j}\somg_{o,j}(X)}
\end{align*}
since $\grandm_s(n)^{k-j}\leq 1$, where $\somg_{o,j}(X)$ is the $j$-th
moment of $\grandm_o(n)$. To the latter sum, we can clearly apply the
arguments used in Section~\ref{sec-cm-resultats} verbatim, and deduce
$$
\somg_{o,j}(X)\ll_j X(\log X)^{\beta(j-1)+\eps}\text{ with }
\beta(j)=2^j-j-2,\text{ for any }\eps>0
$$
hence we have:
\begin{proposition}\label{pr-cm-q}
For all $k\geq 0$ and $X\geq 2$ we have
\begin{align*}
\somg_k(X) & \ll_{\eps} X(\log X)^{\beta(k-1)+\eps}\text{ for $k\geq 1$}\\
\somp_k(X) & \ll_{\eps} X(\log X)^{\beta(k)+\eps},
\end{align*}
with $\beta(k)=2^k-k-2$ for any $\eps>0$, the implied constant
depending only on $k$ and $\eps$.
\end{proposition}
\par
\medskip
Computations were performed for $p\leq 20,000,000$.
Here is a table with values of $\ejum(X)$, $\somg'(X)$
and of the ratio $\somg'(X)/\li(X)$:
\par\medskip
\begin{center}
\begin{supertabular}{c||c|c|c}
$X$ & $\somg'_A(X)$ &  $\somg'_A(X)/\li(X)$ & $\ejum(X)$ \\
\hline
1000 & 67 & 0.37723 & 27\\
10000 & 486 & 0.39000 & 187\\
100000 & 3693 & 0.38349 & 1430\\
1000000 & 29068 & 0.36969 & 11052\\
5000000 & 126445 & 0.36268 & 47674\\
7500000 & 182930 & 0.35975 & 68842\\
10000000 & 238563 & 0.35878 & 89693\\
12500000 & 292994 & 0.35778 & 110021\\
15000000 & 346590 & 0.35692 & 130095\\
17500000 & 399567 & 0.35624 & 149871\\
20000000 & 451562 & 0.35530 & 169294\\
\end{supertabular}
\end{center}
\medskip
Here is a table with values of $\somg_2(X)$ et $\somg_3(X)$, compared
with $\li(X)$ and $X$ respectively:
\par\medskip
\begin{center}
\begin{supertabular}{c||c|c|c|c}
$X$ & $\somg_2(X)$ & $\somg_2(X)/\li(X)$ & $\somg_3(X)$ & 
$\somg_3(X)/X$ \\
\hline
100000 & 16757 & 1.7401 & 43637 & 0.43637 \\
500000 & 73154 & 1.7582 & 198966 & 0.39793 \\
1000000 & 138492 & 1.7613 & 384224 & 0.38422 \\
2500000 & 323992 & 1.7680 & 919320 & 0.36772 \\
5000000 & 618660 & 1.7745 & 1786380 & 0.35727 \\
7500000 & 902363 & 1.7746 & 2635021 & 0.35133 \\
10000000 & 1180791 & 1.7758 & 3469855 & 0.34698 \\
12500000 & 1454892 & 1.7766 & 4285228 & 0.34281 \\
15000000 & 1724899 & 1.7763 & 5098883 & 0.33992 \\
17500000 & 1992562 & 1.7765 & 5897698 & 0.33701 \\
20000000 & 2258677 & 1.7772 & 6714287 & 0.33571 \\
\end{supertabular}
\end{center}
\medskip
Here is the table of values $>1$ taken by $\grandm(n)$ in this range
(those $k$ for which no $n$ satisfies $\grandm(n)=k$ are omitted):
\begin{center}
\begin{supertabular}{c|c|c|c|c|c|c|c|c|c}
2 & 3 & 4 & 5 & 6 & 7 & 8 & 9 & 10 & 11 \\
\hline
106007 & 37191 & 14291 & 6123 & 2835 & 1360 & 670 & 386 & 195 & 108\\
\hline\hline
12 & 13 & 14 & 15 & 16 & 17 & 18 & 19 & 20 & 24\\
\hline
60 & 33 & 13 & 9 & 7 & 1 & 2 & 1 & 1 & 1 \\
\end{supertabular}
\end{center}
The $n$ with $\grandm(n)=24$ is $n=12818000$. Notice that $n=2^4\cdot
5^3 \cdot 13\cdot  17\cdot  29$, each prime $\not=2$ being (of course)
a sum of two squares. We have $r(n)=32$ in this case. In practice, it
is quite easy to find rather large multiplicities without constructing
a complete table: take an integer $n$ divisible by $4$ (because
$A[2]\subset A(\Qq)$) and with many prime factors $\equiv 1\mods{4}$
so that $r(n)$ is large, and look at the primes $p$, $n^-\leq p\leq n^+$,
for those with $n_p=n$. 
\par
For comparison, the integers $n\leq 10^8$ with $\grandm_E(n)=5$ or
$\grandm_F(n)=5$ factorize as follows:
\begin{gather*}
13269240=2^3\cdot 3^2\cdot 5\cdot 29\cdot 31\cdot 41,\quad
14469576=2^3\cdot 3\cdot 11\cdot 23\cdot 2383,\\
20024896=2^6\cdot 139\cdot 2251,\quad
52472068=2^4\cdot 11\cdot 37\cdot 167\cdot 193,\\
64703760=2^4\cdot 3\cdot 5\cdot 11\cdot 24509,\quad
5597128=2^3\cdot 699641,\\
64220836=2^2\cdot 19\cdot 491\cdot 1721,\quad
85004608=2^6\cdot 13\cdot 71\cdot 1439.\quad
86998320=2^4\cdot 3^3\cdot 5\cdot 40277,
\end{gather*}
the prime factors exhibiting no obvious property (?).

\mysection{Conclusion}
\label{sec-concl}

The many questions raised in this paper seem very hard to attack, but
on the other they seem to be very interesting from the point of view
of analytic number theory. Given the extensive experience with the
distribution of primes in arithmetic progressions to large moduli,
and the (much more modest) first results for CM curves obtained here,
one would like to have some kind of sieve method available for the
non-CM curves: roughly speaking, sieve is powerful because it
exploits the embedding of primes inside the integers, and because the
divisibility of integers by a given $d\geq 1$ can be used to recover
primes by inclusion-exclusion, so some of the regularity of the
distribution of integers can be exploited.
\par
For a non-CM curve $E/\Qq$, the function 
$\dun(p)$ has no obvious interpretation as the restriction to primes
of an arithmetic function defined for all $n$, whereas if $E/\Qq$ has
CM, $\dun(p)$ is $b(\pi-1)$, where $\pi$ is the Frobenius at $p$ and
$b(a)$ is defined for any $a\in \End(E)$ as the largest integer 
$b\in\Zz$ with $(b)\mid (a)$.\footnote{
The results of Duke and Toth
(\cite{duke-toth}) can be used to ``lift'' the Frobenius on $E_p$ to a
matrix in $M(2,\Zz)$, well-defined up to $GL(2)$-conjugacy, which
reduced modulo $d$ gives the action of $\frob_p$ on $\Qq(E[d])$ for
any $d$ (prime to the discriminant). But I do not see how to
isolate the conjugacy classes of this type; the set of all matrices is
too big to give information on a single elliptic curve.}
\par
Also, despite the fact that the modularity of elliptic curves would
seem to provide a ``dual view'', similar to that of Dirichlet
characters instead of $1$-dimensional Galois representations, it is
really the Artin $L$-functions attached to the fields $K(E[d])/K$ which
are of importance. Those can have rank as large as $d$ (roughly),
which makes all current analytic techniques incapable of dealing with
them, individually or on average, even assuming the Artin conjecture,
or that they are automorphic $L$-functions.
\par
Thus it seems much work is required to understand those analytic
problems. As for arithmetic progressions however, where the stumbling
block of the Riemann Hypothesis has often been circumvented by
startling new results (Linnik's dispersion method, the
Bombieri-Vinogradov theorem, the results of
Bombieri-Friedlander-Iwaniec, etc...), one may hope that there is much
to discover.

\end{document}